\documentclass{amsart}
\usepackage{amssymb,euscript,amsmath, mathrsfs, amscd, mathabx}
\usepackage[pdftex]{graphicx}
\usepackage[pdftex]{color}

\newcounter{ENUM}
\newcommand{\itm}{\item}
\newenvironment{ilist}[1][0]{\renewcommand{\theENUM}{\roman{ENUM}}\renewcommand{\itm}{\addtocounter{ENUM}{1}\item[(\theENUM)]}\begin{itemize}\setcounter{ENUM}{#1}}{\end{itemize}}
\newenvironment{Ilist}[1][0]{\renewcommand{\theENUM}{\Roman{ENUM}}\renewcommand{\itm}{\addtocounter{ENUM}{1}\item[(\theENUM)]}\begin{itemize}\setcounter{ENUM}{#1}}{\end{itemize}}
\newenvironment{alist}[1][0]{\renewcommand{\theENUM}{\alph{ENUM}}\renewcommand{\itm}{\addtocounter{ENUM}{1}\item[(\theENUM)]}\begin{itemize}\setcounter{ENUM}{#1}}{\end{itemize}}

\newcommand{\margh}[1]{}

\def\risom{\overset{\sim}{\rightarrow}}

\input xy
\xyoption{all}
\CompileMatrices

\newcommand{\el}{$\ell$}
\newcommand{\e}{\varepsilon}
\def\ZZ{{\mathbb Z}}

\def\AA{{\mathbb A}}

\def\QQ{{\mathbb Q}}

\def\cF{{\mathcal F}}
\def\cG{{\mathcal G}}

\def\cS{{\mathcal S}}

\def\cP{{\mathcal P}}

\def\cU{{\mathcal U}}

\def\cM{{\mathcal M}}

\def\cLG{\mathcal{LG}}

\def\cPic{{\mathcal P}ic}
\def\sE{{\mathscr E}}
\def\sF{{\mathscr F}}

\def\sK{{\mathscr K}}
\def\sL{{\mathscr L}}

\def\sO{{\mathscr O}}
\def\sQ{{\mathscr Q}}
\def\sV{{\mathscr V}}
\def\sW{{\mathscr W}}

\def\fg{{\mathfrak g}}

\def\fm{{\mathfrak m}}

\def\vp{\varphi}

\def\Sch{\operatorname{Sch}}
\def\BSch{B\text{-}\Sch}

\def\Hom{\operatorname{Hom}}

\def\GL{\operatorname{GL}}

\def\Spec{\operatorname{Spec}}
\def\Pic{\operatorname{Pic}}

\def\id{\operatorname{id}}

\def\rk{\operatorname{rk}}

\def\ord{\operatorname{ord}}
\def\refn{\operatorname{ref}}

\def\EHT{\operatorname{EHT}}

\def\spn{\operatorname{span}}

\def\LG{\operatorname{LG}}

\def\I{\operatorname{I}}
\def\II{\operatorname{II}}
\def\cl{\operatorname{cl}}

\newtheorem{thm}{Theorem}[subsection]
\newtheorem{prop}[thm]{Proposition}
\newtheorem{lem}[thm]{Lemma}
\newtheorem{cor}[thm]{Corollary}

\theoremstyle{definition}
\newtheorem{defn}[thm]{Definition}

\newtheorem{ex}[thm]{Example}
\newtheorem{sit}[thm]{Situation}
\theoremstyle{remark}
\newtheorem{notn}[thm]{Notation}
\newtheorem{rem}[thm]{Remark}

\numberwithin{equation}{subsection}

\begin{document}
\title{Limit linear series moduli stacks in higher rank}
\author{Brian Osserman}
\begin{abstract} 
In order to prove new existence results in Brill-Noether theory for 
rank-$2$ vector bundles with fixed special determinant, we develop 
foundational definitions and results for limit linear series of higher-rank 
vector bundles. These include two entirely new constructions of ``linked
linear series'' generalizing earlier work of the author for the classical
rank-$1$ case, as well as a new canonical stack structure for the previously
developed theory due to Eisenbud, Harris and Teixidor i Bigas. This last
structure is new even in the classical rank-$1$ case, and yields the first
proper moduli space of Eisenbud-Harris limit linear series for families of 
curves. We also develop results comparing these three constructions.
\end{abstract}

\thanks{The author was partially supported by NSA grant H98230-11-1-0159
during the preparation of this work.}
\maketitle

\tableofcontents

\section{Introduction}

Ever since the main questions of classical Brill-Noether theory were
resolved in the 1980's, the natural generalization to higher-rank vector
bundles has been a subject of attention. At its most basic, this is
the study of how many sections a (semistable) vector bundle of given rank
and degree can have on a general curve. More precisely, one
studies $\fg^k_{r,d}$s, consisting of pairs $(\sE,V)$ of a vector bundle 
$\sE$ of rank $r$ and degree $d$, together with a $k$-dimensional space $V$
of global sections of $\sE$. This generalization is of
fundamental interest due to the basic nature of understanding questions
on vector bundles with sections, and on morphisms from curves to
Grassmannians. However, despite our rather rudimentary understanding thus
far, the subject has found several important applications, including 
work of Mukai \cite{mu4} classifying curves of low genus,
which he in turn applied to the classification of Fano threefolds 
\cite{mu6}, the work of Farkas and Popa \cite{f-p1} giving a 
counterexample to the slope conjecture, and most recently,
work of Bhosle, Brambila-Paz, and Newstead proving the rank-$1$ case
of a conjecture of Butler regarding stability of the kernel of the 
evaluation map of linear series \cite{b-b-n2}, which has in turn
been applied by Brambila-Paz and Torres-Lopez \cite{b-t1} to obtain new 
results on Chow stability of curves.

Despite extensive work and many partial results (see \cite{g-t1} for a
survey), there is as yet not even a comprehensive conjecture for the
case of rank-$2$ vector bundles. The naive generalization of the 
Brill-Noether theorem to higher rank fails in every possible way, with
complications arising in particular from stability conditions and from
the role of vector bundles with special determinant. The latter was
first observed by Bertram and Feinberg \cite{b-f2} and Mukai \cite{mu2}
for the case of bundles of rank $2$ with canonical determinant,
and studied more systematically in \cite{os16} and \cite{os19}. We now
explain it in more detail.

The natural generalization of the classical Brill-Noether number to 
the higher-rank case is given by 
$$\rho:=1+r^2(g-1)-k(k-d+r(g-1)).$$
Because we will consider moduli stacks rather than coarse moduli spaces,
the actual expected dimension will be $\rho-1$ for the varying determinant
case. Thus, the moduli stack $\cG^k_{r,d}(X)$ of $\fg^k_{r,d}$s on a general 
curve $X$ of genus $g$ has every component of dimension at least $\rho-1$. 
In the fixed determinant case, the expected dimension remains $\rho-g$,
because fixing the determinant rigidifies the moduli problem. 
In many cases, $\cG^k_{r,d}(X)$ is known to have components of
dimension $\rho-1$ (or $\rho-g$ in
the fixed determinant case). However, in the case of rank $2$ and fixed 
determinant $\sL$, if $h^1(\sL)>0$, then the dimension of the moduli space
is always at least $\rho-g+\binom{k}{2}$. Bertram, Feinberg and Mukai
conjectured that the locus of $\fg^k_{2,d}$s with canonical determinant
behaves like classical Brill-Noether loci, and the portion of their 
conjecture asserting non-emptiness when $\rho-g+\binom{k}{2} \geq 0$ remains 
open.
 
The primary objective of the present paper, together with \cite{o-t2}, is
to develop the necessary machinery for the use degeneration techniques to 
prove existence results for moduli spaces of $\fg^k_{2,d}$s in the case 
of fixed special determinant. The most powerful tool for studying 
classical Brill-Noether theory is the Eisenbud-Harris theory of limit
linear series \cite{e-h1}. In \cite{te1}, Teixidor i Bigas generalized 
this theory to 
the higher-rank case. The most technical part of the theory is the proof 
of smoothing theorems, in which one shows that the existence of families 
of limit linear series having the expected dimension on a given reducible
curves implies the existence of linear series (or $\fg^k_{r,d}$s) on 
smooth curves. These results are proved by constructing suitable moduli
spaces for families of curves, and proving dimensional lower bounds. 
There is no difficulty in carrying out the Eisenbud-Harris construction
in the higher-rank case using the naive expected dimensions $\rho-1$ and
$\rho-g$. However, this result is useless in the case of special 
determinant, because the dimension will always be strictly larger than 
$\rho-g$ whenever $k \geq 2$. In order to obtain a usable theory in this
setting, is thus necessary to combine the dimensional lower bounds of
the special determinant setting with those of the theory of limit linear
series. Unfortunately, in the context of Eisenbud-Harris-Teixidor limit
linear series, it is not clear how to track the symmetries which cause the
higher expected dimension in the special determinant case. 

In the present paper, we introduce new moduli stacks of higher-rank limit 
linear series (which, for the sake of avoiding confusion, we label ``linked
linear series''). These both generalize the construction of \cite{os8} to
higher rank and to curves with more than two components. These stacks have 
comparison morphisms between them and to the stack of 
Eisenbud-Harris-Teixidor limit linear series. As in the case of \cite{os8},
these morphisms are isomorphisms on certain open loci of interest, but on 
the boundary frequently have positive-dimensional fibers. One advantage
of the linked linear series perspective is that with it, one can see the
necessary symmetries in the case of special determinant. Accordingly, in 
\cite{o-t2}, we specialize to the setting of rank $2$ and special 
determinant, and use these symmetries to prove the necessary dimensional
lower bounds on moduli spaces to get effective smoothing theorems. 
Combining these with the comparison results of the present paper, we are
then able to carry out limit linear series computations to prove 
existence of components of moduli spaces of $\fg^k_{2,d}$s having fixed
special determinant in a large family of examples. See \cite{o-t2} for
details. Zhang \cite{zh2} has proved further existence results in the
canonical determinant case using our smoothing theorem.

We now discuss the contents of the paper in more detail. We begin
in \S \ref{sec:prelim} with background on ``smoothing families'' -- the
families of nodal curves which we consider for our degenerations. In
\S \ref{sec:linked} we introduce the central new constructions, including
the definition of and construction of moduli stacks for type I and type II 
linked linear series. This provides two distinct generalizations of the
construction of \cite{os8}. Type I linked linear series have the advantage 
of being more suited to universal constructions, while it is the type II 
linked linear series which play the crucial role in the study of loci with
fixed special determinants. In \S \ref{sec:eht}, we provide a new 
perspective on Eisenbud-Harris-Teixidor limit linear series, yielding
a construction of a moduli stack which works in families. This is new
even in the rank-$1$ case, where we obtain for the first time a proper 
moduli space; the previous construction for families included only refined 
limit linear series. In the higher-rank case, the same construction
also shows that the Eisenbud-Harris-Teixidor limit linear series are 
locally closed in the natural ambient space, which was not clear from the
original definition; see Remark \ref{rem:eht-no-stack}. We then compare
the linked linear series construction to the limit linear series 
construction, proving that they are isomorphic on certain open loci.
In \S \ref{sec:chains}, we carry out a detailed analysis of the situation
for chains of curves, focusing on the locus of ``chain-adaptable''
limit linear series, which occurs ubiquitously in known families of
higher-rank limit linear series. This culminates in Corollary 
\ref{cor:compare}, which is our main comparison result and the cornerstone
of \cite{o-t2}. In \S \ref{sec:complement} we briefly develop several
complementary directions, including foundations for the fixed determinant
case, behavior of stability conditions, and results on specialization
of $\fg^k_{r,d}$s under degeneration. In Appendix \ref{app:pre-lg} we
define and study ``prelinked Grassmannians,'' generalizing the linked
Grassmannians of \cite{os8}, and finally, in Appendix 
\ref{app:pushforward-detl} we develop a generalization of determinantal
loci to pushforwards of coherent sheaves, which is used in the new
approach to constructing the moduli stack of Eisenbud-Harris-Teixidor
limit linear series.

Thus, although our primary motivation is the machinery and existence results 
for the special determinant case addressed in \cite{o-t2}, the present paper 
includes a number of new ideas which should prove useful more broadly.
For instance, the new description of classical Eisenbud-Harris limit linear
series simplifies existing proofs, and also suggests how one might approach 
a theory of limit linear series for curves not of compact type.

Finally, we briefly discuss the impetus for and implications of the
shift to working with stacks. The most compelling reason to work with
stacks is in order to be able to work with strictly semistable vector
bundles, for which coarse moduli spaces are very poorly behaved. A side
benefit is that even on stable loci, arguments become cleaner, without 
the need to pass to etale covers in order to produce universal families of 
vector bundles. From the point of view of comparison results between
constructions, the situation is not substantially subtler in working with
stacks than it would be for schemes (or even varieties). On loci for which
comparison morphisms are isomorphisms, one can typically deduce the stack
isomorphism from bijectivity on points, with little additional work. 
Conversely, on loci for which the morphisms are not isomorphisms, the fibers
are typically infinite. Dimension theory for stacks is somewhat subtler
than for schemes, but this is addressed by the theory developed in 
\cite{os21}, and consequently is handled entirely transparently in the 
present paper.

\subsection*{Acknowledgements}
I would like to thank Montserrat Teixidor i Bigas, David Eisenbud, and
Johan de Jong for helpful conversations.

\subsection*{Notational conventions}
Because we will use graphs rather extensively for notational purposes, we 
state our conventions. If $G$ is a (possibly directed) graph, we denote by
$V(G)$ and $E(G)$ the sets of vertices and edges of $G$, respectively. In
the directed case, for $e \in E(G)$, we denote by $t(e)$ and $h(e)$ the
tail and head of $e$, respectively.

If $X$ is an $S$-scheme, and $\sF$ a sheaf on $S$, to avoid introducing
notation for structure morphisms we denote by $\sF|_X$ the pullback of
$\sF$ to $X$.

If $X$ is a reducible curve, we will use subscripts to denote sheaves on
all of $X$ (of varying multidegrees) and spaces of global sections, and
we will use superscripts to denote sheaves on individual components of $X$
and spaces of global sections. We use script for vector bundles, and
roman letter for vector spaces.

\section{Preliminaries}\label{sec:prelim}

This section is devoted to a comprehensive treatment of the families of
curves we will consider, called ``smoothing families.'' A substantial
portion of our analysis involves the introduction and development of the 
``almost local'' condition on smoothing families, which guarantee that the 
dual graphs of fibers behave in a relatively simple manner. We conclude 
with a discussion of moduli stacks of vector bundles with prescribed
multidegree.

\subsection{Smoothing families}
The following definition differs only slightly from that of \cite{os8}:
because we work here with stacks rather than coarse moduli spaces, the
hypothesis on existence of sections is unnecessary. Also, because
vanishing conditions play a different role in higher rank, we omit
the choice of smooth sections of the smoothing family along which to 
impose such conditions.

\begin{defn}\label{def:smoothing-family}
A morphism of schemes $\pi: X \rightarrow B$
constitutes a \textbf{smoothing family} if:

\begin{Ilist}
\itm $B$ is regular and connected;
\itm $\pi$ is flat and proper;
\itm The fibers of $\pi$ are genus-$g$ curves of compact type;
\itm Each connected component $\Delta'$ of the singular locus of $\pi$ maps
isomorphically onto its scheme-theoretic image $\Delta$ in $B$, and
furthermore $\pi ^{-1} (\Delta)$ breaks into two (not necessarily
irreducible) components intersecting along $\Delta'$;
\itm Any point in the singular locus of $\pi$ which is smoothed
in the generic fiber is regular in the total space of $X$.
\end{Ilist}
\end{defn}

The following lemma is useful for constructing stacks of higher-rank
limit linear series.

\begin{lem}\label{lem:ample-exists} Locally on the base $B$, a smoothing 
family always carries a $\pi$-ample line bundle.
\end{lem}

\begin{proof} Given $y \in B$, by Corollary 9.6.4 of \cite{ega43}, a line 
bundle is $\pi$-ample in a neighborhood if and only if it has positive 
degree on every component of $X_y$. By hypothesis, $X$ breaks into 
components for every singularity of $X_y$ which is not generically smoothed, 
and line
bundles on $X$ are uniquely determined (up to twisting by pullbacks
from $B$) by their restrictions to these components. In addition, the
properties of a smoothing family are preserved under restriction to such
a component (see Lemma 3.2 (ii) of \cite{os8}), so it is enough
to treat the case that $X$ is irreducible, or equivalently, that 
every singularity of $X_y$ is smoothed generically. In this case, our
hypotheses give us that $X$ is regular. Then let $D$ be the closure of any
height-$1$ point of the generic fiber; we have then $\sO(D)$ is a line
bundle of strictly positive degree. Replacing $D$ by a suitable power, we
may assume the degree is at least as large as the number of components of
$X_y$. Now, for every singularity $\Delta'$ of $\pi$, we have by
hypothesis that $\pi^{-1}(\pi(\Delta'))$ breaks into components intersecting
along $\Delta'$; these components are each divisors on $X$, and twisting
$\sO(D)$ by them, we can redistribute the degree arbitrarily over the
components of $X_y$, and in particular can obtain a line bundle with
positive degree on each component, as desired.
\end{proof}

We take the opportunity to state a somewhat more general version of 
Theorem 3.4 of \cite{os8}. The proof is the same, with appropriate
considerations for working over a non-algebraically closed field as
described in \cite{os4}.

\begin{thm}\label{thm:grd-fam-exist}Let $X_0$ be a curve of compact type over 
a field $k$, and $\bar{P}_1, \dots, \bar{P}_n$ distinct smooth $k$-valued
points. Suppose that each node of $X_0$ is a $k$-valued point, and each
geometric component of $X_0$ is defined over $k$.
Then $X_0$ may be placed into a smoothing family $X/B$ with sections disjoint
smooth sections $P_i$
specializing to the $\bar{P}_i$, where $B$ is a curve over $k$, and where
the generic fiber of $X$ over $B$ is smooth.
\end{thm}

We also have the following standard structural statement:

\begin{prop}\label{prop:nodal-locus}
If $\pi:X \to B$ is a smoothing family, and $\Delta$ is the scheme-theoretic
image of a connected component of the non-smooth locus of $\pi$, then
$\Delta$ is regular.
\end{prop}

\begin{proof} If $\Delta=B$, this is automatic from the hypothesis that
$B$ is regular. Otherwise, let $\Delta' \subseteq X$ be the connected
component of the non-smooth locus of $\pi$ with image $\Delta$; the
proposition will follow from the definition of a smoothing family if
we show that $\Delta'$ is regular. On the other hand, the definition 
states that $X$ is regular along $\Delta'$. The statement being invariant
under completion, we may pass to complete local rings, in which case we
have that $\pi$ looks like $\Spec \hat{A}[u,v]/(uv-t) \to \Spec \hat{A}$,
where $\hat{A}$ is the relevant complete local ring of $B$, and 
$t \in \fm_{\hat{A}}$ (see for instance \S 2.23 of \cite{dj3}). It is then 
clear that $X$ being regular along $\Delta'$ implies that 
$t \not\in \fm_{\hat{A}}^2$, and the ideal of $\Delta'$ is cut out by 
$(u,v)$,
so the complete local ring of $\Delta'$ at the chosen point is
isomorphic to $\hat{A}/(t)$, which is regular, as desired.
\end{proof}

\subsection{The almost local condition}
We now devote some attention to the behavior of dual graphs in families
of curves. We begin with the following:

\begin{lem}\label{lem:closure}
Suppose that $\pi:X \to B$ is a smoothing family, and $y$ specializing
to $y'$ are points of $B$. Then if $\Gamma_y$ and $\Gamma_{y'}$ denote
the dual graphs of the fibers $X_y$ and $X_{y'}$ respectively, there
is a unique contraction map
$$\cl_{y,y'}:\Gamma_{y'} \to \Gamma_y$$
induced on vertices by associating to a component $Y'$ of $X_{y'}$ the 
component $Y$ of $X_y$ containing $Y'$ in its closure. The behavior of
$\cl_{y,y'}$ on edges is as follows: given an edge $e$ in $\Gamma_{y'}$ 
corresponding to a node $\Delta'$ in $X_{y'}$, if there is a node of
$X_y$ specializing to $\Delta'$, then $\cl_{y,y'}$ maps $e$ to the 
corresponding edge of $\Gamma_y$; otherwise, $e$ is contracted.

If also $y'$ specializes to some $y'' \in B$, then we have
$$\cl_{y,y''}=\cl_{y,y'} \circ \cl_{y',y''}.$$
\end{lem}

\begin{proof} 
Given $y$ specializing to $y'$ and
a component $Y'$ of $X_{y'}$, we first need to see that there exists a 
unique component $Y$ of $X_y$ containing $Y'$ in its closure.
Existence follows from flatness of $\pi$. To see uniqueness, let $Y$
and $Z$ be distinct components of $X_y$, and let $\Delta'$ be the
connected component of the non-smooth locus of $\pi$ containing a node
of $X_y$ separating $Y$ and $Z$. Let $\Delta$ be the image of $\Delta'$;
since $y,y' \in \Delta$, we can check uniqueness of specialization after
restricting to $\Delta$.
By definition of a smooothing family, we have that $\pi^{-1}(\Delta)$ breaks 
into (not necessarily irreducible) components $Y_{\Delta}$ and $Z_{\Delta}$
with $Y_{\Delta} \cap Z_{\Delta}=\Delta'$. Then the generic point of any
component of $X_y$ or $X_{y'}$ is contained in precisely one of 
$Y_{\Delta}$ and $Z_{\Delta}$, so it follows that 
$Y'$ can be in the closure of at most one of $Y$ and $Z$, as desired.
This gives $\cl_{y,y'}$ on the vertices of $\Gamma_{y'}$ and $\Gamma_y$;
it is then straightforward to check that we obtain a contraction map,
with the claimed behavior on edges.
Finally, associativity is clear from 
the definition.
\end{proof}

Although the next definition is slightly complicated, the idea behind it
is simple: it captures the condition that a smoothing family is 
combinatorially local, in the sense that there is a maximal dual graph
for the fibers of the family, and the dual graph of every fiber is 
naturally a contraction of the maximal one.

\begin{defn}\label{def:almost-local}
We say a smoothing family $\pi:X \to B$ is \textbf{almost
local} if the following condition is satisfied: if $\Delta'_1,\dots,\Delta'_m$ 
are the connected components of the non-smooth locus of $\pi$, with images 
$\Delta_1,\dots,\Delta_m$ in $B$, then there exists a $y_0 \in B$ such 
that for all $S,S' \subseteq \{1,\dots,m\}$, we have $\cap_{i \in S} \Delta_i$
non-empty, and for any irreducible components
$Z$ of $\cap_{i \in S} \Delta_i$ and $Z'$ of $\cap_{i \in S'} \Delta_i$,
every irreducible component of $Z \cap Z'$ contains $y_0$. 
\end{defn}

\begin{rem}\label{rem:almost-local}
Observe that the almost local hypothesis is trivially satisfied in the case 
that $\pi:X \to B$ has connected non-smooth locus, as is the situation in 
\cite{os8}. In addition, it is always satisfied locally on $B$: indeed,
given $y \in B$, we construct an open neighborhood on which $\pi$ is almost
local simply by removing any irreducible component of each 
$\cap_{i\in S} \Delta_i$ which does not contain $y$, and doing likewise
for intersections of pairs of components. 
\end{rem}

The following proposition says that an almost local smoothing family
admits a maximal dual graph in a natural way; this is the reason for
the hypothesis, as it will greatly simplify keeping track of dual graphs
and multidegrees.

\begin{prop}\label{prop:almost-local-maxl-graph} Suppose that $\pi:X \to B$
is an almost local smoothing family. Then there exists a graph $\Gamma$,
occuring as the dual graph of some fiber of $\pi$,
and, for every $y \in B$, a contraction 
$$\cl_y:\Gamma \to \Gamma_y,$$
where $\Gamma_y$ is the dual graph of the fiber $X_y$, satisfying
the following condition: if $y$ 
specializes to $y'$, we have 
\begin{equation}\label{eq:cl-factor} 
\cl_y = \cl_{y,y'} \circ \cl_{y'}.
\end{equation}
Moreover, up to automorphism of $\Gamma$, we have that $\Gamma$ and the
contractions $\cl_y$ are unique.
\end{prop}

Note that the uniqueness assertion in particular implies that the data
of $\Gamma$ and the $\cl_y$ contractions is independent of the choice of
$y_0$ from Definition \ref{def:almost-local}.

\begin{proof} We first fix some notation. Let $\Delta'_1,\dots,\Delta'_m$ 
be the connected components of the non-smooth locus, and 
$\Delta_1,\dots,\Delta_m$ their images in $B$. For $y \in B$, set 
$S_y \subseteq \{1,\dots,m\}$ to be the subset of $i$ such that 
$y \in \Delta_i$.
Observe that given $y$ specializing to $y'$, if $S_y=S_{y'}$ the contraction 
$\cl_{y',y}$ is a surjection of
trees with the same number of edges, and hence is necessarily an isomorphism.

Let $y_0$ be as in Definition \ref{def:almost-local}. Set 
$\Gamma=\Gamma_{y_0}$, with $\cl_{y_0}$ being the identity. Now, for any
$y\in B$, let $\tilde{y}$ be a generic point of $\cap_{i \in S_y} \Delta_i$
which specializes to $y$, so that $\cl_{\tilde{y},y}$ is an isomorphism.
It follows from the definition of almost local
that $\tilde{y}$ specializes to $y_0$, so we can set
\begin{equation}\label{eq:cl-formula}
\cl_y= \cl_{\tilde{y},y}^{-1} \circ \cl_{\tilde{y},y_0} \circ \cl_{y_0}.
\end{equation}
It remains to check that given $y$ specializing to $y'$, we have
\eqref{eq:cl-factor}. Let $\tilde{y}'$ be the generization of $y'$ used
to define $\cl_{y'}$, and let $Z,Z'$ be the closures of $\tilde{y}$ and
$\tilde{y}'$, respectively. Let $\bar{Z}$ be a component of $Z \cap Z'$
containing $y'$, and $\bar{y}$ its generic point. Then according to the
almost local hypothesis, we have that $\bar{y}$ specializes to $y_0$. 
Using associativity, we then have
\begin{align*}
\cl_y & = \cl_{\tilde{y},y}^{-1} \circ \cl_{\tilde{y},y_0} \circ \cl_{y_0}\\
& = \cl_{\tilde{y},y}^{-1} \circ \cl_{\tilde{y},\bar{y}} \circ
\cl_{\bar{y},y_0} \circ \cl_{y_0}\\
& = \cl_{\tilde{y},y}^{-1} \circ \cl_{\tilde{y},\bar{y}} \circ
\cl_{\bar{y},y_0} \circ \cl_{y_0}\\
& = \cl_{\tilde{y},y}^{-1} \circ \cl_{\tilde{y},y'} \circ
\cl_{\bar{y},y'}^{-1} \circ \cl_{\bar{y},y_0} \circ \cl_{y_0}\\
& = \cl_{y,y'} \circ \cl_{\bar{y},y'}^{-1} \circ 
\cl_{\bar{y},y_0} \circ \cl_{y_0}\\
& = \cl_{y,y'} \circ \cl_{\tilde{y}',y'}^{-1} \circ 
\cl_{\tilde{y}',y_0} \circ \cl_{y_0}\\
& = \cl_{y,y'} \circ \cl_{y'},
\end{align*}
as desired.

For the uniqueness assertion, observe that $\Gamma_{y_0}$ has $m$
edges, the maximal number possible among any $\Gamma_y$, so we must
have $\Gamma \cong \Gamma_{y_0}$. If we fix a choice of $\cl_{y_0}$,
which is defined precisely up to automorphism of $\Gamma$, we then 
have that \eqref{eq:cl-factor} implies that \eqref{eq:cl-formula} must
hold, so we have that $\cl_y$ is uniquely determined for all $y$,
as asserted.
\end{proof}

\begin{cor}\label{cor:almost-local-notn} Let $\pi:X \to B$ be an almost
local smoothing family, and $\Gamma$ and $\cl_y$ for $y \in B$ as given
by Proposition \ref{prop:almost-local-maxl-graph}. Then there is a 
bijection from $E(\Gamma)$ to the connected components of the non-smooth
locus of $\pi$ induced by sending an edge $e \in E(\Gamma)$ to
$$\Delta'_e:=\bigcup_{y:\cl_y\text{ does not contract $e$}} 
\Delta'_{\cl_y(e)},$$
where $\Delta'_{\cl_y(e)}$ denotes the node of $X_y$ corresponding to
$\cl_y(e)$.

Furthermore, if $\Delta_e \subseteq B$ denotes the image of $\Delta'_e$,
then given $v \in V(\Gamma)$ adjacent to $e$ there is a unique closed
subset $Y_{(e,v)} \subseteq \pi^{-1} (\Delta_e)$ such that for each
$y \in \Delta_e$, the fiber $(Y_{(e,v)})_y$ is equal to the union of
the components of $X_y$ corresponding to the vertices of $\Gamma_y$ lying
in the same connected component as $v$ in $\Gamma_y \smallsetminus \{e\}$.
\end{cor}

Thus, if $v,v'$ are the two vertices adjacent to an edge $e$, then
$Y_{(e,v)} \cup Y_{(e,v')}=\pi^{-1} (\Delta_e)$, and 
$Y_{(e,v)} \cap Y_{(e,v')}=\Delta'_e$.

\begin{proof} For the first assertion, we need to check that $\Delta'_e$
is in fact a connected component of the non-smooth locus of $\pi$, and
that the induced map is a bijection. Let $y_0$ be as in the definition
of almost local, and fix $e \in E(\Gamma)$. Let $\Delta'$ be the 
connected component of the non-smooth locus containing 
$\Delta'_{\cl_{y_0}(e)}$. We want to see that $\Delta'_e=\Delta'$. Then
given $y \in B$, let $\tilde{y}$ be a point generizing $y$ and $y_0$, and 
such that $\cl_{\tilde{y},y}$ is an isomorphism. First suppose that
$\cl_y$ does not contract $e$. Then we have that 
$\Delta'_{\cl_{\tilde{y}}(e)}$ specializes to both 
$\Delta'_{\cl_y(e)}$ and to $\Delta'_{\cl_{y_0}(e)}$. Thus, since $\Delta'$
is a connected component of the non-smooth locus, we conclude that 
$\Delta'_{\cl_y(e)} \in \Delta'$, so $\Delta'_e \subseteq \Delta'$. For
the opposite containment, suppose that $y$ is in the image of $\Delta'$.
Then since $\Delta'$ is a section of $\pi$ over its image, we see that
$\Delta'$ has a point of $X_{\tilde{y}}$ which specializes to
$\Delta'_{\cl_{y_0}(e)}$. It follows that $\cl_{\tilde{y}}$ does not
contract $e$, and then neither does $\cl_y$. Thus, $y$ is in the image of
$\Delta'_e$, and since $\Delta'$ is a section over its image and
$\Delta'_e \subseteq \Delta'$, we conclude that $\Delta'_e = \Delta'$,
as desired.

Next, to see that we have a bijection, we note that injectivity is trivial,
since for $e \neq e'$ we have that $\Delta'_e|_{X_{y_0}}$ and 
$\Delta'_{e'}|_{X_{y_0}}$ are distinct points. On the other hand,
the definition of almost local imposes that every connected component 
$\Delta'$ of the non-smooth locus of $\pi$ meets $X_{y_0}$, so meets some 
$\Delta'_e$. Since $\Delta'_e$ is also a connected component of the 
non-smooth locus, we conclude $\Delta'=\Delta'_e$, giving surjectivity.

For the assertion on $Y_{(e,v)}$, we have by hypothesis that 
$\pi^{-1}(\Delta_e)$ decomposes as $Y \cup Z$, with $Y \cap Z =\Delta'_e$;
exactly one of $Y$ or $Z$ contains the component of $X_{y_0}$ corresponding
to $\cl_{y_0}(v)$, so if we set $Y_{(e,v)}$ equal to this subset, we need
only verify that it has the desired form on every fiber. It is evident
from the construction that for each $y \in \Delta_e$, we have $Y_{(e,v)}|_y$
equal to the union of components of $X_y$ corresponding to a connected 
component of $\Gamma_y\smallsetminus\{e\}$, so we need only verify that
the connected component in question is the one containing $\cl_y(v)$. 
If $y$ specializes to $y_0$, this is immediate from the fact that
$\cl_y=\cl_{y,y_0} \circ \cl_{y_0}$, and in the general case it follows
by choosing some $\tilde{y}$ specializing to both $y$ and $y_0$ as above.
\end{proof}

\subsection{Multidegrees of vector bundles}
We conclude this section with some background on multidegrees of vector
bundles for smoothing families, and the associated moduli stacks.

\begin{defn}\label{def:multidegree}
Given an almost-local smoothing family $\pi:X \to B$, let
$\Gamma$ and $\cl_y$ be as in Proposition \ref{prop:almost-local-maxl-graph}.
Given also a $B$-scheme $S$, a vector bundle $\sE$ on $X \times_B S$
has \textbf{multidegree} $w=(i_v)_{v \in V(\Gamma)}$ if, for 
every $s\in S$ and every component $Y$ of the fiber $(X \times_B S)_s$,
we have 
$$\deg \sE|_Y = \sum_{v \in V(\Gamma):\cl_y (v)=v_Y} i_v,$$
where $y$ is the image of $s$ in $B$, and $v_Y$ is the vertex of $\Gamma_y$
corresponding to $Y$.

We denote the groupoid of vector bundles of rank $r$ and multidegree 
$w$ on $X/B$ by $\cM_{r,w}(X/B)$.
\end{defn}

Note that our hypotheses on smoothing families imply in particular that
irreducible components of fibers are geometrically irreducible, so in 
the above definition, the dual graphs of $(X\times_B S)_s$ and $X_y$ are 
canonically identified. 

For lack of a suitable reference, we include the proof of the following
basic fact.

\begin{prop}\label{prop:multideg-stack} We have that 
$\cM_{r,w}(X/B)$ is an open substack of the moduli stack
$\cM_{r,d}(X/B)$ of vector bundles of rank $r$ and degree 
$d:=\sum_{v \in V(\Gamma)} i_v$. In particular, it is an Artin stack,
locally of finite type over $B$.
\end{prop}

\begin{proof}
Given $w=(i_v)_v$, we can associate to every pair
$(e,v)$ with $e \in E(\Gamma)$ and $v \in V(\Gamma)$ adjacent to $e$
an integer $i_{(e,v)}$: let $C$ be the connected component of 
$\Gamma \smallsetminus \{e\}$ containing $v$, and set 
$i_{(e,v)}=\sum_{v \in C} i_v$. We then observe that $\sE$ has multidegree 
$w$ if
and only if for every $(e,v)$ as above, and every fiber of $X \times _B S$ 
supported over $\Delta_e$, the bundle $\sE$ has degree $i_{(e,v)}$ on the 
preimage of $Y_{(e,v)}$.
But now, it is easy to see that this is an open 
condition. Indeed, we have that $Y_{(e,v)}$ is flat over $\Delta_e$ by the 
argument for Lemma 3.2 (ii) of \cite{os8}, so degree is locally constant. 
Thus, if the condition is satisfied on a single fiber $X_s$ of 
$X \times_B S$, it is necessarily satisfied at all points in the connected
component of $S$ except possibly for a finite union of 
preimages of subsets $\Delta_e$, each of which is closed. 
\end{proof}

\section{The stacks of higher-rank linked linear series}\label{sec:linked}

We now come to the foundational definitions and constructions for our
new spaces of limit linear series. In order to distinguish our terminology
from that of Eisenbud-Harris-Teixidor, we refer to our objects as
``linked linear series.'' This is consistent with \cite{o-t1}, but deviates
from \cite{os8}. The idea of a linked linear series is to consider a
collection of vector bundles obtained by various twists from a single
vector bundle, together with spaces of global sections on each such twist
satisfying a linkage condition under certain natural maps. We will introduce 
two variants of this idea, each generalizing the single construction of
\cite{os8}, but for each one the process is the same: we introduce a directed 
graph to parametrize the twists we wish to consider, describe twisting line
bundles and associated morphisms, and then use these to define the linked
linear series we wish to consider.

For the second construction, we will consider infinite graphs. This makes
the definitions more transparent because we do not have to consider edges 
near the boundary as special cases, but it will come at the cost that in
the foundational theorem, we will have to prove that the resulting 
construction can be described equivalently on a finite graph.

\subsection{The underlying graphs}
The basic situation we will consider is the following:

\begin{sit}\label{sit:basic} Let $\pi:X \to B$ be an almost local smoothing 
family. Let $\Gamma$ be the graph obtained from Proposition
\ref{prop:almost-local-maxl-graph}, with associated contractions 
$\cl_y$ for $y \in B$. Let $\Delta'_e$, $\Delta_e$ and $Y_{(e,v)}$ be as
in Corollary \ref{cor:almost-local-notn}.

Let $r,d,k$ be positive integers, and fix also 
integers $b$ and $d_v$ for each $v \in V(\Gamma)$, satisfying
\begin{equation}\label{eq:deg-sum} 
\sum_{v \in V(\Gamma)} d_v - |E(\Gamma)|rb=d.
\end{equation}
\end{sit}

We will associate two directed graphs $G_{\I}$ and $G_{\II}$ to $\Gamma$, 
corresponding to our two generalizations of the construction of \cite{os8}.
The vertices of both graphs will correspond to choices of multidegrees, and
hence lie in $\ZZ^{V(\Gamma)}$. The choice of the $d_v$ (which together
determine $b$) serves two purposes: it determines a congruence class modulo
$r$ for the multidegrees in question, and it also determines in some sense
where our ``extremal'' vertices lie. This was not necessary in the classical
rank-$1$ case because no line bundle of negative degree has nonzero global
sections. We will see in Proposition \ref{prop:increase-b}
below that -- at least in the type II case -- we can always increase the 
$d_v$, and this will have the effect of imbedding the original moduli space 
as an open substack into the newer one.

\begin{defn}\label{def:graphs} Suppose we are in the above situation.
Let $V(G_{\II}) \subseteq \ZZ^{V(\Gamma)}$ consist of vectors
$(i_v)_{v \in V(\Gamma)}$ satisfying:
\begin{Ilist}
\itm $\sum_v i_v=d$, and
\itm $i_v \equiv d_v \pmod{r}$ for all $v \in V(\Gamma)$,
\end{Ilist}
and let $V(G_{\I}) \subseteq V(G_{\II})$ be the subset on which we further have
\begin{Ilist}[2]
\itm $i_v \geq d_v-rb$ for all $v \in V(\Gamma)$, and
\itm $i_v = d_v-rb$ for all but at most two $v \in V(\Gamma)$, with the 
vertices on which we have strict inequality required to be adjacent in
$\Gamma$.
\end{Ilist}

Let $G_{\I}$ be the directed graph with vertex set $V(G_{\I})$, and with an 
edge from $h$ to $h'$ in $V(G_{\I})$ precisely when there exist
adjacent vertices $v,v' \in V(\Gamma)$ such that $h-h'$ is $\pm r$
in index $v$, $\mp r$ in index $v'$, and $0$ elsewhere.

Let $G_{\II}$ be the directed graph with vertex set $V(G_{\II})$, and an edge
from $h$ to $h'$ if there is a vertex $v \in V(\Gamma)$ of valence $\ell$
such that $h'-h$ is $-\ell r$ in index $v$, is $r$ in index $v'$ for
each $v'$ adjacent to $v$, and is $0$ elsewhere.
\end{defn}

\begin{figure}
\centering
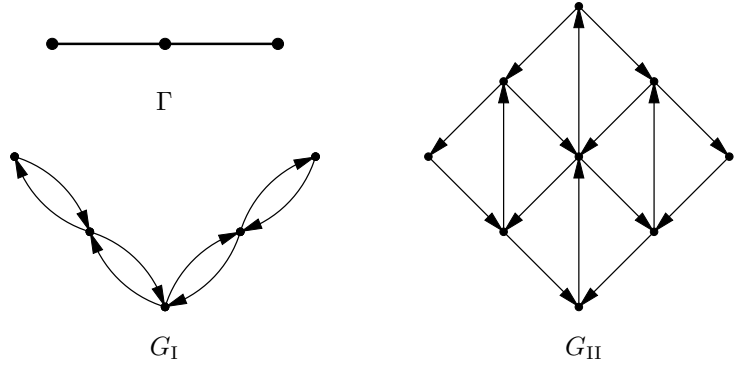
\caption{An example of Definition \ref{def:graphs}, with $b=2$. We have
restricted $G_{\II}$ to the finite region used in the definition of
$\bar{G}_{\II}$ in Definition \ref{def:bar-g} below.}
\end{figure}

Note that although $G_{\I}$ is taken to be directed, for any edge from 
$h$ to $h'$ there is a corresponding edge from $h'$ to $h$. Thus, $G_{\I}$
has no more information than the underlying undirected graph, but we use the
directed version for convenience of notation.

By construction, the vertices of either $G_{\I}$ and $G_{\II}$ naturally
induce multidegrees on fibers of $\pi$, and we will sometimes treat them
as multidegrees without further comment.

For both $G_{\I}$ and $G_{\II}$, the edges (or more generally paths) coming 
out of a given vertex can be described simply in terms of $\Gamma$, and it 
will be convenient to introduce notation expressing this relationship.

\begin{notn}\label{not:v-e-P}
Given an edge $\e \in G_{\I}$, denote by $v_{\I}(\e)$ and 
$e_{\I}(\e)$ the vertex and (adjacent) edge of $\Gamma$ such that 
following $\e$ decreases the coordinate in index $v_{\I}(\e)$, and 
increases the coordinate indexed by the other vertex adjacent to 
$e_{\I}(\e)$.
Given $w \in V(G_{\I})$, and a sequence 
$(e_1,v_1),\dots,(e_m,v_m)$
of edges together with adjacent vertices in $\Gamma$,
denote by $P(w,(e_1,v_1),\dots,(e_m,v_m))$ the path 
in $G_{\I}$ (if it exists) starting at $w$ and consisting of edges 
$\e_1,\dots,\e_m$ such that $v_{\I}(\e_i)=v_i$ and $e_{\I}(\e_i)=e_i$ for 
each $i$.

Given an edge $\e \in G_{\II}$, denote by $v_{\II}(\e)$ the corresponding
vertex of $\Gamma$. Given $w \in V(G_{\II})$ and a sequence 
$v_1,\dots,v_m$
of (not necessarily distinct) vertices of $\Gamma$, denote by 
$P(w,v_1,\dots,v_m)$ the path in $G_{\II}$ starting at $w$ and 
consisting of edges $\e_1,\dots,\e_m$ such that $v_{\II}(\e_i)=v_i$ 
for each $i$.
\end{notn}

Observe that for any given collection of $(e_1,v_1),\dots,(e_m,v_m)$, we
may not obtain
a path $P(w,(e_1,v_1),\dots,(e_m,v_m))$ in $G_{\I}$, due to our conditions
on the vertices of $G_{\I}$.
Because $G_{\I}$ is a tree, minimal paths between any two vertices in it
are unique. On the other hand, the endpoint of $P(w,v_1,\dots,v_m)$ is 
independent of the ordering of $v_1,\dots,v_m$. In fact, we can be very 
precise about the extent to which (minimal) paths are
determined by their endpoints in $G_{\I}$ and $G_{\II}$, as follows.

\begin{prop}\label{prop:min-paths} Minimal paths are unique in $G_{\I}$,
and $P(w,(e_1,v_1),\dots,(e_m,v_m))$ is minimal if and only if 
no edge $e$ of $\Gamma$ with adjacent vertices $v,v'$ has both $(e,v)$ 
and $(e,v')$ appearing among the $(e_i,v_i)$. 

A path $P(w,v_1,\dots,v_m)$ in $G_{\II}$ is minimal if and only if not
every vertex of $\Gamma$ occurs as one of the $v_i$. In this case, the
resulting path is unique up to reordering the $v_i$.

More generally, paths 
$P(w,(e_1,v_1),\dots,(e_m,v_m))$ and 
$P(w,(e'_1,v'_1),\dots,(e'_{m'},v'_{m'}))$ in $G_{\I}$
have the same endpoint if and only if the multisets 
$\{(e_1,v_1),\dots,(e_m,v_m)\}$ and 
$\{(e'_1,v'_1),\dots,(e'_{m'},v'_{m'})\}$ 
differ by unions of sets of the form 
$\{(e,v),(e,v')\}$, where $v$ and $v'$
are the two edges adjacent to $e$.

Similarly, two paths $P(w,v_1,\dots,v_m)$ and 
$P(w,v'_1,\dots,v'_{m'})$ in $G_{\II}$
have the same endpoint if and only if the multisets 
$\{v_1,\dots,v_m\}$
and $\{v'_1,\dots,v'_{m'}\}$ differ by a multiple of $V(\Gamma)$.
\end{prop}

The proof is straightforward, and left to the reader.

Next, observe that conditions (I) and (III) together imply the following 
description of the extremal vertices of $G_{\I}$:

\begin{prop}\label{prop:extremal-degs} 
For any $(i_v)_{v \in V(\Gamma)} \in G_{\I}$ we 
have $i_v \leq d_v$, with
equality if and only if $i_{v'}=d_{v'}-rb$ for all $v' \neq v$ in $V(\Gamma)$.
\end{prop}

\begin{notn}\label{not:wv}
Given $v \in V(\Gamma)$, denote by $w_v$ the vertex of $G_{\I}$ 
with coordinates $i_v=d_v$, and $i_{v'}=d_{v'}-rb$ for all $v' \neq v$. 
\end{notn}

Although there are no explicit extremal vertices in the definition of 
$G_{\II}$, the vertices $w_v$ are in some sense extremal for $G_{\II}$
as well; see Proposition \ref{prop:type-ii-finite}. We use the convention
of the infinite graph $G_{\II}$ in order to keep definition as simple as
possible.

\begin{rem}\label{rem:bad-defn}
Perhaps the most natural generalization of the construction 
of \cite{os8} is a composite of the type I and type II constructions,
where the multidegrees are allowed to be arbitrary as in type II, but
the edges between vertices are defined in terms of nodes of the curve
as in type I. However, this approach turns out not to be well behaved,
for the reason that the maps are simply too degenerate for linkage to
impose enough requirements. Indeed, with this approach we see that the 
sort of behavior described in Example \ref{ex:bad-compare} below 
occurs even for refined $g^0_d$s on curves with three components. 
\end{rem}

\subsection{Linked linear series: type I}
In the following, we will define, given a vector bundle $\sE$ of multidegree
$w_0 \in V(G_{\I})$, natural twists $\sE_w$ for all $w \in V(G_{\I})$, along
with maps (defined up to scalar, and existing locally on $B$)
$$f_{\e}:\pi_* \sE_w \to \pi_* \sE_{w'}$$ 
for any edge $\e$ of $G_{\I}$ going from $w$ to $w'$. Given these, we
define type-I linked linear series as follows.

\begin{defn}\label{def:grd-space-i} In Situation \ref{sit:basic}, let 
$G_{\I}$ be as in Definition \ref{def:graphs}, and choose a 
vertex 
$w_0 \in V(G_{\I})$. The moduli groupoid $\cG^{k,\I}_{r,d,d_{\bullet}}(X/B)$ 
of \textbf{type-I linked linear series} is the category fibered 
in groupoids over $\BSch$ whose objects consist of tuples 
$(S,\sE,(\sV_w)_{w \in V(G_{\I})})$, where $S$ is
a $B$-scheme, $\sE$ is a vector bundle of rank $r$ and multidegree
$w_0$ on $X \times_B S$, and $\sV_w$ is a rank-$k$ subbundle (in the
sense of Definition \ref{def:subbundle})
of $\pi_* (\sE_w)$, satisfying the following conditions:
\begin{Ilist}
\itm for every pair $(e,v)$ of an edge and adjacent vertex of $\Gamma$, 
and every $z \in S$ with image $y \in B$ lying in $\Delta_e$,
we have
$$H^0(Y,\sE_{w_v}|_{Y}(-(b+1)(\Delta')))=0,$$
where $Y$ denotes the component of the fiber $X_z$ corresponding to $\cl_y(v)$,
and $\Delta'$ denotes the node corresponding to $\cl_y(e)$;
\itm for every edge $\e$ in $G_{\I}$, let 
$w$ be the tail and $w'$ the head. Then we require that
$$f_{\e}(\sV_w) \subseteq \sV_{w'}.$$
\end{Ilist}
\end{defn}

By definition of smoothing family, the irreducible components of $X_y$
remain irreducible in $X_z$, so condition (I) makes sense.
Note that, as in \cite{os8}, although $f_{\e}$ is defined only locally on 
$B$, the condition that $f_{\e}(\sV_w) \subseteq \sV_{w'}$ is invariant 
under scalar multiplication by $\sO_B^*$, so the resulting closed condition 
is defined on all of $B$.
Similarly, while the definition depends \textit{a priori} on the choice of 
$w_0$, this is largely a matter of convenience: different choices of $w_0$
will yield equivalent groupoids.

In order to complete Definition \ref{def:grd-space-i}, we 
begin by describing the relevant twisting line bundles. 

\begin{notn}\label{not:twisting-bundles} 
For every pair $(e,v)$ of an edge $e \in E(\Gamma)$ and an adjacent vertex
$v$, denote by $\sO_{(e,v)}$ the line bundle on $X$ obtained as follows:
write
$$X|_{\Delta_e}=Y_{(e,v)} \cup Z_{(e,v)},$$ 
where $Z_{(e,v)}=Y_{(e,v')}$ with $v'$ the other vertex adjacent to $e$.

Now, if $\Delta_e\neq B$, we have
$Y_{(e,v)}$ a (necessarily Cartier) divisor in $X$, and we set 
$\sO_{(e,v)}=\sO_X(Y_{(e,v)})$.

On the other hand, if $\Delta_e=B$, then line bundles on $X$ are uniquely
determined by their restrictions to $Y_{(e,v)}$ and $Z_{(e,v)}$, and
we define $\sO_{(e,v)}$ to be $\sO_{Y_{(e,v)}}(-\Delta'_e)$ on
$Y_{(e,v)}$ and to be $\sO_{Z_{(e,v)}}(\Delta'_e)$ on $Z_{(e,v)}$.

Given $w,w' \in V(G_{\I})$, let $P$ be a minimal directed path in $G_{\I}$
from $w$ to $w'$. Write $P=P(w,(e_1,v_1),\dots,(e_m,v_m))$, and set
$$\sO_{w,w'}=\bigotimes_{i=1}^m \sO_{(e_i,v_i)}.$$

Finally, given $w_0 \in V(G_{\I})$, and $\sE$ an $S$-valued point of
$\cM_{r,w_0}(X/B)$, for $w \in V(G_{\I})$ write
$$\sE_w:=\sE \otimes \sO_{w_0,w}|_{X \times_B S}.$$
\end{notn}

Exactly as in \cite{os8}, locally on $B$ we also have maps, unique up to 
non-zero scalar, defined as follows:

\begin{notn}\label{not:fe}
Given a pair $(e,v)$ of an edge $e \in E(\Gamma)$ and an adjacent vertex
$v$, we first observe we have a canonical morphism 
$$\theta_{(e,v)}:\sO_X \to \sO_{(e,v)}.$$
When $\Delta_e \neq B$, 
this is defined to be the canonical inclusion. On the other hand,
when $\Delta_e = B$, this is defined to be $0$ on $Y_{(e,v)}$, and the
canonical inclusion on $Z_{(e,v)}$.
Now, let $v'$ be the other vertex adjacent to $e$. Then (locally on $B$,
in the case $\Delta_e \neq B$) we have that 
$\sO_{(e,v)} \otimes \sO_{(e,v')} \cong \sO_X$. If we fix such an isomorphism
(unique up to an element of $\sO_B^*$),
we obtain a morphism
$$\theta'_{(e,v)}:\sO_{(e,v)} \overset{\theta_{(e,v')}}{\to} 
\sO_{(e,v)} \otimes \sO_{(e,v')} \to \sO_X.$$

Finally, fix $w_0 \in V(G_{\I})$. For any edge $\e$ in $G_{\I}$, let 
$w$ be the tail and $w'$ the head, and let $e=e_{\I}(\e) \in E(\Gamma)$.
If $\sE$ is an $S$-valued point of $\cM_{r,w_0}(X/B)$,
then we either have $\sE_w = \sE_{w'} \otimes \sO_{(e,v)}|_{X \times_B S}$
or $\sE_{w'}=\sE_w \otimes \sO_{(e,v)}|_{X \times_B S}$, where $v$ is a 
vertex adjacent
to $e$. Thus, using $\theta_{(e,v)}$ or $\theta'_{(e,v)}$, and pushing forward
under $\pi$, we obtain a morphism
$$f_{\e}: \pi_* \sE_w \to \pi_* \sE_{w'}.$$
\end{notn}

This completes the definition of type-I linked linear series.
Observe that by construction, $\theta_{(e,v)} \circ \theta'_{(e,v)}$ is
equal to scalar multiplication by a scalar cutting out $\Delta_e$, and
similarly for $f_{\e} \circ f_{\e'}$, if $\e'$
is the edge with the same adjacent vertices, but going in the opposite
direction as $\e$.

\subsection{Linked linear series: type II}
We now move on to the second construction, which differs from the first
construction in that it both requires more restrictive hypotheses,
and also depends on additional choices. The additional hypothesis is the 
following.

\begin{sit}\label{sit:basic-ii} In Situation \ref{sit:basic}, suppose 
further that for all $e \in E(\Gamma)$, the closed subschemes 
$\Delta_e \subseteq B$ agree. Denote this common closed subscheme by 
$\Delta$, and for every $v \in V(\Gamma)$, let $Y_v$ be the irreducible 
component of $\pi^{-1}(\Delta)$ corresponding to $v$.
\end{sit}

Note that the additional condition of Situation \ref{sit:basic-ii} is
tautologically satisfied for the case of curves with at most two components,
or if $B=\Spec F$ for some field $F$.
In addition, Proposition \ref{prop:nodal-locus} implies that if the 
$\Delta_e$ agree set-theoretically, they also agree
scheme-theoretically, so in particular if $B$ is the spectrum of a DVR
and the generic fiber is smooth, the additional condition of Situation
\ref{sit:basic-ii} is likewise satisfied. 

As in the type-I case, under the additional hypothesis of Situation
\ref{sit:basic-ii} we will define, given a vector bundle $\sE$ of 
multidegree $w_0 \in V(G_{\II})$, natural twists $\sE_w$ for all 
$w \in V(G_{\II})$, along with maps (defined up to scalar, and existing 
locally on $B$)
$$f_{\e}:\pi_* \sE_w \to \pi_* \sE_{w'}$$ 
for any edge $\e$ of $G_{\II}$ going from $w$ to $w'$. The additional
data required to define the maps $f_{\e}$ will be a suitable collection
of sections of our twisting bundles $\sO_v$. Given the above definitions, 
we define type-II linked linear series as follows.

\begin{defn}\label{def:grd-space-ii} In Situation \ref{sit:basic-ii}, let 
$G_{\II}$ be as in Definition \ref{def:graphs}, and choose a vertex 
$w_0 \in V(G_{\II})$. In addition, for each $v \in V(\Gamma)$, choose
a morphism $\theta_v: \sO_X \to \sO_v$ which vanishes precisely on $Y_v$.
The moduli groupoid $\cG^{k,\II}_{r,d,d_{\bullet}}(X/B,\theta_{\bullet})$ 
of \textbf{type-II linked linear series} is the category fibered 
in groupoids over $\BSch$ whose objects consist of tuples 
$(S,\sE,(\sV_w)_{w \in V(G_{\II})})$, where $S$ is
a $B$-scheme, $\sE$ is a vector bundle of rank $r$ and multidegree
$w_0$ on $X \times_B S$, and $\sV_w$ is a rank-$k$ subbundle (in the
sense of Definition \ref{def:subbundle})
of $\pi_* (\sE_w)$, satisfying the following conditions:
\begin{Ilist}
\itm for every pair $(e,v)$ of an edge and adjacent vertex of $\Gamma$,
and every $z \in S$ with image $y \in B$ lying in $\Delta_e$,
we have
$$H^0(Y,\sE_{w_v}|_{Y}(-(b+1)(\Delta')))=0,$$
where $Y$ denotes the component of the fiber $X_z$ corresponding to $\cl_y(v)$,
and $\Delta'$ denotes the node corresponding to $\cl_y(e)$;
\itm for every edge $\e$ in $G_{\II}$, let 
$w$ be the tail and $w'$ the head. Then we require that
$$f_{\e}(\sV_w) \subseteq \sV_{w'}.$$
\end{Ilist}
\end{defn}

In the context of this definition, the parameters $d_v$ and $b$ are
somewhat artificial and included mainly for consistency. In particular,
it is always possible to increase them and obtain open immersions on the
resulting moduli spaces; see \S \ref{sec:b}. As in the type I case, the
choice of $w_0$ doesn't change the resulting groupoid.

We now describe how to construct our twisting bundles and maps in order to
complete Definition \ref{def:grd-space-ii}.

\begin{notn}\label{not:twisting-bundles-ii} 
For every vertex $v \in V(\Gamma)$, denote by $\sO_v$ the line bundle on 
$X$ obtained as follows: 

If $\Delta\neq B$, we have
$Y_{v}$ a Cartier divisor in $X$, and we set 
$\sO_{v}=\sO_X(Y_v)$.

On the other hand, if $\Delta=B$, let $Z_1,\dots,Z_m$ be the closures in $X$
of the connected components of $X\smallsetminus Y_v$, and for $i=1,\dots,m$,
let $\Delta'_i$ be the node $Z_i \cap Y_v$. 
We then define $\sO_{v}$ to be $\sO_{Y_{v}}(-\sum_i \Delta'_i)$ on
$Y_{v}$ and to be $\sO_{Z_i}(\Delta'_i)$ on each $Z_i$.

Given $w,w' \in V(G_{\II})$, let $P=(e_1,\dots,e_m)$ be a minimal directed 
path in $G_{\II}$ from $w$ to $w'$. Let $v_i = v_{\II}(e_i)$ for
$i=1,\dots,m$, and set
$$\sO_{w,w'}=\bigotimes_{i=1}^m \sO_{v_i}.$$

Finally, given $w_0 \in V(G_{\II})$, and an $S$-valued point $\sE$ of
$\cM_{r,w_0}(X/B)$, for $w \in V(G_{\II})$ write
$$\sE_w:=\sE \otimes \sO_{w_0,w}|_{X \times_B S}.$$
\end{notn}

\textit{A priori}, the notation $\sO_{w,w'}$ and $\sE_w$ is ambiguous,
since it is used in both the type-I and type-II constructions. However, 
we will see in Proposition
\ref{prop:two-types-compat} below that there is in fact no ambiguity.

\begin{notn}\label{not:type-ii-maps}
For each $v \in V(\Gamma)$, suppose we fix a morphism
$$\theta_{v}:\sO_X \to \sO_{v}$$
vanishing precisely on $Y_v$.

Next, observe that $\bigotimes_{v \in V(\Gamma)} \sO_{v} \cong \sO_X$. Fixing 
such an isomorphism (unique up to an element of $\sO_B^*$),
we obtain a induced morphism
$$\theta'_{v}:\bigotimes_{v' \neq v} \sO_{v'} \overset{\theta_v}{\to} 
\bigotimes_{v' \in V(\Gamma)} \sO_{v'} \to \sO_X.$$

Finally, fix $w_0 \in V(G_{\II})$. For any edge $\e$ in $G_{\II}$, let 
$w$ be the tail and $w'$ the head, and let $v$ be the associated edge of
$\Gamma$. If $\sE$ is an $S$-valued point of $\cM_{r,w_0}(X/B)$
then we either have $\sE_{w'} = \sE_{w} \otimes \sO_{v}|_{X \times_B S}$
or $\sE_{w}=\sE_{w'} \otimes \bigotimes_{v' \neq v} \sO_{v'}|_{X\times_B S}$.
Thus, using $\theta_{v}$ or $\theta'_{v}$, and pushing forward
under $\pi$, we obtain a morphism
$$f_{\e}: \pi_* \sE_w \to \pi_* \sE_{w'}.$$
\end{notn}

Note that in the case $\Delta \neq B$, we have that each $\theta_v$ and hence
each $f_{\e}$ is unique up to scalar. However, this is not the case when
$\Delta = B$, as, in the notation of Notation \ref{not:twisting-bundles-ii},
we can scale $\theta_v$ independently on each of the $Z_i$.

\subsection{Foundational results}
The following proposition is straightforward from the definitions,
and implies that our notation is consistent, and
more importantly, that we have a forgetful morphism from type II linked
linear series to type I linked linear series.

\begin{prop}\label{prop:two-types-compat}
In Situation \ref{sit:basic-ii}, given $w,w' \in V(G_{\I})$, we have that 
the line bundles $\sO_{w,w'}$
given in Notation \ref{not:twisting-bundles} and \ref{not:twisting-bundles-ii}
are isomorphic. Consequently, for any $w \in V(G_{\I})$ the vector bundles
$\sE_w$ in Notation \ref{not:twisting-bundles} and 
\ref{not:twisting-bundles-ii} are likewise isomorphic.

In addition, if $\e \in E(G_{\I})$ is an edge with tail $w$ and head $w'$,
for any choice of $\theta_{\bullet}$ as in Definition \ref{def:grd-space-ii},
and any minimal path $P=(\e_1,\dots,\e_n)$ from $w$ to $w'$ in $G_{\II}$, we 
have that $f_{\e}$ agrees with $f_{\e_n} \circ \dots \circ f_{\e_1}$ up to
scalar multiplication by an element of $\sO_B^*$.
\end{prop}

Note in particular that in the last part of the proposition, each $f_{\e_i}$
depends in general on the choice of $\theta_{\bullet}$, but we see that the
composition $f_{\e_n} \circ \dots \circ f_{\e_1}$ does not.

We conclude:

\begin{cor}\label{cor:forget} There is a forgetful morphism
$$\cG^{k,\II}_{r,d,d_{\bullet}}(X/B,\theta_{\bullet}) \to 
\cG^{k,\I}_{r,d,d_{\bullet}}(X/B).$$
\end{cor}

In order to state the foundational theorem, it is helpful to introduce
the following open substack of $\cM_{r,w_0}(X/B)$.

\begin{notn}\label{not:bundles-substack} 
Let $X/B$ be an almost local smoothing 
family, and $k,r,d_{\bullet}$ as in Situation \ref{sit:basic}.
Then $\cM_{r,w_0,d_{\bullet}}(X/B)$ denotes the subgroupoid of
$\cM_{r,w_0}(X/B)$ consisting of vector bundles $\sE$ on $X \times_B S$ 
such that for every pair $(e,v)$ of an edge and adjacent vertex of $\Gamma$,
and every $z \in S$ with image $y \in B$ lying in $\Delta_e$,
we have
$$H^0(Y,\sE_{w_v}|_{Y}(-(b+1)(\Delta')))=0,$$
where $Y$ denotes the component of the fiber $X_z$ corresponding to $\cl_y(v)$,
and $\Delta'$ denotes the node corresponding to $\cl_y(e)$;
\end{notn}

\begin{lem}\label{lem:substack-open} Let $\sE$ on $X_y$ be a $K$-valued point
of $\cM_{r,w_0,d_{\bullet}}(X/B)$. Then for all $v \in V(\Gamma)$, if $Y$
denotes the component of $X_y$ corresponding to $v$, the restriction
map 
$$H^0(X_y,\sE_{w_v}) \to H^0(Y,\sE_{w_v}|_Y)$$
is injective.

The groupoid
$\cM_{r,w_0,d_{\bullet}}(X/B)$ is an open substack of $\cM_{r,w_0}(X/B)$.
\end{lem}

\begin{proof} 
For the first assertion, suppose we have a section in the kernel of 
restriction to $Y$, which we wish to show is zero on all of $X_y$. We
induct on the number of components of $X_y$, with the base case that $Y=X_y$
being trivial. Now, if $Y' \neq Y$ is another component of $X_y$, with $v'$ a
corresponding vertex of $\Gamma$ being a leaf, then there is a unique node
$\Delta'$ on $Y'$, and
$\sE_{w_v}|_{Y'}=\sE_{w_{v'}}|_{Y'}(-b\Delta')$.
The situation being compatible with restriction to connected
subcurves, we may assume by induction that our section vanishes on all
components other than $Y'$, and in particular at $\Delta'$,
and thus the restriction to $Y'$ is a section of
$\sE_{w_{v'}}|_{Y'}(-(b+1)\Delta')$, which must be zero by definition of 
$\cM_{r,w_0,d_{\bullet}}(X/B)$.

In order to prove openness, it is enough to check that the
condition is constructible, and closed under generization. Constructibility
is clear: considering dual graphs of fibers gives a stratification of $B$ 
by locally closed subschemes, and within each stratum, the condition in
question fails on a closed subset, by semicontinuity of $h^0$. Note that
on each stratum, each component of $X$ has irreducible fibers and is proper
and flat (for flatness, see the argument for Lemma 3.2 (ii) of \cite{os8}).

In order to prove closure under generization, it is enough to consider 
the case that the base is the spectrum of a DVR. In this case, suppose
the condition is satisfied on the special fiber $X_0$; we will show it is
also satisfied on the generic fiber $X_{\eta}$. Given
$(e,v)$, let $Y_{\eta}$ and $\Delta'_{\eta}$ (respectively, $Y_0$ and
$\Delta'_0$) be the relevant component and node on $X_{\eta}$ 
(respectively, on $X_0$). Let $Z$ denote the restriction to $X_0$ of the 
closure of $Y_{\eta}$ in $X$; this contains $Y_0$, but may also contain 
additional components of $X_0$. Since we have assumed our condition is 
satisfied on $X_0$, we have that 
$H^0(Y_0,\sE_{w_v}|_{Y_0}(-(b+1)(\Delta'_0))) = 0$,
and also, by the first statement of the lemma, that the restriction map
$$H^0(Z,\sE_{w_v}|_{Z}) \to H^0(Y_0,\sE_{w_v}|_{Y_0})$$
is injective. This implies that
$$H^0(Z,\sE_{w_v}|_{Z}(-(b+1)(\Delta'_0))) = 0,$$
which by semicontinuity implies that
$$H^0(Y_{\eta},\sE_{w_v}|_{Y_{\eta}}(-(b+1)(\Delta'_{\eta}))) = 0$$
as well, as desired. 
\end{proof}

Thus, the natural forgetful morphism
$$\cG^{k,\I}_{r,d,d_{\bullet}}(X/B) \to \cM_{r,w_0}(X/B)$$
factors through $\cM_{r,w_0,d_{\bullet}}(X/B)$, by definition.

It will be convenient to have notation for compositions of morphisms $f_{\e}$,
as follows.

\begin{notn}\label{notn:path-maps} Given a path $P=(\e_1,\dots,\e_m)$ in
$G_{\I}$ (respectively, $G_{\II}$), set
$$f_{P}:=f_{\e_m} \circ \cdots \circ f_{\e_1}.$$
\end{notn}

The following definition will play an important role in the type-II
case.

\begin{defn}\label{def:simple} Given a field $K$, a $K$-valued point 
$(\sE, (V_w)_{w \in V(G_{\II})})$ of $\cG^{k,\II}_{r,d,d_{\bullet}}(X/B)$ is
\textbf{simple} if there exist $w_1,\dots,w_k \in V(G_{\II})$ (not 
necessarily distinct) and $v_i \in \sV_{w_i}$ for $i=1,\dots,k$ such that 
for all $w \in V(G_{\II})$, if $P_1,\dots,P_k$ are minimal paths from
$w_i$ to $w$ in $G_{\II}$, the vectors $f_{P_1}(v_1),\dots,f_{P_k}(v_k)$
form a basis of $\sV_w$.
\end{defn}

It almost follows from Nakayama's lemma that the simple points are an open 
substack of $\cG^{k,\II}_{r,d,d_{\bullet}}(X/B,\theta_{\bullet})$, except 
that the infinite nature of $G_{\II}$ requires some additional argument;
nonetheless, openness follows from Proposition \ref{prop:type-ii-finite} 
below. 

The foundational theorem on our constructions is the following:

\begin{thm}\label{thm:foundation} Let $X/B$ be an almost local smoothing 
family, and $k,r,d_{\bullet}$ as in Situation \ref{sit:basic}.
Then $\cG^{k,\I}_{r,d,d_{\bullet}}(X/B)$ is an Artin stack,
and the natural map
$$\cG^{k,\I}_{r,d,d_{\bullet}}(X/B) \to \cM_{r,w_0,d_{\bullet}}(X/B)$$
is relatively representable by schemes which are projective, at least
locally on the target.
Moreover, formation of $\cG^{k,\I}_{r,d,d_{\bullet}}(X/B)$
is compatible with any base change $B' \to B$ which preserves the almost 
local smoothing family hypotheses.
In particular, if $y \in B$ is a point with $X_y$ smooth, then
the base change to $y$ parametrizes pairs $(\sE,V)$ of
a vector bundle $\sE$ of rank $r$ and degree $d$ on $X_y$ together with
a $k$-dimensional vector space $V \subseteq H^0(X_y,\sE)$.

Under the further hypothesis of Situation \ref{sit:basic-ii}, all the above 
statements also hold for
$\cG^{k,\II}_{r,d,d_{\bullet}}(X/B,\theta_{\bullet})$.
Moreover, the simple locus of
$\cG^{k,\II}_{r,d,d_{\bullet}}(X/B,\theta_{\bullet})$ has universal relative dimension at 
least $k(d-k-r(g-1))$ over $\cM_{r,w_0,d_{\bullet}}(X/B)$,
and hence has universal relative dimension at least $\rho-1$ over $B$.
In particular, if some fiber 
$\cG^{k,\II}_{r,d,d_{\bullet}}(X_y/y,\theta_{\bullet})$ has dimension exactly $\rho-1$ at a 
simple point $z$, then 
$\cG^{k,\II}_{r,d,d_{\bullet}}(X/B,\theta_{\bullet})$ is universally open over $B$ in a 
neighborhood of $z$, with pure fiber dimension $\rho-1$.
\end{thm}

In the above, the notion of having universal relative dimension at least a 
given number is as introduced in Definition 7.1 of \cite{os21}. While the
same dimensional statements hold also for
$\cG^{k,\I}_{r,d,d_{\bullet}}(X/B)$, they are largely vacuous in this case,
since for the most part simple points occur only for type II linked linear
series.

In order to prove the theorem, it will be useful to introduce additional
notation.

\begin{notn}\label{notn:twists} In Situation \ref{sit:basic-ii},
given $w=(i_v)_v \in V(G_{\II})$ and a pair $(e,v)$ of an edge of 
$\Gamma$ together with an adjacent vertex, let $S \subseteq V(\Gamma)$ 
be the vertices not in the same connected component as $v$ of 
$\Gamma \smallsetminus e$. Then we write
$$t_{(e,v)}(w) = |S|b-\frac{1}{r}\sum_{v' \in S} (d_{v'}-i_{v'}).$$
\end{notn}

Then $t_{(e,v)}(w)$ is the number of times we twist down by $\Delta'_e$
on $Y_{(e,v)}$ in order to go from multidegree $w_v$ to multidegree $w$. 

Then we see easily from \eqref{eq:deg-sum}
that if $v'$ is the other vertex adjacent to $e$, we have
$$t_{(e,v)}(w)+t_{(e,v')}(w)=b.$$
We also see that if we start at $w$ and go along an edge $\e$ of
$G_{\II}$, the value of $t_{(e,v)}$ increases by $1$ whenever 
$v=v_{\II}(\e)$, decreases by $1$ whenever $v$ is adjacent to
$v_{\II}(\e)$ in $\Gamma$, and remains the same otherwise.

We now verify that in the definition of type II linked linear series, we
could have restricted our attention to a suitable finite graph.

\begin{defn}\label{def:bar-g} Consider the finite graph
$\bar{G}_{\II}$ described as follows: the vertices of $\bar{G}_{\II}$
are the subset of the vertices $w$ of $G_{\II}$ satisfying 
$t_{(e,v)}(w) \leq b$ for all $(e,v)$,
and the edges consist of all edges of $G_{\II}$ between
vertices of $\bar{G}_{\II}$. 

Finally, let $\bar{\cG}^{k,\II}_{r,d,d_{\bullet}}(X/B,\theta_{\bullet})$ be the
groupoid defined in the same manner as $\cG^{k,\II}_{r,d,d_{\bullet}}(X/B,\theta_{\bullet})$,
but with $G_{\II}$ replaced by $\bar{G}_{\II}$.
\end{defn}

\textit{A priori}, one might have to add extra edges to $\bar{G}_{\II}$
to account for paths between vertices of $\bar{G}_{\II}$ which are not
entirely contained in $\bar{G}_{\II}$. However, the following proposition
shows that this is not necessary.

\begin{prop}\label{prop:type-ii-min-paths} For any 
$w,w' \in V(\bar{G}_{\II})$, there exists a minimal path from $w$ to $w'$
containing only vertices of $\bar{G}_{\II}$.
\end{prop}

\begin{proof} Let $P=P(w,v_1,\dots,v_m)$ be a minimal path from $w$ to
$w'$ in $G_{\II}$; it suffices to show that we can reorder the $v_i$ so
that $P$ remains in $\bar{G}_{\II}$. In fact, it is enough to show that
for some $v_i$, we have $P(w,v_i)$ in $\bar{G}_{\II}$, since then if $w''$
is the endpoint of $P(w,v_i)$, we have $P(w'',v_1,\dots,v_{i-1},v_{i+1},v_m)$
a minimal path from $w''$ to $w'$, and we can repeat the process inductively.
Now, if $P(w,v_1)$ is not in $\bar{G}_{\II}$, this means necessarily that
$t_{(e,v_1)}(w)=b$ for some $e$ adjacent to $v_1$. Now, because 
$w' \in \bar{G}_{\II}$, if $v'$ is the other vertex adjacent to $e$, we see
that we must have $v'=v_i$ for some $i \neq 1$. Then again, if $P(w,v_i)$
is not in $\bar{G}_{\II}$, we must have some $t_{(e',v_i)}(w)=b$ for some
$e$ adjacent to $v_i$, and since $t_{(e,v_i)}(w)=0$, we have $e' \neq e$.
Thus, if we continue in this way, we traverse the graph $\Gamma$ without
repeating any vertices, and the process must conclude with some $v_{\ell}$
having $P(w,v_{\ell})$ in $\bar{G}_{\II}$, as desired.
\end{proof}

The following is a trivial consequence of the first part of Lemma 
\ref{lem:substack-open}.

\begin{cor}\label{cor:simple-support} 
Suppose $\sE$ on $X_y$ is a $K$-valued point of 
$\cM_{r,w_0,d_{\bullet}}(X/B)$, and we have $w \in V(G_{\II})$, 
$s \in H^0(X_y,\sE_w)$, and $v \in V(\Gamma)$ such that there is a 
path $P_v$ in $G_{\II}$ from $w$ to $w_v$ with $f_{P_v}(s) \neq 0$.
Then $s$ does not vanish uniformly on the component of $X_y$ corresponding
to $v$.

In particular, for any simple point of $\cG^{k,\II}_{r,d,d_{\bullet}}(X/B,\theta_{\bullet})$,
if $w_1,\dots,w_k$ are as in Definition \ref{def:simple}, then necessarily
each $w_i \in V(\bar{G}_{\II})$.
\end{cor}

We are now ready to prove the equivalence of the finite version of 
the type II construction.

\begin{prop}\label{prop:type-ii-finite} The canonical forgetful morphism
$$\cG^{k,\II}_{r,d,d_{\bullet}}(X/B,\theta_{\bullet}) \to 
\bar{\cG}^{k,\II}_{r,d,d_{\bullet}}(X/B,\theta_{\bullet})$$
is an equivalence of groupoids.

In particular, the simple points of 
$\cG^{k,\II}_{r,d,d_{\bullet}}(X/B,\theta_{\bullet})$ form an open substack.
\end{prop}

\begin{proof}
We visibly have a forgetful morphism 
$\cG^{k,\II}_{r,d,d_{\bullet}}(X/B,\theta_{\bullet})\to 
\bar{\cG}^{k,\II}_{r,d,d_{\bullet}}(X/B,\theta_{\bullet})$, so in order
to verify that this is an equivalence, it suffices to show that an 
$S$-valued point
$(\sE, (\bar{\sV}_{\bar{w}})_{\bar{w} \in V(\bar{G}_{\II})})$ of 
$\bar{\cG}^{k,\II}_{r,d,d_{\bullet}}(X/B,\theta_{\bullet})$ extends uniquely to a 
$S$-valued point $(\sE, (\sV_{w})_{w \in V(G_{\II})})$ of 
$\cG^{k,\II}_{r,d,d_{\bullet}}(X/B,\theta_{\bullet})$. 

For each $w \in V(G_{\II})$,
choose $\bar{w}_w \in V(\bar{G}_{\II})$ admitting the shortest possible path
$P_w$ in $G_{\II}$ from $\bar{w}_w$ to $w$. We claim that $\bar{w}_w$ is
uniquely determined by this condition, and that moreover $\bar{w}_w$ is
characterized by the property that for any minimal path $P$ from $\bar{w}_w$
to $w$, no vertex of $P$ other than $\bar{w}_w$ lies in $\bar{G}_{\II}$.
For both claims, the main observation is the following: if 
$w \not\in \bar{G}_{\II}$, and $\bar{w}_w$ and $\bar{w}'_w$ are any vertices 
of $\bar{G}_{\II}$, and $P$ and $P'$ minimal paths from $\bar{w}_w$ and 
$\bar{w}'_w$ to $w$, respectively, then we can modify $P$ and $P'$
without changing their length (or start or endpoints) so that they pass
through a common intermediate vertex $w' \neq w$. Indeed, if
$P=P(\bar{w}_w,v_1,\dots,v_m)$ and $P'=P(\bar{w}'_w,v'_1,\dots,v'_m)$,
we must have some $(e,v)$ such that $t_{(e,v)}(w)>b$, and then we see
that $v$ must occur both as one of the $v_i$ and one of the $v'_i$. 
Reordering the $v_i$ and $v'_i$ so that $v_m=v'_m=v$ then clearly has 
the desired effect. 

If now we suppose that $P$ and $P'$ both have minimal length among all
paths from vertices of $\bar{G}_{\II}$ to $w$, we conclude the first 
claim by induction on the length of $P$: if the length is $0$, there is
nothing to show, while if the length is positive, the above observation
allows us to consider the strictly shorter paths to $w'$, and the
induction hypothesis then implies that $\bar{w}_w=\bar{w}'_w$. For the
second claim, it is immediate that if $\bar{w}_w$ admits a minimal path
$P$ to $w$ among all vertices of $\bar{G}_{\II}$, then no minimal path
from $\bar{w}_w$ to $w$ can contain a vertex of $\bar{G}_{\II}$ other
than $\bar{w}_w$. Conversely, suppose that $\bar{w}_w$ has the property
that no minimal path to $w$ contains another vertex of $\bar{G}_{\II}$,
let $P$ be any minimal path to $w$, and choose 
$\bar{w}'_w \in V(\bar{G}_{\II})$ with $P'$ a path from $\bar{w}'_w$ to $w$.
We need to see that $P'$ is at least as long as $P$, and we argue once 
again by induction on the length of $P$. If $P$ has length zero, there is
nothing to show. On the other hand, if $P$ has positive length, then our
observation implies that it and $P'$ may be modified to pass through a 
common vertex $w'$. By hypothesis, we either have $w' = \bar{w}_w$, in
which case the length of $P$ is certainly at most the length of $P'$,
or $w' \not\in V(\bar{G}_{\II})$, in which case we may apply the induction
hypothesis to conclude again that the length of $P$ is at most the length
of $P'$, as desired.

Now, returning to the question of lifting points 
$(\sE, (\bar{\sV}_{\bar{w}})_{\bar{w} \in V(\bar{G}_{\II})})$ of 
$\bar{\cG}^{k,\II}_{r,d,d_{\bullet}}(X/B,\theta_{\bullet})$ to
$\cG^{k,\II}_{r,d,d_{\bullet}}(X/B,\theta_{\bullet})$,
set $\sV_w=f_{P_w} \left(\sV_{\bar{w}_w}\right)$ for each $w$. 
We first verify that this is 
in fact a subbundle, and of the correct rank $k$. Given $w$, 
write 
$P_w=P(\bar{w}_w,v_1,\dots,v_m)$. We need to show that $f_{P_w}$
is universally injective,
which may be checked over each 
point $y \in B$. Now, if we let $V' \subseteq V(\Gamma)$ be the subset 
consisting of vertices equal to $v_i$ for some $i$, the kernel of $f_{P_w}$
(over $y \in \Delta$)
consists of sections vanishing on every $Y_v$ for $v \not \in V'$. It thus
suffices to see that such a section must also vanish on every $v \in V'$.
By our second claim above, no reordering of ${P_w}$ can meet $\bar{G}_{\II}$
away from $\bar{w}_w$. It follows that for any $v \in V'$, there is
some $e \in E(\Gamma)$ adjacent to $v$ such that 
$t_{(e,v)}(\bar{w}_w) \geq b$,
and because $\bar{w}_w \in \bar{G}_{\II}$, we conclude that 
$t_{(e,v)}(\bar{w}_w)=b$. Let $v'$ be the other vertex adjacent to $e$.
Any section vanishing on $Y_{v'}$ then vanishes to order at least $b+1$
at $\Delta_e$ as a section of $\sE_{w_v}|_{Y_v}$, so by condition (I)
of Definition \ref{def:grd-space-ii},
such a section vanishes also on $v$. Thus, if $v'$ is
not in $V'$, we conclude that any section in the kernel of $f_{P_w}$ must
vanish on $Y_v$, as desired. On the other hand, if $v' \in V'$, we can
repeat the process, finding another edge $e'$ such that $t_{(e',v')}=b$.
Traversing the tree $\Gamma$ in this way, since we cannot repeat vertices
we must eventually come to a vertex not in $V'$, and we ultimately conclude
that any section in the kernel of $f_{P_w}$ must vanish on every component
we have traversed, and in particular on $Y_v$, giving us the desired
injectivity. We thus conclude that the $\sV_w$ are subbundles of rank $k$.

It is then clear that each $\sV_w$ is uniquely
determined by the linkage conditions, and it remains to check that our
construction satisfies the linkage condition. Given $w,w' \in V(G_{\II})$
connected by an edge $\e \in E(G_{\II})$ from $w$ to $w'$, we wish to verify 
that $f_{\e}(\sV_w) \subseteq \sV_{w'}$. With $\bar{w}_w$ and ${P_w}$ as above,
let $P'$ be the path from $\bar{w}_w$ to $w'$ obtained by composing ${P_w}$
with $\e$; say $P'=P(\bar{w}_w,v_1,\dots,v_m)$, where $v_m=v_{\II}(\e)$. Now, 
reorder the $v_i$ so that $P'=P'' \circ Q$, where $Q$ ends in $\bar{G}_{\II}$
and has maximum possible length (here we allow either $P''$ or $Q$ to have
length zero, as required). Then the second claim above implies that
the starting point of $P''$ is equal to $\bar{w}_{w'}$, and $P''=P_{w'}$.
By definition, we have 
$\sV_{w'}=f_{P''}(\sV_{\bar{w}_{w'}})$, and 
$f_{\e}(\sV_w)=f_{P'}(\sV_{\bar{w}_w})$, and the reordering of $P'$ doesn't
affect the latter because it only changes $f_{P'}$ by an invertible scalar. 
It thus suffices to see that 
$f_{Q}(\sV_{\bar{w}_w}) \subseteq \sV_{\bar{w}_{w'}}$, which follows from
our initial linkage hypothesis.
\end{proof}

We can now easily prove our foundational theorem on linked linear series.

\begin{proof}[Proof of Theorem \ref{thm:foundation}]
The statement being local on $B$,
we may suppose by Lemma
\ref{lem:ample-exists} that we have a divisor $D$ on $X$ which is
$\pi$-ample. We construct $\cG^{k,\I}_{r,d_{\bullet}}(X/B)$ (respectively,
$\cG^{k,\II}_{r,d_{\bullet}}(X/B,\theta_{\bullet})$) as a 
projective relative scheme over $\cM_{r,w_0,d_{\bullet}}(X/B)$, yielding
the first assertion of the theorem. Let $\cU$ be a quasicompact open 
substack of $\cM_{r,w_0,d_{\bullet}}(X/B)$, and $\widetilde{\sE}$ the
universal bundle on $X \times_B \cU$; then, replacing $D$ by a 
sufficiently high multiple, we may suppose that $R^i p_{2*} \sE_w(D)=0$ for
all $i>0$ and all $w \in V(G_{\I})$. Let $\cG$ be the prelinked
Grassmannian (see Definition \ref{def:pre-lg}) associated to $G_{\I}$ and 
$p_{2*} \sE_w(D)$, with maps $p_{2*}(f_{\e})$ for $\e \in G_{\I}$. Note
that the condition of Definition \ref{def:pre-lg} is clearly satisfied
in this case, due to Proposition \ref{prop:min-paths}.
Then, according to Corollary \ref{cor:lg-stacks}, we have that
$\cG$ is relatively representable by a projective scheme over $\cU$, and
it thus suffices to show that (the preimage of $\cU$ in)
$\cG^{k,\I}_{r,d_{\bullet}}(X/B)$ is a closed substack of $\cG$. But, making
use of Lemma \ref{lem:grd-sub} (ii) to ensure compatibility of notions of
subbundles, we see that
$\cG^{k,\I}_{r,d_{\bullet}}(X/B)$ is cut out by the closed condition that
for all $w \in V(G_{\I})$, we have $\sV_w$ in the kernel of the morphism
$$p_{2*} \sE_w(D) \to p_{2*} \left(\sE_w(D)|_D\right).$$

This proves the first statement of the theorem, and the same construction
works also for $\cG^{k,\II}_{r,d_{\bullet}}(X/B,\theta_{\bullet})$, using $\bar{G}_{\II}$ 
in place of $G_{\I}$ by virtue of Proposition \ref{prop:type-ii-finite}.

Compatibility with base change is straightforward, although we note
that the graph $\Gamma$ will typically be contracted under base change
$B' \to B$. In the type-I case, arbitrary contractions are possible, 
depending on which of the $\Delta_e$ miss the image of $B'$. For these
$e$, the corresponding morphisms $f_{\e}$ are all isomorphisms over the
image of $B'$, so we obtain a natural isomorphism 
$$\cG^{k,\I}_{r,d_{\bullet}}(X/B) \times_B B' \to
\cG^{k,\I}_{r,d'_{\bullet}}(X'/B'),$$ 
where $d'_{\bullet}$ is induced by $d_{\bullet}$.
Compatibility with base change is similar in the type-II case, except that
$\Gamma$ is either preserved or contracted to a point, so the situation
is simpler.

Finally, in the type II case, we compute the relative dimension of the 
simple locus. The bundles $p_{2*} \sE_w(D)$ all have rank 
$d+r \deg D+r(1-g)$, so Corollary \ref{cor:lg-stacks} implies that
the simple locus of $\cG$ is smooth of relative dimension
$k(d+r \deg D+r(1-g)-k)$ over $\cU$, which is in turn smooth of dimension
$r^2(g-1)$. We next claim that, at least etale locally on the base, the 
substack $\cG^{k,\II}_{r,d_{\bullet}}(X/B,\theta_{\bullet})$ is cut out by 
$k r \deg D$ equations. Etale locally, we may suppose that 
$D=\sum_{v \in V(\Gamma)} D_v$, where the $D_v$ are disjoint divisors
supported in the smooth locus of $\pi$, and with $D_v|_{X_y}$ contained 
in the component of $X_y$ corresponding to $\cl_y v$ for every $y \in B$.
We will show that 
$\cG^{k,\II}_{r,d_{\bullet}}(X/B,\theta_{\bullet})$ may be cut out by the conditions that
for all $v \in V(\Gamma)$, we have $\sV_{w_v}$ in the kernel of the morphism
$$p_{2*} \sE_{w_v}(D) \to p_{2*} \left(\sE_{w_v}(D)|_{D_v}\right).$$
Each $v$ will then impose $k r \deg D_v$ equations, so the claim will follow.
Since the $D_v$ are disjoint, it suffices to prove that the above condition
implies that for all $w \in V(\bar{G}_{\II})$ and all $v \in V(\Gamma)$,
we have $\sV_{w}$ in the kernel of the morphism
$$p_{2*} \sE_{w}(D) \to p_{2*} \left(\sE_{w}(D)|_{D_v}\right).$$
We observe that we can find a path $P=(\e_1,\dots,\e_n)$ from $w$ to $w_v$ in 
$G_{\II}$ such that $v_{\II}(\e_i) \neq v$ for each $i$. Indeed,
we can construct such a path by working inductively outwards from $v$
in $\Gamma$: for each edge $e$ adjacent to $v$, if $t_{(e,v)}(w)>0$,
and $v'$ is the other vertex adjacent to $e$, we add edges $\e_i$
with $v_{\II}(\e_i)=v'$ until the new $t_{(e,v)}(w')=0$. For
each such $v'$, we then have $t_{(e,v')}(w')=b$, but if
$t_{(e',v')}(w')>0$ for some $e' \neq e$, we add edges $\e_i$ with
$v_{\II}(\e_i)=v''$, where $v''$ is the other vertex
adjacent to $e'$, and continuing in this way we will finally reach $w_v$.
Thus, if we take the morphism 
$p_{2*} \sE_{w}(D) \to p_{2*} \sE_{w_v}(D)$ induced by a minimal path
from $w$ to $w_v$, because 
$D_v$ is disjoint from $Y_{v'}$ for any $v' \neq v$
we see that the induced morphism 
$p_{2*} \left(\sE_{w}(D)|_{D_v}\right) \to 
p_{2*} \left(\sE_{w_v}(D)|_{D_v}\right)$
is an isomorphism.
Now, the composed map $\sV_w \to p_{2*} \left(\sE_{w_v}(D)|_{D_v}\right)$
factors through $\sV_{w_v} \to p_{2*} \left(\sE_{w_v}(D)|_{D_v}\right)$,
which vanishes by hypothesis, so we conclude that
$\sV_w \to p_{2*} \left(\sE_{w}(D)|_{D_v}\right)$ is likewise zero, as 
desired.
This proves the claim, and by Corollary 7.7 of \cite{os21},
we conclude that 
the simple locus of $\cG^{k,\II}_{r,d_{\bullet}}(X/B,\theta_{\bullet})$ has universal 
relative dimension at least $k(d+r(1-g)-k)$ over 
$\cM_{r,w_0,d_{\bullet}}(X/B)$, yielding the desired statement.
\end{proof}

\begin{rem} In the classical rank-$1$ case, it should be possible to
restrict to a simpler (and smaller) bounded collection of multidegrees than 
the one used for $\bar{G}_{\II}$. Indeed, because no line bundle of negative 
degree can have a global section on a smooth curve, it is more natural to
restrict to
$w=(i_v)_{v \in V(\Gamma)} \in V(G_{\II})$ such that $i_v \geq 0$ for all $v$.
However, in this case Proposition \ref{prop:type-ii-min-paths} no longer 
holds, 
so in fact the definitions and arguments involve additional complications.
\end{rem}

We conclude with a brief study of simplicity, showing that it is enough
to check the basis condition at the vertices $w_v$.

\begin{lem}\label{lem:simple-crit} 
Given a field $K$ and a $K$-valued point $(\sE, (V_w)_{w \in V(G_{\II})})$
of the stack $\cG^{k,\II}_{r,d,d_{\bullet}}(X/B,\theta_{\bullet})$, suppose 
there exists 
$w_1,\dots,w_k \in V(G_{\II})$ (not necessarily distinct) and 
$v_i \in V_{w_i}$ for $i=1,\dots,k$ such that for all 
$v \in V(\Gamma)$, there is a path $P_{i,v}$ such that 
$f_{P_{1,v}}(v_1),\dots,f_{P_{k,v}}(v_k)$ form a 
basis of $V_{w_v}$. Then $(\sE, (V_w)_{w \in V(G_{\II})})$ is simple.
\end{lem}

\begin{proof} 
Let $\Delta\subseteq B$ be the image of the singular locus of $\pi$. If 
the point in question lies outside of $\Delta$, then all points are simple,
and there is nothing to show. On the other hand, according to the hypothesis 
of Situation \ref{sit:basic-ii}, all fibers of $\pi$ over $\Delta$ have dual 
graphs equal to $\Gamma$, so we might as well assume that $B$ is a point.

Now, we observe that if we have a path $P$ from some $w_v$ to $w'$ such
that $f_{P}$ is injective (equivalently, non-zero) on $Y_v$, then
$f_P(V_{w_v})=V_{w'}$. Indeed, this is immediate from the fact that the 
restriction map to $Y_v$ is injective on $V_{w_v}$, by Lemma 
\ref{lem:substack-open}.

To prove the lemma, it suffices to prove that for any $w \in V(G_{\II})$, if
$P_1,\dots,P_k$ are minimal paths from $w_1,\dots,w_k$ respectively to $w$,
then the $f_{P_i}(v_i)$ are linearly independent in $V_w$. Accordingly, 
suppose that we have $\sum_i a_i f_{P_i}(v_i)=0$. For a given $i_0$, if we
write $P_{i_0}=P(w_{i_0},u_1,\dots,u_m)$ with $u_i \in V(\Gamma)$, by
minimality of $P_{i_0}$ and Proposition \ref{prop:min-paths}
there some $v\in V(\Gamma)$ not appearing among the $u_i$,
and then $f_{P_{i_0}}$ is injective on $Y_v$.
If $P$ is a minimal path from $w$ to $w_v$, then $f_{P}$ may vanish on 
$Y_v$, and similarly if $P$ is a minimal path from $w_v$ to $w$. However,
there always exists some $w'$ such that, if $P$ and $P'$
are minimal paths from $w$ to $w'$ and from $w_v$ to $w'$ respectively,
then $f_P$ and $f_{P'}$ are both nonvanishing on $Y_v$. By the above
observation, we have that $f_{P'}(V_{w_v})=V_{w'}$. If $P'_1,\dots,P'_k$
are minimal paths from $w_1,\dots,w_k$ respectively to $w_v$, then since
the $f_{P'_i}(v_i)$ are assumed linearly independent, we find that
the $f_{P' \circ P'_i}(v_i)$ are likewise linearly independent. On the
other hand, by hypothesis $\sum_i a_i f_{P \circ P_i}(v_i)=0$.
Now, for each $i$ we have that $f_{P' \circ P'_i}(v_i)$ and 
$f_{P \circ P_i}(v_i)$ differ by a scalar, so we conclude that for each
$i$, either $a_i = 0$ or $f_{P \circ P_i}(v_i)=0$. Since 
$f_{P \circ P_{i_0}}$ was by construction injective on $Y_v$, and each
$v_i$ is non-zero on every component of $X$ by Corollary 
\ref{cor:simple-support}, we have that 
$f_{P \circ P_{i_0}} (v_{i_0}) \neq 0$, and
hence $a_{i_0}=0$. Since $i_0$ was arbitrary, we obtain the desired
linear independence.
\end{proof}

\begin{rem}\label{rem:two-comps} We observe that our two constructions
are in fact completely equivalent in the two-component case, so that they
each provide a generalization of \cite{os8}. Indeed, we have already
observed that the additional condition of Situation \ref{sit:basic-ii}
is tautologically satisfied in the two-component case. It is clear that
in this case, we have $G_{\I}$ and $\bar{G}_{\II}$ equal, each a chain 
with edges in both directions. 
Proposition \ref{prop:two-types-compat} then implies that the bundles are
isomorphic in both constructions, and the maps agree up to invertible
scalar, so the resulting constructions agree.

This is consistent with the nature of the additional choice of 
$\theta_{\bullet}$ in the second construction, since a given $\theta_v$
is determined up to scalar unless $Y_v$ is disconnecting, which never
occurs in the two-component case.
\end{rem}

\section{Stack structures on the Eisenbud-Harris-Teixidor construction}
\label{sec:eht}

We now examine in some detail the generalization by Teixidor i Bigas of the 
construction of Eisenbud and Harris. We present a new point of view on
this construction, which will be important for two reasons: it is an
important step in our comparison theorems, and it yields a canonical
stack structure on the space of Eisenbud-Harris-Teixidor limit linear
series which can be placed into smoothing families. The latter is new
even for the classical case of rank $1$. After introducing this stack
structure, we compare it to the previously constructed stacks of
linked linear series.

In this section, we consider the case that $B=\Spec F$ is a point,
and we therefore omit $B$ from the notation. Thus, we have a 
single reducible curve $X$ over $\Spec F$, with dual graph $\Gamma$. 
In particular, the notation $\pi_*$ simply corresponds to taking global
sections, but we will frequently use $\pi_*$, especially for constructions
which generalize to families.
As before, for $v \in V(\Gamma)$, we will denote by $Y_v$ the corresponding 
component of $X$. Because much of this section
is rather technical, our approach is, as much as possible, to prove 
general results by carrying out key calculations over fields, and then
reducing the general statements to this situation.

\subsection{Statements on points}
We begin by recalling the definition of Eisenbud-Harris-Teixidor limit
linear series, given as a subset of a certain stack which we now describe.
Working on the level of points, we then give an alternate characterization 
of this subset which will lend itself to introducing a stack structure.

Given $v \in V(\Gamma)$, we have the stack
$\cG^k_{r,d_v}(Y_v)$. Given $e \in E(\Gamma)$ adjacent to $v$,
we have also the stack $\cM_{r}(\Delta'_e)$ of rank-$r$ vector
bundles on $\Delta'_e$, where $\Delta'_e$ is the node corresponding to 
$e$. There is a natural restriction map
$\cG^k_{r,d_v}(Y_v) \to \cM_r(\Delta'_e)$, but for reasons which will
soon become apparent, we instead consider the map induced first by
twisting the universal bundle by $\sO_{w_v,w_0}|_{Y_v}$,
and then restricting to $\Delta'_e$. Here $w_0$ is as in Definition
\ref{def:grd-space-i}.

\begin{notn}\label{not:pkrd}
Denote by $\cP^k_{r,d_{\bullet}}(X)$ the product of all the stacks
$\cG^k_{r,d_v}(Y_v)$ fibered over the stacks
$\cM_r(\Delta'_e)$ via the above maps.
\end{notn}

Note that while $\cM_r(\Delta'_e)$ has a single point, that point has
automorphism group $\GL_r$, so the stack is non-trivial. The purpose
of fibering over it is that, due to the definition of $2$-fibered products,
this process precisely introduces a choice of gluing map at each node.
In fact, the moduli stack $\cM_{r,w_0}(X)$ can be realized as the 
$2$-fibered products of the stacks $\cM_{r,d_v}(Y_v)$ over the stacks
$\cM_r(\Delta'_e)$ (twisting to common multidegree $w_0$ as before), so
$\cP^k_{r,d_{\bullet}}(X)$ has a natural forgetful morphism to
$\cM_{r,w_0}(X)$.

We now recall the definition of Eisenbud-Harris-Teixidor limit linear
series. Following Teixidor's description as closely as possible, we give a 
set-theoretic description of the conditions inside $\cP^k_{r,d_{\bullet}}(X)$,
but we will endow this set with a locally closed substack structure in Lemma
\ref{lem:teixidor-translate} and \S \ref{sec:eht-stack} below. For 
convenience in gluing map notation, we choose directions for all
edges of $\Gamma$. We will also make use of the twisting line bundles
we constructed in Notation \ref{not:twisting-bundles} and
\ref{not:twisting-bundles-ii}.

\begin{defn}\label{def:eht-lls}
Let $((\sE^v,V^v)_{v \in V(\Gamma)},(\vp_e)_{e \in E(\Gamma)})$ be a 
$K$-valued point of 
$\cP^k_{r,d_{\bullet}}(X)$, where $(\sE^v,V^v)$ is the corresponding
point of $\cG^k_{r,d_v}(Y_v)$, and if $e$ is an edge from $v$ to $v'$,
$$\vp_e:
\left.\left(\sE^v \otimes \left(\sO_{w_v,w_0}|_{Y_v}\right)
\right)\right|_{\Delta'_e} 
\risom
\left.\left(\sE^{v'} \otimes \left(\sO_{w_{v'},w_0}|_{Y_{v'}}\right)
\right)\right|_{\Delta'_e}$$
is the corresponding gluing isomorphism. Then $((\sE^v,V^v)_v,(\vp_e)_e)$
is an \textbf{Eisenbud-Harris-Teixidor} (or \textbf{EHT}) \textbf{limit 
linear series}
if for each $e \in E(\Gamma)$ from $v$ to $v'$, we have:
\begin{Ilist}
\itm
$H^0(Y_v,\sE^v(-(b+1)\Delta'_e))=0$ and similarly for $v'$;
\itm
if $a^{e,v}_1,\dots,a^{e,v}_k$ and $a^{e,v'}_1,\dots,a^{e,v'}_k$
are the vanishing sequences at $\Delta'_e$ for $(\sE^v,V^v)$ and 
$(\sE^{v'},V^{v'})$ respectively, then for every $i$ we have
\begin{equation}\label{eq:eht-compat} a^{e,v}_i+a^{e,v'}_{k+1-i} \geq b;
\end{equation}
\itm
there exist bases $s^{e,v}_1,\dots,s^{e,v}_k$ and
$s^{e,v'}_1,\dots,s^{e,v'}_k$ of $V^v$ and $V^{v'}$ respectively such
that $s^{e,v}_i$ has vanishing order $a^{e,v}_i$ at $\Delta'_e$ for each $i$,
and similarly for $s^{e,v'}_i$, and if we have $a^{e,v}_i+a^{e,v'}_{k+1-i}=b$
for some $i$, then 
\begin{equation}\label{eq:eht-gluing}\vp_e(s^{e,v}_i)=s^{e,v'}_{k+1-i};
\end{equation}
\end{Ilist}

We say that $((\sE^v,V^v)_v,(\vp_e)_e)$ is \textbf{refined} if equality
always holds in \eqref{eq:eht-compat}.
\end{defn}

\begin{notn}\label{not:eht}
We denote by $\cG^{k,\EHT}_{r,d,d_{\bullet}}(X)$ the set of
Eisenbud-Harris-Teixidor limit linear series, and by 
$\cG^{k,\EHT,\refn}_{r,d,d_{\bullet}}(X)$ the refined locus.
\end{notn}

\begin{rem}\label{rem:eht-gluing} Condition (III) of Definition 
\ref{def:eht-lls} requires some explanation. As we have already 
remarked, the gluings $\vp_e$ together with the $\sE^v$ yield a vector
bundle $\sE$ of multidegree $w_0$. As in the linked linear series 
constructions, our twisting bundles then yield
vector bundles $\sE_w$ of multidegree $w$ for any $w \in V(G_{\II})$,
and this in turn induces gluing maps
$$\left.\left(\sE^v \otimes \left(\sO_{w_v,w}|_{Y_v}\right)
\right)\right|_{\Delta'_e} 
\risom
\left.\left(\sE^{v'} \otimes \left(\sO_{w_{v'},w}|_{Y_{v'}}\right)
\right)\right|_{\Delta'_e}$$
for each $w$. When $w \in V(G_{\I})$, we have
$\sE^v \otimes \left(\sO_{w_v,w}|_{Y_v}\right)$ and
$\sE^{v'} \otimes \left(\sO_{w_{v'},w}|_{Y_{v'}}\right)$ described
explicitly in terms of twisting down by a multiple of the node $\Delta'_e$.

Now, when $a^{e,v}_i+a^{e,v'}_{k+1-i}=b$ we have an induced multidegree 
$w \in V(G_{\I})$, and considering 
$s^{e,v}_i$ as a section of $\sE^v(-a^{e,v}_i \Delta'_e)$ and
$s^{e,v'}_{k+1-i}$ as a section of $\sE^{v'}(-a^{e,v'}_{k+1-i} \Delta'_e)$,
we can apply the above gluing map to compare their values at $\Delta'_e$.

Note that this formulation of the gluing condition differs slightly from
that of \cite{te1}, replacing a projectivized gluing condition with
a gluing condition on (twists of) the vector bundles themselves. However,
our formulation is consistent with the way in which dimension counts are 
carried out in \cite{te1} and related work of Teixidor i Bigas.
\end{rem}

The conditions of Definition \ref{def:eht-lls} are invariant under
field extension, so yield a well-defined subset of
$\cP^k_{r,d_{\bullet}}(X)$. However, obtaining a stack structure is
subtler; see Remark \ref{rem:eht-no-stack}.

We observe that the refined condition is in fact open.

\begin{prop}\label{prop:refined-open} We have
$\cG^{k,\EHT,\refn}_{r,d,d_{\bullet}}(X)$ open in
$\cG^{k,\EHT}_{r,d,d_{\bullet}}(X)$.
\end{prop}

\begin{proof} Given a point of 
$\cG^{k,\EHT,\refn}_{r,d,d_{\bullet}}(X)$, let $\{a^{e,v}\}_{(e,v)}$ be
the vanishing sequences at the nodes. Then because of condition (II)
of Definition \ref{def:eht-lls}, we can find a refined open neighborhood
of the given point in $\cG^{k,\EHT}_{r,d,d_{\bullet}}(X)$ by imposing
that each vanishing sequence be at most equal to the given
$\{a^{e,v}\}_{(e,v)}$, which is an open condition.
\end{proof} 

Given $w \in V(G_{\II})$ and $\theta_{\bullet}$ as in Notation 
\ref{not:type-ii-maps}, if we fix a minimal path $P_v$ from $w$ to
$w_v$ for each $v \in \Gamma$, then the maps $f_{P_v}$ together with
restriction to $Y_v$ and modding out by $V^v$ induce a map
\begin{equation}\label{eq:gluing-ker}
\pi_* \sE_w \to \bigoplus_{v \in V(\Gamma)} (\pi_* \sE^v)/V^v.
\end{equation}
The kernel of this map can be thought of as sections of $\pi_* \sE_w$
gluing from sections lying in the spaces $V^v$ on each $Y_v$.
If $w \in V(G_{\I})$, then according to Proposition 
\ref{prop:two-types-compat} the maps $f_{P_v}$ and hence 
\eqref{eq:gluing-ker} will be independent of the choice of 
$\theta_{\bullet}$.

In order to introduce our stack structures, we examine the relationship 
between the Eisenbud-Harris-Teixidor conditions, and the space of sections 
of the $V^v$ available to glue in arbitrary multidegrees.
The key result is the following, which says the vanishing and gluing 
conditions for an EHT limit series may be reinterpreted in terms of
existence, for each $w$, of $k$-dimensional subspaces of the kernel
of \eqref{eq:gluing-ker}.

\begin{lem}\label{lem:teixidor-translate} For a $K$-valued point 
$((\sE^v,V^v)_{v \in V(\Gamma)},(\vp_e)_{e \in E(\Gamma)})$
of $\cP^k_{r,d_{\bullet}}(X)$ with the $\sE^v$ satisfying (I) of 
Definition \ref{def:eht-lls}, the following are equivalent.
\begin{alist}
\itm 
$((\sE^v,V^v)_{v \in V(\Gamma)},(\vp_e)_{e \in E(\Gamma)})$
lies in $\cG^{k,\EHT}_{r,d,d_{\bullet}}(X)$.
\itm for every $w \in V(G_{\II})$, the map \eqref{eq:gluing-ker}
has kernel of dimension at least $k$.
\itm for every $w \in V(G_{\I})$, the map \eqref{eq:gluing-ker}
has kernel of dimension at least $k$.
\end{alist}

Furthermore, under these equivalent conditions, the kernel of
\eqref{eq:gluing-ker} always has dimension exactly $k$ if $w=w_v$ for some
$v \in V(\Gamma)$, and if $((\sE^v,V^v)_v,(\vp_e)_e)$ is a refined EHT
limit series, then \eqref{eq:gluing-ker} has kernel of dimension
exactly $k$ for all $w \in V(G_{\I})$.
\end{lem}

\begin{rem} The converse of the last statement of Lemma 
\ref{lem:teixidor-translate} is not true even in the case of rank $1$
with a two-component curve. That is, non-refined limit linear series
may have kernel dimension exactly $k$ in \eqref{eq:gluing-ker} for all
$w \in V(G_{\I})$; see Lemma 5.9 of \cite{os13}.
\end{rem}

Before proving the lemma, we introduce some notation, giving a divisor 
variant of Notation \ref{not:twisting-bundles}.

\begin{notn}\label{not:twisting-divisors} 
Given $w \in V(\bar{G}_{\II})$ and $v \in V(\Gamma)$, let
$$D_{w,v}= \sum_{e\text{ adjacent to }v} t_{(e,v)}(w)
\Delta'_e.$$
\end{notn}

The significance of $D_{w,v}$ is that it describes the amount of vanishing
at the nodes required of a section of of $V^v$ in order for it to be possible
that it is the restriction to $Y_v$ of a section of $\pi_* \sE_w$. Put
differently, it describes how much vanishing we have to impose on sections
of $V^v$ in order to be able to use them to construct sections of a space
$V_w \subseteq \pi_* \sE_w$.
By definition, we see that because $w \in V(\bar{G}_{\II})$, we have
$t_{(e,v)}(w) \geq 0$ for all edges
$e$ adjacent to a given $v$, so the divisor $D_{w,v}$ is effective.

\begin{notn}\label{not:twisted-subspaces} Given 
$((\sE^v,V^v)_v,(\vp_e)_e) \in \cG^{k,\EHT,}_{r,d,d_{\bullet}}(X)$, and
$w \in V(\bar{G}_{\II})$, for given $v \in V(\Gamma)$ denote by $\sE^v(w)$ 
the sheaf $\sE^v(-D_{w,v})$, and by $V^{v}(w)$ the subspace of $V^{v}$ 
consisting of sections lying in $\sE^v(w)$.
\end{notn}

It will be useful to consider a component-wise variant of the map
\eqref{eq:gluing-ker}. Namely, using our choice of directions for each edge
of $\Gamma$, given 
$((\sE^v,V^v)_{v \in V(\Gamma)},(\vp_e)_{e \in E(\Gamma)})$ as in
Lemma \ref{lem:teixidor-translate}, for any $w \in V(\bar{G}_{\II})$ we 
obtain a map
\begin{equation}\label{eq:gluing-ker2}
\bigoplus_{v \in V(\Gamma)} V^v(w) \to \bigoplus_{e \in E(\Gamma)} 
\left.\left(\sE^{h(e)}(w)\right)\right|_{\Delta'_e},
\end{equation}
given by the natural restriction map on each $h(e)$, and by the negative of 
the restriction map induced by $\vp_e$ on each $t(e)$.
We will be interested in the kernel of this map, which is independent of
the choice of directions on $E(\Gamma)$.

Our first observation is the following:

\begin{lem}\label{lem:kers-same} The kernels of \eqref{eq:gluing-ker}
and \eqref{eq:gluing-ker2} are identified via the maps
$\pi_* \sE_w \to \pi_* \sE^v$ for $v \in V(\Gamma)$.
\end{lem}

\begin{proof} $\pi_* \sE_w$ maps isomorphically onto the kernel of
$$\bigoplus_{v \in V(\Gamma)} 
\pi_* \left(\sE^v(w)\right)
\to \bigoplus_{e \in E(\Gamma)} 
\left.\left(\sE^{h(e)}(w)\right)\right|_{\Delta'_e}$$
in general, since this is simply the restriction exact sequence.
By construction, we then have that the kernel of \eqref{eq:gluing-ker}
is precisely the subspace of $\pi_* \sE_w$ which maps into the spaces
$V^v(w)$ for each $v$, and is thus identified with the kernel of
\eqref{eq:gluing-ker2}, as claimed. 
\end{proof}

We next give a description of what happens under restriction to a subcurve
$X'$, when non-minimal parts of the multidegree are supported on $X'$.

\begin{lem}\label{lem:restrict-kers} Let $X'$ be a connected subcurve of
$X$, and denote by $\Gamma' \subseteq \Gamma$ the associated dual graph.
Let $((\sE^v,V^v)_{v \in V(\Gamma)},(\vp_e)_{e \in E(\Gamma)})$
be a $K$-valued point of $\cP^k_{r,d_{\bullet}}(X)$, fix 
$w \in V(\bar{G}_{\II})$ which is 
minimal (i.e., equal to $d_v-rb$) on all $v$ not in $\Gamma'$. Finally,
let $\vp$ be the restriction map from the kernel of \eqref{eq:gluing-ker2}
to the kernel of the corresponding map with $\Gamma'$ in place of
$\Gamma$. Then:
\begin{ilist}
\itm $\vp$ is injective if the $\sE^v$ satisfy (I) of Definition 
\ref{def:eht-lls};
\itm $\vp$ is an isomorphism if $((\sE^v,V^v)_{v},(\vp_e)_{e})$
constitutes an EHT limit series.
\end{ilist}
\end{lem}

\begin{proof} The first observation is that by the hypothesis on the form of
$w$, for any $v \not\in\Gamma'$, then $D_{w,v}=b\Delta'$, where $\Delta'$ is 
the unique node of $X$ joining $Y_v$ to a component of $X$ strictly closer 
to $X'$. We prove both statements of the lemma inductively, noting that
$\vp$ is obtained as a composition of restriction maps which drop one
component at a time. It thus suffices to consider the case that $X'$ is 
obtained from $X$ by dropping a single component $Y_v$, meeting $X'$ at
a single node $\Delta'$.

First, suppose condition (I) of Definition 
\ref{def:eht-lls} is satisfied. Then the kernel of $\vp$ consists
of tuples of sections of $V^{v'}(w)$ for $v' \in \Gamma$ which glue at all 
nodes, and vanish uniformly away from $Y_v$. But these conditions imply that 
the section over $Y_v$ must vanish at $\Delta'$ as a section of 
$\pi_* \sE^v(w)=\pi_* \sE^v(-b\Delta')$, so is a section of 
$\pi_* \sE^v(-(b+1)\Delta')$, and must be equal to zero by hypothesis.
This proves the desired injectivity assertion.

For surjectivity, suppose also that conditions (II) and (III) of 
Definition \ref{def:eht-lls} are satisfied, and we are given a tuple
of sections of $V^{v'}(w)$ for $v' \in \Gamma'$ which glue at nodes on $X'$.
Let $v''$ be the vertex adjacent to $v$ in $\Gamma$, so that $\Delta'$ is
the node joining $Y_v$ to $Y_{v''}$.
We wish to show that there is a section of $V^v(w)$ which glues to the
given section $s \in V^{v''}(w)$. If $s$ vanishes at $\Delta'$, there
is nothing to show, since we can choose the zero section of $V^v(w)$. On
the other hand, if $s$ is nonvanishing at $\Delta'$, we claim that
conditions (II) and (III) imply that there is a section of $V^v(w)$ 
gluing to it. Let $s^1_j$ and $s^2_j$ be bases for
$V^v$ and $V^{v''}$ respectively as in condition (III), and suppose 
$s^2_1,\dots,s^2_{\ell}$ are the sections nonvanishing at $\Delta'$; thus,
the restriction of $s$ to $\Delta'$ agrees with some linear combination
of the restrictions to $\Delta'$ of the $s^2_j$ with $j \leq \ell$. Now,
by condition (II) we must have that $s^1_{k+1-\ell},\dots,s^1_k$ all
vanish to order at least $b$ at $\Delta'$, and by condition (I) they
must vanish to order exactly $b$. Thus by condition (III), all the
$s^2_j$ with $\ell \leq j$ glue to sections of $V^v(w)$, and we conclude
that $s$ likewise glues to give a section of $V^v(w)$, as desired.
\end{proof}

\begin{proof}[Proof of Lemma \ref{lem:teixidor-translate}]
Our first observation is that we may replace condition (b) by the same
condition with $\bar{G}_{\II}$ in place of $G_{\II}$. Indeed, if 
condition (b) is satisfied for all $w \in V(\bar{G}_{\II})$, and we
have $w \in V(G_{\II})$, then, following the notation of the proof of 
Proposition \ref{prop:type-ii-finite}, let $P_w$ be a minimal length
path from $\bar{G}_{\II}$ to $w$, with starting point $\bar{w}_w$. Then
we have that $f_{P_w}$ is universally injective, but it clearly maps the
kernel of \eqref{eq:gluing-ker} in multidegree $\bar{w}_w$ to the kernel
in multidegree $w$, so we conclude that the kernel in multidegree $w$ also
has dimension at least $k$, as desired.

We now prove the stated equivalence by induction on the number of components 
of $X$, using Lemma \ref{lem:kers-same} to replace the kernel of
\eqref{eq:gluing-ker} with that of \eqref{eq:gluing-ker2}. The base 
case for the induction is the case of two components.
Since in this case $\Gamma$
consists of two vertices and a single edge, to simplify notation we denote
the unique node by $\Delta'$, and
by $d_i$, $V^i$ and $a^i_{\ell}$ (with ${\ell}=1,\dots,k$, and $i=1,2$) 
respectively 
the degrees, spaces of sections, and vanishing sequences at $\Delta'$ of the 
$i$th component. Now, in this case conditions (b) and (c) are trivially 
equivalent. We then have $w=(d_1-rj, d_2-r(b-j))$ for some $j$, and 
$$\dim V^1(w)=\#\{\ell:a_\ell^1 \geq j\}, \quad \text{and} \quad
\dim V^2(w)=\#\{\ell:a_\ell^2 \geq b-j\}.$$
Let $s$ be the rank of \eqref{eq:gluing-ker2}. Then the dimension of
the kernel is equal to $\dim V^1(w)+\dim V^2(w) -s$.
Say $a^1_{\ell_1},\dots,a^1_{\ell_2}$ are the vanishing indices equal to
$j$, and $a^2_{\ell_3},\dots,a^2_{\ell_4}$ are the vanishing indices equal
to $b-j$. Here, if no $a^1_{\ell}$ is equal to $j$, we set $\ell_2=\ell_1-1$
to be maximal with $a^1_{\ell_2} \leq j$, and similarly if no $a^2_{\ell}$ 
is equal to $b-j$. Then we have $\dim V^1(w)=k+1-\ell_1$, and 
$\dim V^2(w)=k+1-\ell_3$, so the condition that the kernel of 
\eqref{eq:gluing-ker2} have dimension at least $k$ is equivalent to the
inequality 
\begin{equation}\label{eq:gluing-ineq}
k+2 \geq \ell_1+\ell_3+s.
\end{equation}

Now, suppose (a) is satisfied.
We first observe that condition (II) of an EHT limit series implies that
$\ell_1+\ell_4 \leq k+1$, and $\ell_2+\ell_3 \leq k+1$.  Indeed, if
$\ell_1>1$, then $a^1_{\ell_1-1}<j$, so we have $a^2_{k+2-\ell_1}>b-j$,
and hence $k+2-\ell_1>\ell_4$. On the other hand, if $\ell_1=1$, then
$\ell_4 \leq k$ immediately gives the first inequality. The second
inequality follows by considering $\ell_3$ in place of $\ell_1$.
To prove \eqref{eq:gluing-ineq},
first consider the case that $k+1 > \ell_2+\ell_4$. We necessarily have
$$s \leq (\ell_2+1-\ell_1)+(\ell_4+1-\ell_3)<k+3-\ell_1-\ell_3,$$
giving the desired inequality. On the other hand, if $k+1 \leq \ell_2+\ell_4$,
we obtain a range of indices 
$a^2_{k+1-\ell_2},\dots,a^2_{\ell_4}$ complementary to
$a^1_{k+1-\ell_4},\dots,a^1_{\ell_2}$, and the gluing condition then gives us
$$s \leq (\ell_2+1-\ell_1)+(\ell_4+1-\ell_3)-(\ell_4+\ell_2-k)
= k+2-\ell_1-\ell_3,$$
as desired. We thus conclude that (a) implies (b).

Observe also that if (a) is satisfed and we have the refinedness condition,
then following the above argument, we have have
$a^2_{k+1-\ell_1}=b-j$, so $k+1-\ell_1 \leq \ell_4$, so we conclude
that $\ell_1+\ell_4=k+1$, and similarly $\ell_2+\ell_3=k+1$. Thus
$$\ell_2+1-\ell_1=\ell_4+1-\ell_3=\ell_4+\ell_2-k=k+2-\ell_1-\ell_3,$$
and since $s \geq \#\{\ell:a^1_{\ell} = j\}=\ell_2+1-\ell_1$, we conclude 
that \eqref{eq:gluing-ineq} is satisfied with equality, and
hence the kernel of \eqref{eq:gluing-ker2} has dimension exactly $k$.

Conversely, suppose that (b) is satisfied, and choose any $j$ occurring as
$a^1_\ell$ for some $\ell$; setting $\ell_1,\dots,\ell_4$ and $s$ as above,
by \eqref{eq:gluing-ineq} we have $s \leq k+2-\ell_1-\ell_3$. Given
$\ell$ with $\ell_1 \leq \ell \leq \ell_2$, we first wish to see that
$a^1_{\ell}+a^2_{k+1-\ell} \geq b$.
Since $s \geq \ell_2+1-\ell_1$, we obtain $k+1-\ell_2 \geq \ell_3$.
We conclude that
$$a^1_{\ell}+a^2_{k+1-\ell}\geq a^1_{\ell}+a^2_{k+1-\ell_2} \geq
a^1_{\ell}+a^2_{\ell_3}=b,$$
so condition (II) of an EHT limit series is satisfied, as desired.
We next wish to verify the gluing condition. Again applying
\eqref{eq:gluing-ineq}, the
images of $V^1(w)$ and $V^2(w)$ at the node $\Delta'$ must intersect in
dimension at least 
$$(\ell_2+1-\ell_1)+(\ell_4+1-\ell_3)-(k+2-\ell_1-\ell_3)
= \ell_2+\ell_4-k.$$
If this is non-positive, there are no indices $\ell$ between $\ell_1$ and
$\ell_2$ such that $k+1-\ell$ is between $\ell_3$ and $\ell_4$, and we
claim that if it is positive, it is equal to the number of such indices.
Indeed, we already know that condition (II) of an EHT limit series is
satisfied, so by the above observation, we have $\ell_1+\ell_4 \leq k+1$
and $\ell_2+\ell_3 \leq k+1$, and the claim follows easily. We
conclude that an appropriate choice of basis elements will also 
satisfy the gluing condition, and hence that we have an EHT limit series.

For the induction step, assume we have at least three components.
To see that (a) implies (b), the
basic observation is that the restriction of an EHT limit linear series
to a connected subcurve is still an EHT limit linear series. 
Choose a component $Y_{v_1}$ meeting the rest of $X$ at only a single node
$\Delta'_{e_1}$. Let $X'$ be the complement of $Y_{v_1}$ in $X$, and let 
$Y$ be the union of $Y_{v_1}$ and $Y_{v_2}$, where $v_2$ is the other edge 
adjacent to $e_1$. Denote by $\Gamma'$ the dual graph of $X'$. 
Restricting to $X'$ and $Y$, we have multidegrees $w'$ and $w''$ which
differ from $w$ only in index $v_2$: $w'-w$ in index $v_2$ is
$r t_{(e_1,v_2)}(w)$, while $w''-w$ in index $v_2$ is 
$r \sum_{e'} t_{(e',v_2)}(w)$, where $e'$ ranges over edges other than
$e_2$ which are adjacent to $v_2$.
Then, by the induction hypothesis we have that the kernels of the maps
$$
\bigoplus_{v \in V(\Gamma')} V_v(w') \to \bigoplus_{e \in E(\Gamma')} 
\left(\sE^{h(e)} \otimes (\sO_{w_v,w'})|_{Y_{h(e)}}\right)|_{\Delta'_e}
$$
and
$$
V_{v_1}(w'')\oplus V_{v_2}(w'') \to 
\left(\sE^{h(e_1)} \otimes 
(\sO_{w_{v_1},w''})|_{Y_{v_1}}\right)|_{\Delta'_{e_1}}
$$
each have dimension at least $k$. However, the kernel of 
\eqref{eq:gluing-ker2} is the precisely the fibered product of the above 
two kernels over $V_{v_2}$.
Since $\dim V_{v_2} = k$, we conclude 
that the kernel of \eqref{eq:gluing-ker2} has dimension at least $k$, as 
desired.

It thus suffices to see that (c) implies (a). Since (a) may be checked node 
by node, it is enough
to show that for any $e_1 \in E(\Gamma)$ from $v_1$ to $v_2$,
the restriction of $((\sE^v,V^v)_v,(\vp_e)_e)$ to 
$Y_{v_1} \cup Y_{v_2}$
lies in $\cG^{k,\EHT}_{r,d,(d_{v_1},d_{v_2})}(Y_{v_1} \cup Y_{v_2})$. 
But for any $w \in V(G_{\I})$ which is equal to $d_{v'}-rb$ for all 
$v' \neq v_1,v_2$, by Lemma \ref{lem:restrict-kers} we have that
restriction to $Y_{v_1} \cup Y_{v_2}$ induces an injection
of the kernel of \eqref{eq:gluing-ker2} into the kernel of
\begin{equation}\label{eq:two-comp-gluing-ker}
V^{v_1}(w)\oplus V^{v_2}(w) \to 
\left.\left(\sE^{v_1}(w)\right)\right|_{\Delta'_{e_1}},
\end{equation}
so by the case of two components we conclude the desired statement.

We next observe that if we have an EHT limit series, and $w=w_v$ for
some $v \in V(\Gamma)$, then setting $X'=Y_v$ it follows immediately
from Lemma \ref{lem:restrict-kers} that the kernel of 
\eqref{eq:gluing-ker} always has dimension exactly $k$. Finally, we
wish to show that if an EHT limit series is refined, then
the dimension of the kernel of \eqref{eq:gluing-ker} is equal to
$k$ for all $w \in V(G_{\I})$.
We have already seen this in the case of two components.
Now, suppose we have a refined EHT limit series, and $w \in V(G_{\I})$ which
is equal to $d_{v''}-rb$ for all $v''$ except for two adjacent vertices
$v,v'$. Let $X'=Y_v \cup Y_{v'}$. Then the restriction to $X'$ is itself
a refined EHT limit series, so by the two-component case together with
Lemma \ref{lem:restrict-kers}, we conclude that the kernel of 
\eqref{eq:gluing-ker} has dimension $k$, as desired.
\end{proof}

\begin{rem}\label{rem:eht-no-stack}
As stated, the definition of $\cG^{k,\EHT}_{r,d,d_{\bullet}}(X)$
does not yield any obvious stack structure. Condition (I) is open, while
condition (II) may be expressed as a finite union of closed conditions,
but the gluing in condition (III) is subtler. A significant difficulty is 
that this condition cannot be applied uniformly, but only on individual
strata with given vanishing sequence. However, as a consequence of Lemma
\ref{lem:teixidor-translate}, we will see in \S \ref{sec:eht-stack}
below that in fact 
$\cG^{k,\EHT}_{r,d,d_{\bullet}}(X)$ does have a canonical stack structure.
\end{rem}

\subsection{Stack structures}\label{sec:eht-stack}
We now start introducing and relating various natural stack structures
which arise in connection with Eisenbud-Harris-Teixidor limit linear
series. First, from Lemma \ref{lem:teixidor-translate}, we may use the 
constructions of Appendix \ref{app:pushforward-detl} to define a stack 
structure on $\cG^{k,\EHT}_{r,d,d_{\bullet}}(X)$.

\begin{defn}\label{defn:eht-stack}
We endow $\cG^{k,\EHT}_{r,d,d_{\bullet}}(X)$ with the structure
of a locally closed substack of $\cP^k_{r,d_{\bullet}}(X)$ as follows:
an $S$-valued point 
$((\sE^v,\sV^v)_{v \in V(\Gamma)},(\vp_e)_{e \in E(\Gamma)})$
of $\cP^k_{r,d_{\bullet}}(X)$ is a point
of $\cG^{k,\EHT}_{r,d,d_{\bullet}}(X)$ if it lies in the open substack
satisfying (I) of Definition \ref{def:eht-lls}, and 
for each $w \in V(G_{\I})$,
it satisfies the closed condition that the $k$th vanishing locus 
(Definition \ref{def:pushforward-locus}) of
$$\pi_* \sE_w \to \bigoplus_{v \in V(\Gamma)} \pi_* \sE^v/\sV^v$$
is equal to all of $S$.
\end{defn}

\begin{rem} As a locally closed substack of $\cP^k_{r,d_{\bullet}}(X)$, it 
is evident that $\cG^{k,\EHT}_{r,d,d_{\bullet}}(X)$ is an Artin stack,
locally of finite type over $\Spec F$, and the morphism
$$\cG^{k,\EHT}_{r,d,d_{\bullet}}(X) \to \cM_{r,w_0,d_{\bullet}}(X)$$
is relatively representable by schemes which are projective, at least
locally on the target.
\end{rem}

\begin{rem} According to Lemma \ref{lem:teixidor-translate}, we could have
used all $w \in V(G_{\II})$ instead of $V(G_{\I})$ to obtain a (possibly
different) stack structure on $\cG^{k,\EHT}_{r,d,d_{\bullet}}(X)$; we chose
$V(G_{\I})$ largely because it seems closer in spirit to the original
Eisenbud-Harris-Teixidor definition.
\end{rem}

In order to carry out comparisons between stack structures, we now move 
from considering $K$-valued points to $S$-valued points. One of the
pleasant aspects of EHT limit series is that, due to the 
component-by-component and node-by-node nature of the definition, the
restriction to a subcurve of an EHT limit series is still an EHT limit 
series. We begin by observing that this property carries over to our
stack structures.

\begin{prop}\label{prop:eht-restrict} Let $X'$ be a connected subcurve
of $X$. Then restriction to $X'$ induces a morphism
$$\cG^{k,\EHT}_{r,d,d_{\bullet}}(X) \to 
\cG^{k,\EHT}_{r,d',d'_{\bullet}}(X'),$$
where $d'_{\bullet}$ denotes the restriction of $d_{\bullet}$ to the
components lying in $X'$, and $d'$ is determined by keeping $b$ fixed.
\end{prop}

\begin{proof} Let $\Gamma'$ be the dual graph of $X'$, and let $G'_{\I}$
be the associated graph as in Definition \ref{def:graphs}. Since we know the 
restriction morphism is defined on the set-theoretic level, the open 
condition (I) of Definition \ref{def:eht-lls} is 
necessarily preserved under restriction, and it is enough to prove that 
given an $S$-valued point
$((\sE^v,\sV^v)_v,(\vp_e)_e)$ of $\cG^{k,\EHT}_{r,d,d_{\bullet}}(X)$,
for each $w' \in V(G'_{\I})$,
the $k$th vanishing locus of
$$\pi_* \sE_{w'} \to \bigoplus_{v \in V(\Gamma')} \pi_* \sE^v/\sV^v$$
is equal to all of $S$. Let $w \in V(G_{\I})$
be the extension of $w'$ by $d_v-rb$ for all 
$v \not\in \Gamma'$. Then by hypothesis, the $k$th vanishing locus of
$$\pi_* \sE_w \to \bigoplus_{v \in V(\Gamma)} \pi_* \sE^v/\sV^v$$
is equal to all of $S$. Now, we have the natural restriction morphism
$\sE_w \to \sE_{w'}$, and the projection morphism
$\bigoplus_{v \in V(\Gamma)} \sE^v \to \bigoplus_{v \in V(\Gamma')} \sE^v$, 
and these clearly satisfy (I) and (II) of Corollary 
\ref{cor:pushforward-morphisms}. Moreover, Lemma \ref{lem:restrict-kers}
implies that (III) of Corollary
\ref{cor:pushforward-morphisms} is likewise satisfied, so we conclude the
desired statement.
\end{proof}

Next, recall that on 
the refined locus, the Eisenbud-Harris-Teixidor conditions could also be
used to give a stack structure. Our first comparison statement is 
that on the refined locus, the two stack structures agree. We begin by
making precise the latter stack structure.

\begin{defn}\label{def:eht-var}
Let $a^{\Gamma}$ consist of, for each pair $(e,v)$ of an edge 
of $\Gamma$ together with an adjacent vertex, a sequence 
$a_1^{e,v} \leq \dots \leq a_k^{e,v}$, such that for every edge $e$ with 
$v,v'$ the adjacent vertices, the sequences satisfy \eqref{eq:eht-compat} 
with equality. Define the locally closed substack 
$\cG^{k,\EHT}_{r,d,d_{\bullet},a^{\Gamma}}(X)$ of
$\cP^k_{r,d_{\bullet}}(X)$ to be the stack defined by the following 
conditions:
\begin{Ilist}
\itm the open condition (I) of Definition
\ref{def:eht-lls};
\itm the locally closed condition that for each adjacent pair $(e,v)$ of
$\Gamma$, the vanishing sequence on $Y_v$ at $\Delta'_e$ is exactly
equal to $a_1^{e,v},\dots,a_k^{e,v}$;
\itm the closed gluing condition that for each $w \in G_{\I}$ equal to
$d_{v''}-rb$ for all $v''$ other than some $v,v'$ adjacent to an edge $e$
of $\Gamma$, we have that the restrictions $\sV^v(w)|_{\Delta'_e}$ and 
$\sV^{v'}(w)|_{\Delta'_e}$
agree under the gluing map $\vp_e$.
\end{Ilist}
\end{defn}

\begin{prop}\label{prop:eht-refd-agrees} 
The open substack of $\cG^{k,\EHT}_{r,d,d_{\bullet}}(X)$ consisting of
refined EHT limit series is isomorphic to the (disjoint) union of
the substacks $\cG^{k,\EHT}_{r,d,d_{\bullet},a^{\Gamma}}(X)$ as $a^{\Gamma}$
ranges over allowable collections of sequences.
\end{prop}

\begin{proof} Because restricting to a locus with given (refined) vanishing
sequences $a^{\Gamma}$ yields an open substack of 
$\cG^{k,\EHT}_{r,d,d_{\bullet}}(X)$, it is enough to fix a choice of 
$a^{\Gamma}$ and work with the corresponding open locus. 

First, suppose we are given an $S$-valued point $((\sE^v,\sV^v)_v,(\vp_e)_e)$
of $\cG^{k,\EHT}_{r,d,d_{\bullet},a^{\Gamma}}(X)$; we wish to show that it is
also an $S$-valued point of $\cG^{k,\EHT}_{r,d,d_{\bullet}}(X)$, meaning that
for each $w \in V(G_{\I})$,
it satisfies the closed condition that the $k$th vanishing locus of
$$\pi_* \sE_w \to \bigoplus_{v \in V(\Gamma)} \pi_* \sE^v/\sV^v$$
is equal to all of $S$. According to Proposition 
\ref{prop:pushforward-detl-subbundles}, it is enough to show that
$\pi_* \sE_w$ has a rank-$k$ subbundle contained in the kernel of the map
to $\bigoplus_{v \in V(\Gamma)} \pi_* \sE^v/\sV^v$. We observe that due to the
form of $w$, we have $D_{w,v}$ supported at a single point for every
$v \in V(\Gamma)$. Thus, the condition on vanishing sequences implies that
we may view $\sV^v(w)$ as a subbundle of $\sV^v$ for each $v$. Similarly,
the gluing condition implies that for each $e$ and with $v,v'$ the adjacent 
vertices, $\sV^v(w)$ and $\sV^{v'}(w)$ have the same image (with constant 
rank) in 
$\left.\left(\sE^{h(e)}(w)\right)\right|_{\Delta'_e}$, so we conclude 
by inductively traversely $\Gamma$ that the map
$$\bigoplus_{v \in V(\Gamma)} \sV^v(w) \to \bigoplus_{e \in E(\Gamma)} 
\left.\left(\sE^{h(e)}(w)\right)\right|_{\Delta'_e}$$
has constant rank, and hence its kernel is a subbundle. We know that on 
points, the dimension of the kernel is equal to $k$, so we conclude that 
the kernel is a subbundle of rank $k$, which yields a rank-$k$ subbundle 
of $\pi_* \sE_w$ contained in the 
kernel of the map to $\bigoplus_{v \in V(\Gamma)} \pi_* \sE^v/\sV^v$, as 
desired.

For the opposite containment, suppose that $((\sE^v,\sV^v)_v,(\vp_e)_e)$
is an $S$-valued point of $\cG^{k,\EHT}_{r,d,d_{\bullet}}(X)$ lying
set-theoretically in $\cG^{k,\EHT}_{r,d,d_{\bullet},a^{\Gamma}}(X)$. As before
the open condition is automatically satisfied, so we need only verify that
(II) and (III) of 
Definition \ref{def:eht-lls} are verified (stack-theoretically), for 
which it suffices to work with the case $S$ local. 
Given $e$, let $v,v'$ be the adjacent vertices, and restrict to the
subcurve $Y_v \cup Y_{v'} \subseteq X$. By Proposition 
\ref{prop:eht-restrict}, we obtain an $S$-valued point of
$\cG^{k,\EHT}_{r,d',(d_v, d_{v'})}(Y_v \cup Y_{v'})$, which is to say that
for any $j$, if we set $w=(d_v-jr,d_{v'}-(b-j)r)$, the $k$th vanishing
locus of
$$\pi_* \sE_w \to \pi_* \sE^v/\sV^v \oplus \pi_* \sE^{v'}/\sV^{v'}$$
is all of $S$. On the other hand, since we are in the refined locus, by
Lemma \ref{lem:teixidor-translate} the $(k+1)$st vanishing locus is empty,
so by Proposition \ref{prop:pushforward-detl-subbundles}, we conclude that
the kernel $\sK_w$ is locally of free of rank $k$, and commutes with
base change. We know that set-theoretically, the rank of
$$\sV^v \to \left.\left(\sE^{v}\right)\right|_{j \Delta'_e}$$
is constant, equal to some $b^{v}_w$, and similarly the rank of
$$\sV^{v'} \to \left.\left(\sE^{v'}\right)\right|_{(b-j) \Delta'_e}$$
is constant, equal to some $b^{v'}_w$; to prove that
(II) of Definition \ref{def:eht-lls} is satisfied,
we wish to show that the ranks are scheme-theoretically constant, or 
equivalently, that the closed subschemes of $S$ on which they are less
than or equal to $b^{v}_w$ (respectively, $b^{v'}_w$) are
equal to all of $S$. We describe the argument for the first map, the second
being similar. It suffices 
to produce a subbundle of $\sV^v$ of rank $k-b^{v}_w$ contained in the 
kernel of the map in question. But $\sK_w$ maps into this kernel. At the 
closed point, we have that $\sK_w$ maps to $\sV^v$ with rank $k-b^{v}_w$
so using Nakayama's lemma we can construct a subbundle of $\sK_w$ of
rank $k-b^{v}_w$ which yields a subbundle of $\sV^v$ of the same rank,
giving the desired rank bound and proving that (II)
of Definition \ref{def:eht-lls} is satisfied. Similarly, we wish to 
prove that the rank of 
$$\sV^v(w) \oplus \sV^{v'}(w) \to \sE_w|_{\Delta'_e}$$
is at most $\#\{\ell:a^{(e,v)}_{\ell}=j\}=\#\{\ell:a^{(e,v')}_{\ell}=b-j\}$
scheme-theoretically everywhere on $S$. The rank of the source is
$2k-b^v_w-b^{v'}_w=k+\#\{\ell:a^{(e,v)}_{\ell}=j\}$, so it is enough
to produce a subbundle of the kernel of rank $k$.
But again, $\sK_w$ maps to the kernel, and does 
so injectively at the closed point, so we obtain the desired subbundle, and
conclude the asserted equality of substacks.
\end{proof}

The stack structures we have defined thus far are in most ways fully
satisfactory, but one facet is still missing: because the definition of
$\cP^k_{r,d_{\bullet}}(X)$ involves bundles and spaces of sections on
individual components of $X$, it does not make sense in smoothing families.
We thus conclude our examination of stacks structures on
$\cG^{k,\EHT}_{r,d,d_{\bullet}}(X)$ by giving a closely related construction
which is isomorphic to our original stack structure, but which works well 
in smoothing families. 

\begin{notn}\label{not:tildePkrd} 
Denote by $\widetilde{\cP}^k_{r,d_{\bullet}}(X)$ the stack over
$\cM_{r,w_0}(X)$ parametrizing, for each $v \in V(\Gamma)$, a
choice of rank-$k$ subbundle $\sV_{w_v}$ of $\pi_* \sE_{w_v}$.
\end{notn}

Thus, the only difference between $\widetilde{\cP}^k_{r,d_{\bullet}}(X)$
and $\cP^k_{r,d_{\bullet}}(X)$ is that instead of considering collections
$\sV^v$ of spaces of global sections of $\sE^v=\sE_{w_v}|_{Y_v}$ on each 
component $Y_v$, we take spaces of global sections of $\sE_{w_v}$ on all
of $X$. Because the degrees in question are extremal, we will see that
on the loci of interest to us, nothing is changed. We can then define
our alternate substack structure in a manner precisely parallel to the
previous one.

\begin{defn}\label{def:tildeG} We define the locally closed substack 
$\widetilde{\cG}^{k,\EHT}_{r,d,d_{\bullet}}(X)$ 
of $\widetilde{\cP}^k_{r,d_{\bullet}}(X)$ as follows:
an $S$-valued point 
$(\sE,(\sV_{w_v})_{v \in V(\Gamma)})$
of $\widetilde{\cP}^k_{r,d_{\bullet}}(X)$ is a point
of $\widetilde{\cG}^{k,\EHT}_{r,d,d_{\bullet}}(X)$ if it lies in the preimage
of the open substack 
$\cM_{r,w_0,d_{\bullet}}(X) \subseteq \cM_{r,w_0}(X)$,
and for each $w \in G_{\I}$, it satisfies the closed condition that the 
$k$th vanishing locus of
$$\pi_* \sE_w \to \bigoplus_{v \in V(\Gamma)} \pi_* \sE_{w_v}/\sV_{w_v}$$
is equal to all of $S$.
\end{defn}

\begin{prop}\label{prop:restriction-isom} Restriction of the $\sV_{w_v}$
to $Y_v$ induces a morphism from the preimage of
$\cM_{r,w_0,d_{\bullet}}(X)$ in
$\widetilde{\cP}^{k}_{r,d_{\bullet}}(X)$
to the stack
$\cP^{k}_{r,d_{\bullet}}(X)$. Moreover, this morphism induces an isomorphism
$$\widetilde{\cG}^{k,\EHT}_{r,d,d_{\bullet}}(X) \to
\cG^{k,\EHT}_{r,d,d_{\bullet}}(X).$$
\end{prop}

\begin{proof} All that needs to be checked to see that we obtain a
morphism to $\cP^k_{r,d_{\bullet}}(X)$ is that on the open subset in 
question, the for each $v$ we have that $\sV_{w_v}$ yields a rank-$k$ 
subbundle of $\pi_* \sE^v$ under composition with the restriction map
$\pi_* \sE_{w_v} \to \pi_* \sE^v$. For this, it is enough by
Lemma \ref{lem:grd-sub} (iii)
to check injectivity on points, and we obtain the desired
statement from Lemma \ref{lem:restrict-kers} with $X'=Y_v$, since the
kernels of the maps in question are identified with $H^0(X,\sE_{w_v})$ and
$H^0(Y_v,\sE^v)$, respectively. This proves the first assertion. For 
convenience, we denote the induced subbundle of $\pi_* \sE^v$ by
$\sV^v$.

Next, the image of this morphism is evidently contained in the open 
substack of $\cP^{k}_{r,d_{\bullet}}(X)$ satisfying (I)
of Definition \ref{def:eht-lls}, so in order to see
we get an induced morphism
$$\widetilde{\cG}^{k,\EHT}_{r,d,d_{\bullet}}(X) \to
\cG^{k,\EHT}_{r,d,d_{\bullet}}(X),$$
it suffices to see that if the $k$th vanishing locus of
$$\pi_* \sE_w \to \bigoplus_{v \in V(\Gamma)} \pi_* \sE_{w_v}/\sV_{w_v}$$
is equal to all of $S$, then the same is true of
$$\pi_* \sE_w \to \bigoplus_{v \in V(\Gamma)} \pi_* \sE^v/\sV^v.$$
But because the second map factors through the first, this is immediate
from Corollary \ref{cor:pushforward-morphisms}. 

It remains to see that 
this induced morphism is an isomorphism, which amounts to the assertion
that for a given $S$-valued collection of subbundles 
$\sV^v \subseteq \pi_* \sE^v$, there exists a unique collection of
subbundles $\sV_{w_v} \subseteq \pi_* \sE_{w_v}$ satisfying our vanishing
condition for all $w$ and restricting to the given ones.
But by definition the morphism 
\begin{equation}
\pi_* \sE_{w_v} \to \bigoplus_{v \in V(\Gamma)} \pi_* \sE^v/\sV^v
\end{equation}
has $k$th vanishing locus equal to $S$, and by Lemma 
\ref{lem:teixidor-translate} its $(k+1)$st vanishing locus is empty,
so Proposition \ref{prop:pushforward-detl-subbundles} implies that
the kernel $\sK_v$ is a rank-$k$ subbundle, which we claim is the desired
$\sV_{w_v}$. But the restriction of $\sK_{v}$ to $Y_v$ is
contained in $\sV^v$ by construction, and it follows from the injectivity of 
$\pi_* \sE_{w_v} \to \pi_* \sE^v$ on points together with Lemma
\ref{lem:grd-sub} (iii) and (iv) 
that in fact $\sK_v$ restricts to $\sV^v$. The same injectivity also implies
that $\sV_{w_v}$ is uniquely determined by the condition that it restricts
to $\sV^v$.
Finally, it follows that the morphism 
$\pi_* \sE_{w_v}/\sV_{w_v} \to \pi_* \sE^v/\sV^v$ is injective, even
after arbitrary base change. Thus, for any $w \in V(G_{\I})$, the last
statement of Proposition \ref{prop:pushforward-detl-defined} together with 
the hypothesis that the $k$th vanishing locus of
$$\pi_* \sE_w \to \bigoplus_{v \in V(\Gamma)} \pi_* \sE^v/\sV^v$$
is equal to all of $S$ implies that the same is true of
$$\pi_* \sE_w \to \bigoplus_{v \in V(\Gamma)} \pi_* \sE_{w_v}/\sV_{w_v}.$$
We have thus shown the existence and uniqueness of the desired lift.
\end{proof}

\subsection{Comparison results}
We now compare our linked linear series stacks to the stack of EHT limit
series. We begin with the following obvious observation.

\begin{prop}\label{prop:forgetful-map} There is a forgetful morphism
\begin{equation}\label{eq:forget} \cG^{k,\I}_{r,d,d_{\bullet}}(X) \to
\cP^k_{r,d_{\bullet}}(X),
\end{equation}
given by
sending a tuple $(S,\sE,(\sV_w)_{w \in V(G_{\I})})$ to the tuple
$((\sE\otimes \sO_{w_0,w_v})|_{Y_v},\sV_{w_v}|_{Y_v})_{v \in V(\Gamma)}$,
together with the natural gluing maps at the nodes induced by $\sE$.
\end{prop}

The following definition will play an important role in our comparison
results. The purpose of \S \ref{sec:compatible} below is to develop
a robust criterion for when the conditions of the definition are satisfied.

\begin{defn}\label{def:constrained} 
A $K$-valued point $((\sE^v,V^v)_v,(\vp_e)_e)$ of 
$\cG^{k,\EHT}_{r,d,d_{\bullet}}(X)$ is said to be \textbf{constrained}
if there exist and $w_1,\dots,w_k \in V(\bar{G}_{\II})$
(not necessarily distinct) and for $i=1,\dots,k$ a vector
$v_i$ in the kernel of 
\eqref{eq:gluing-ker2} (in index $w_i$) such that the following conditions 
are satisfied:
\begin{Ilist}
\itm for each $w_i$, the kernel of \eqref{eq:gluing-ker2} has dimension
exactly $k$;
\itm for each $v \in V(\Gamma)$, the images of $v_1,\dots,v_k$ in $V^v$
form a basis.
\end{Ilist}
\end{defn}

We then have:

\begin{cor}\label{cor:constrained-open} The constrained points form an
open subset of $\cG^{k,\EHT}_{r,d,d_{\bullet}}(X)$, and in the definition,
the kernel of \eqref{eq:gluing-ker2} may be replaced by the kernel of
\eqref{eq:gluing-ker}.
\end{cor}

\begin{proof} The latter statement is Lemma \ref{lem:kers-same}. We then
conclude openness of condition (I) because the kernels have dimension at
least $k$ everywhere by Lemma \ref{lem:teixidor-translate}, so condition
(I) describes the complement of a closed substack by 
Proposition \ref{prop:pushforward-detl-defined}. On the other hand,
condition (II) is open by Nakayama's lemma.
\end{proof}

The main comparison theorem is the following:

\begin{thm}\label{thm:compare} The morphism
\eqref{eq:forget} factors through and surjects onto the locally closed 
substack $\cG^{k,\EHT}_{r,d,d_{\bullet}}(X)$.
Moreover, on the preimage of the refined locus 
$\cG^{k,\EHT,\refn}_{r,d,d_{\bullet}}(X)$ we have that \eqref{eq:forget} is
an isomorphism.

Under the additional hypothesis of Situation \ref{sit:basic-ii}, if we 
compose with the forgetful morphism of Corollary \ref{cor:forget}, the
resulting morphism
\begin{equation}\label{eq:forget-comp} \cG^{k,\II}_{r,d,d_{\bullet}}(X,\theta_{\bullet}) \to
\cG^{k,\EHT}_{r,d,d_{\bullet}}(X),
\end{equation}
is an isomorphism on the preimage of the constrained locus
of $\cG^{k,\EHT}_{r,d,d_{\bullet}}(X)$.
\end{thm}

\begin{proof} 
That \eqref{eq:forget} factors through $\cG^{k,\EHT}_{r,d,d_{\bullet}}(X)$
follows immediately from the definitions together with Proposition 
\ref{prop:pushforward-detl-subbundles}. 

To check surjectivity, given a point $((\sE^v,V^v)_v,(\vp_e)_e)$ of 
$\cG^{k,\EHT}_{r,d,d_{\bullet}}(X)$, we will construct a point of 
$\cG^{k,\I}_{r,d,d_{\bullet}}(X)$ mapping to it as follows:
let $\sE$ be the vector bundle of multidegree $w_0$ determined by 
gluing (the appropriate twists of) the $\sE^v$ via the maps $\vp_e$;
for $w \in V(G_{\I})$, let $\widetilde{V}_w$ be the kernel of
\eqref{eq:gluing-ker}. Then, according to Lemma 
\ref{lem:teixidor-translate}, we have that $\widetilde{V}_w$ has 
dimension at least $k$ for all $w$, and has dimension exactly $k$ if
$w=w_v$ for some $v \in V(\Gamma)$. The spaces $\widetilde{V}_w$ are
visibly linked, so it is enough to show that we can produce a system
of $k$-dimensional subspaces $V_w \subseteq \widetilde{V}_w$ which remain 
linked. Because we will have $V_{w_v}=\widetilde{V}_{w_v}$ for all
$v \in V(\Gamma)$, it is enough to work with one pair $(v_1,v_2)$ of adjacent
edges at a time, considering only $w=(i_v)_v$ with $i_v=d-rb$ for all
$v \neq v_1,v_2$. We will prove the desired statement by showing that we
can always reduce the dimension of some $\widetilde{V}_w$ while maintaining
linkage, until we reach the situation that all spaces have dimension $k$.

Accordingly, with $(v_1,v_2)$ fixed, let $e \in E(\Gamma)$ be the edge 
joining them, and take the $w$ with $i_{v_1}$ minimal such 
that $\dim \widetilde{V}_w> k$. We will construct a subspace 
$\widetilde{V}'_w \subseteq \widetilde{V}_w$ of dimension one less, and such
that if $w',w''$ are the vertices adjacent to $w$ in $G_{\I}$, we still have
$\widetilde{V}'_w$ containing the images of $\widetilde{V}_{w'}$ and
$\widetilde{V}_{w''}$. Iterating this process will thus yield the desired
system of subspaces. There are two cases to consider. First, if 
$\widetilde{V}_w$ contains a section $s$ which is nonvanishing at 
$\Delta'_e$, let $\ell \subseteq \sE_w|_{\Delta'_e}$ be the line generated
by $s$, and $\sE'_w \subseteq \sE_w$ the kernel of the map
$\sE_w \to \left(\sE_w|_{\Delta'_e}\right)/\ell$. Then set 
$\widetilde{V}'_w$ to be the subspace of $\widetilde{V}_w$ consisting of
sections lying in $\sE'_w$. 
Otherwise, if $\widetilde{V}_w$ consists entirely of sections vanishing
at $\Delta'_e$, suppose we have ordered $w'$ and $w''$ so that the index
$i'_{v_1}$ for $w'$ is less than $i_{v_1}$. 
Let $\widetilde{V}_w|_{v_i}$ (respectively,
$\widetilde{V}_{w'}|_{v_i}$) denote the image of $\widetilde{V}_w$
(respectively, $\widetilde{V}_{w'}$) under restriction to $Y_{(e,v_i)}$ for
$i=1,2$. Then, by minimality of $i_{v_1}$, we see that we must have
$\dim \widetilde{V}_w|_{v_1} > \dim \widetilde{V}_{w'}|_{v_1}$,
and in particular $\widetilde{V}_w|_{v_1}$ strictly contains the image of 
$\widetilde{V}_{w'}|_{v_1}$ under the natural inclusion; let $V' \subseteq
\widetilde{V}_w|_{v_1}$ be a codimension-$1$ subspace containing this image;
then we may set $\widetilde{V}'_w$ to be the preimage of $V'$ under the 
restriction map.  
This completes the proof of surjectivity.

Next, to see that we get an isomorphism on the refined locus, suppose that
$((\sE^v,\sV^v)_v,(\vp_e)_e)$ is an $S$-valued point of 
$\cG^{k,\EHT,\refn}_{r,d,d_{\bullet}}(X)$; we wish to show that
there exists a unique $S$-valued point 
$(\sE,(\sV_w)_{w \in V(G_{\I})})$ of $\cG^{k,\I}_{r,d,d_{\bullet}}(X)$ 
mapping to it. By the last part of Lemma \ref{lem:teixidor-translate},
we have that for every $w \in V(G_{\I})$,
the kernel of \eqref{eq:gluing-ker} has dimension exactly $k$. It follows
from the definition of the stack $\cG^{k,\EHT}_{r,d,d_{\bullet}}(X)$
together with Proposition \ref{prop:pushforward-detl-subbundles} that the 
kernel is a subbundle
of rank $k$, so we must have $\sV_w$ equal to this kernel. This gives
uniqueness, and also existence, since it is evident that the kernels
of \eqref{eq:gluing-ker} as $w$ varies must satisfy the required linkage
condition.

For \eqref{eq:forget-comp}, if we have an $S$-valued point 
$((\sE^v,\sV^v)_v,(\vp_e)_e)$ of 
$\cG^{k,\EHT}_{r,d,d_{\bullet}}(X)$ contained in the constrained locus, we
need to show that there exists a unique $S$-valued point 
$(\sE,(\sV_w)_{w \in V(G_{\II})})$ of 
$\cG^{k,\II}_{r,d,d_{\bullet}}(X,\theta_{\bullet})$ 
mapping to it. It is enough to prove this locally,
so we reduce to the case that there is a fixed collection of $w_i$ 
as in Definition \ref{def:constrained} which works for every point of $S$.
For each $i$, we first claim that the $k$th vanishing
locus of \eqref{eq:gluing-ker} in multidegree $w_i$ is equal to all of $S$.
The argument for this is a variant
of the induction argument that (a) implies (b) in Lemma 
\ref{lem:teixidor-translate}. Indeed, if we only have two components the
statement is true tautologically, since necessarily 
$w_i \in V(\bar{G}_{\II})=V(G_{\I})$. If we have at least three components,
choose a component $Y_{v_1}$ meeting the rest of $X$ at only a single node
$\Delta'_{e_1}$. Let $X'$ be the complement of $Y_{v_1}$ in $X$, and let 
$Y$ be the union of $Y_{v_1}$ and $Y_{v_2}$, where $v_2$ is the other edge 
adjacent to $e_1$. Denote by $\Gamma'$ the dual graph of $X'$.  
Restricting to $X'$ and $Y$, and letting $w'$ and $w_{12}$ be the multidegrees
induced by $w_i$, by the induction hypothesis we have that the 
$k$th vanishing loci of the maps
$$
\pi_* \sE_{w'} \to \bigoplus_{v \in V(\Gamma')} (\pi_* \sE^v)/\sV^v
$$
and
$$
\pi_* \sE_{w_{12}} \to  
(\pi_* \sE^{v_1})/\sV^{v_1}  \oplus (\pi_* \sE^{v_2})/\sV^{v_2}
$$
both contain all of $S$. However, the kernel of 
\eqref{eq:gluing-ker} is the precisely the fibered product of the above 
two kernels over $\sV^{v_2}$. Since $\sV^{v_2}$ has rank $k$, we conclude 
from Corollary \ref{cor:pushforward-fibered-prod} that the $k$th vanishing
locus of \eqref{eq:gluing-ker} is all of $S$, as claimed.
On the other hand, the definition of constrained point 
implies that the
$(k+1)$st vanishing locus of \eqref{eq:gluing-ker} (still in multidegree
$w_i$) is empty on the set of constrained points,
so from Proposition
\ref{prop:pushforward-detl-subbundles} we conclude that the kernel of
\eqref{eq:gluing-ker} in multidegree $w_i$ is a subbundle of rank $k$.
Thus if the desired lift exists, we must have $\sV_{w_i}$ equal
to the kernel of \eqref{eq:gluing-ker}. 

Next, working locally around a given point $y \in S$ and letting $v_i$ be as 
in the definition of constrained, we may lift each $v_i$ to a 
$\tilde{v}_i \in \sV_{w_i}$, and we claim that for any $w \in V(G_{\II})$, 
locally around $y$ we have that $\sum_i f_{w_i,w} \tilde{v}_i$ is a rank-$k$ 
subbundle of $\pi_* \sE_w$. According to Lemma \ref{lem:grd-sub} (iii),
it is enough to see that the $f_{w_i,w} v_i$ are linearly independent
in $\pi_* \sE_w|_y$. 
This follows from the same argument as Lemma \ref{lem:simple-crit}.
Thus, if we set $\sV_w= \sum_i f_{w_i,w} \tilde{v}_i$, we see that we 
obtain the desired point of 
$\cG^{k,\II}_{r,d,d_{\bullet}}(X,\theta_{\bullet})$. 
Indeed, because we are working with a fixed choice of basis vectors,
the linkage condition is immediate from the fact that for two
paths $P,P'$ with the same start and endpoints, the maps $f_P$ and $f_{P'}$
are related by scalars. In addition, since our second description
of each $\sV_{w_i}$ must be contained in the first, Lemma \ref{lem:grd-sub}
(iv) implies that they agree. Finally, since any choice of $\sV_w$ must 
contain the chosen one, again using Lemma \ref{lem:grd-sub} (iv) we conclude 
the desired uniqueness of our lift.
\end{proof}

\begin{ex}\label{ex:bad-compare}
We give a simple example showing that with three or more components, 
the seemingly rather restrictive ``constrained'' hypothesis is indeed
necessary. Specifically, we see that
the forgetful map \eqref{eq:forget-comp}
is not an isomorphism even in the rank-$1$ case on the locus 
of refined and simple limit series.

Indeed, consider a curve with
three rational components $Y_1,Y_2,Y_3$, with $Y_1$ glued to $Y_2$ at a
node $P_1$, and $Y_2$ glued to $Y_3$ at a node $P_2$. Thus, in this case
all the line bundles are uniquely determined. We consider simple
linked $g^1_2$'s generated by a section $s_1$ in multidegree $(1,0,1)$ and
a section $s_2$ in multidegree $(0,2,0)$. Recall that in order for $s_1$ and
$s_2$ to serve as generators for a simple point, neither of them can vanish
uniformly on any of the $Y_i$. The choices of $s_1$ may be 
described as follows: if $x_1$ denotes the section (unique up to scaling)
which is nonzero on $Y_1$ but identically zero on $Y_2$ and $Y_3$, and $x_3$
denotes the section nonzero on $Y_3$ but identically zero on $Y_1$ and $Y_2$,
then $s_1$ is unique up to scaling and adding multiples of $x_1$ and $x_3$.
On the other hand, for a fixed choice of $s_2$, the image of $s_2$ in 
multidegree $(1,0,1)$ vanishes on $Y_2$, and is thus of the 
form $c_1 x_1 + c_3 x_3$ for some non-zero $c_1,c_3$. Thus, we see that
for a fixed choice of $s_2$, modifying $s_1$ by multiples of $x_1$ and 
$x_3$ changes the resulting space in multidegree $(1,0,1)$, and we obtain
a one-parameter family of $2$-dimensional subspaces in this way. On the 
other hand, such 
modifications do not affect the subspaces in any other multidegree: the
image of $s_1$ in multidegree $(0,2,0)$ vanishes uniformly on $Y_1$ and
$Y_3$, so it
is enough to check the multidegrees $(1,1,0)$ and $(0,1,1)$, and 
we see that modifying $s_1$ by multiples of $x_1$ and $x_3$ changes the
image in multidegree $(1,1,0)$ by a multiple of (the image of) $x_1$. But
the image of $s_2$ in multidegree $(1,1,0)$ is also a multiple of the image
of $x_1$, so the resulting $2$-dimensional subspace is unchanged. The same
holds for $(0,1,1)$, and we see that in this case the comparison morphism
\eqref{eq:forget-comp} contracts curves.
\end{ex}

\section{Limit series on chains of curves}\label{sec:compatible}
\label{sec:chains}

As we have discussed in Example \ref{ex:bad-compare},
unlike the case of two components treated in \cite{os8} our 
comparison morphism from type II linked linear series to EHT limit linear 
series is not in general an isomorphism over the refined locus. Thus, for
our applications it becomes very important to identify loci over which
the comparison morphism is an isomorphism. It turns out that for chains
of curves, there is a natural open subset, which we refer to as the
``chain-adaptable'' locus, over which we obtain the desired behavior.
This locus contains all examples considered by Teixidor i Bigas in the
existence arguments of \cite{te1}, \cite{te3}, \cite{te5}.

As in the previous section, here we assume throughout that $B=\Spec F$,
and $X$ is a projective curve over $B$.

\subsection{Pairs of vanishing sequences on smooth curves}
We begin with some observations on smooth curves with pairs of marked points.

If $(\sE,V)$ is a $\fg^r_d$ on a smooth projective curve, then we denote
the vanishing sequence at a point $P$ by 
$a_1^{(\sE,V)}(P),\dots,a_k^{(\sE,V)}(P)$.

\begin{defn}\label{def:adapted}
Let $X$ be a smooth projective curve over $\Spec F$, and $(\sE,V)$ a 
$\fg^k_{r,d}$ on $X$. Given points $P,Q \in X(F)$, we say that a basis 
$s_1,\dots,s_k \in V$
is \textbf{$(P,Q)$-adapted} if $\ord_P s_i = a_i^{(\sE,V)}(P)$ and 
$\ord_Q s_i = a_{k+1-i}^{(\sE,V)}(Q)$ for $i=1,\dots,k$. We say that $(\sE,V)$
is \textbf{$(P,Q)$-adaptable} if there exists a $(P,Q)$-adapted basis of $V$.
\end{defn}

The following lemma will be useful:

\begin{lem}\label{lem:vanishing-seq}
Suppose $X$ is smooth over $\Spec F$, and $(\sE,V)$ is a $\fg^k_{r,d}$ on
$X$. Given $P \in X(F)$, let $a_1,\dots,a_k$ be the 
vanishing sequence of $(\sE,V)$ at $P$, and $s_1,\dots,s_k \in V$ such that
the order of vanishing of $s_i$ at $P$ is $a_i$. Then the following are
equivalent:
\begin{alist}
\itm the $s_i$ form a basis for $V$;
\itm for each $a$, the set of $s_i$ vanishing to order $a$ at $P$ is a
basis for the fiber $V(-aP)/V(-(a+1)P)$;
\itm for each $a$, the set of $s_i$ vanishing to order $a$ at $P$ is 
linearly independent in the fiber $V(-aP)/V(-(a+1)P)$;
\itm for each $a$, the set of $s_i$ vanishing to order $a$ at $P$ spans
the fiber $V(-aP)/V(-(a+1)P)$.
\end{alist}
\end{lem}

\begin{proof} By the definition of vanishing sequence, we have that (b),
(c), and (d) are equivalent. On the other hand, given (c), and a 
nonzero linear combination $s=\sum_i c_i s_i$, if $i_0$ is minimal with 
$c_{i_0} \neq 0$, we see that $s$ must vanish to order precisely $a_{i_0}$
at $P$, and in particular is not the zero section. Thus, (c) implies (a).
Finally, suppose we have (a). We prove (d) by induction on $a$. The
statement is trivial for $a=0$, so suppose $a$ is general, and we know
the desired statement for all orders strictly less than $a$. Then again
using the definition of vanishing sequence, the $s_i$ vanishing to orders
$a'<a$ form bases of $V(-a'P)/V(-(a'+1)P)$, so no non-zero linear combination
of them is in $V(-aP)$. Since the $s_i$ are a basis of $V$, linear 
combinations of them span $V(-aP)/V(-(a+1)P)$, and we conclude that it
suffices to take linear combinations of $s_i$ vanishing to order exactly
$a$, as desired.
\end{proof}

\begin{prop}\label{prop:adapted} Suppose $X$ is smooth over $\Spec F$, and
we are given $P,Q \in X(F)$. Given also a pair $(\sE,V)$, for any $a,b \geq 0$
we have
\begin{equation}\label{eq:adaptable-ineq}
\dim V(-aP-bQ) \geq \#\{i:a_i^{(\sE,V)}(P) \geq a \text{ and } 
a_{k+1-i}^{(\sE,V)} (Q) \geq b\}.
\end{equation}
Moreover, $(\sE,V)$ is $(P,Q)$-adaptable if and only if for all $a,b \geq 0$ 
we have that \eqref{eq:adaptable-ineq} is satisfied with equality.
In particular, the $(P,Q)$-adaptable pairs are open in the moduli stack of
all pairs with given vanishing sequences at $P$ and at $Q$.
\end{prop}

\begin{proof} We first verify that the asserted inequality holds in
general. The main observation is that if the right-hand side of 
\eqref{eq:adaptable-ineq} is positive, then we have
$$\{i:a_{k+1-i}^{(\sE,V)}(Q) < b \} \subseteq \{i:a_i^{(\sE,V)}(P) \geq a\}.$$
Indeed, positivity implies there exists $i_0$ such that
$a_{i_0}^{(\sE,V)}(P) \geq a$ and $a_{k+1-{i_0}}^{(\sE,V)} (Q) \geq b$.
Then if $a_{k+1-i}^{(\sE,V)}(Q)<b$, we have $i>i_0$, so we get the desired
statement.
Now, if the right-had side of \eqref{eq:adaptable-ineq} is zero, there is
nothing to prove. Otherwise, because $V(-aP-bQ)=V(-aP) \cap V(-bQ)$, we have
\begin{align*} 
\dim V(-aP-bQ) & \geq \dim V(-aP) + \dim V(-bQ)-k \\
& = \#\{i:a_i^{(\sE,V)}(P) \geq a\} + \# \{i:a_i^{(\sE,V)}(Q) \geq b\} -k \\
& = \#\{i:a_i^{(\sE,V)}(P) \geq a\} - \# \{i:a_i^{(\sE,V)}(Q) < b\} \\
& = \#\{i:a_i^{(\sE,V)}(P) \geq a \text{ and } 
a_{k+1-i}^{(\sE,V)} (Q) \geq b\},
\end{align*}
where the last equality follows from our observation above.

It remains to check that equality in \eqref{eq:adaptable-ineq} is 
equivalent to $(P,Q)$-adaptability. First suppose that we have a 
$(P,Q)$-adapted basis $s_1,\dots,s_k$.
Then by Lemma \ref{lem:vanishing-seq}, there cannot
be any cancellation in orders of vanishing of linear combinations of the
$s_i$, so an element of $V(-aP-bQ)$ must be a linear combination of 
basis vectors $s_i$ with $i$ satisfying $a_{i}^{(\sE,V)} \geq a$ and 
$a_{k+1-i}^{(\sE,V)} \geq b$. We thus conclude equality in 
\eqref{eq:adaptable-ineq}.

Conversely, if we have equality in \eqref{eq:adaptable-ineq}, we wish 
to construct a $(P,Q)$-adapted basis $s_1,\dots,s_k$. We proceed via a
double downward induction on pairs $(a,b)$, with $a,b \geq 0$, and bounded 
above by $a_k^{(\sE,V)}(P)$ and $a_k^{(\sE,V)}(Q)+1$, respectively. Given
$(a,b)$, denote by $I_{(a,b)}$ the set of $i$ such that
$a_i^{(\sE,V)}(P)\geq a$ and $a_{k+1-i}^{(\sE,V)}(Q) \geq b$.
For each such pair $(a,b)$, and every $i$ with either $a_i^{(\sE,V)}(P)>a$ 
or $a_i^{(\sE,V)}(P)=a$ and $a_{k+1-i}^{(\sE,V)}(Q) \geq b$, 
we will construct $s_i \in V$ 
having the orders of vanishing at $P$ and $Q$ required for a $(P,Q)$-adapted 
basis, and satisfying further the condition that for all $(a',b')$ with
either $a'>a$ or $a'=a$ and $b' \geq b$, the set of $s_i$ with 
$i \in I_{(a',b')}$ is a basis of
$V(-a'P-b'Q)$. If we have constructed such $s_i$ for the pair $(0,0)$, then
in particular we have a $(P,Q)$-adapted basis, as desired. We start the
induction with $(a,b)$ at their maximal values, so that $V(-aP-bQ)=0$, 
and the desired conditions are tautologically satisfied. For each value
of $a$, we decrease $b$ by $1$ until we reach $b=0$. We then decrease $a$
by $1$, and reset $b$ to its maximal value. At each step, we use the $s_i$
we have constructed in the previous step, adding new ones as necessary.

In the case that $a$ remains constant and $b$ has decreased, by hypothesis
we already have $s_i$ for all $i$ with either $a_i^{(\sE,V)}(P)>a$ or 
$a_i^{(\sE,V)}(P)=a$ and $a_{k+1-i}^{(\sE,V)}(Q) > b$, 
and satisfying the desired basis condition for $(a',b')$ with $a'>a$
or $a'=a$ and $b'>b$. Thus, we need only extend the $s_i$ to include
$i$ with $a_i^{(\sE,V)}(P)=a$ and $a_{k+1-i}^{(\sE,V)}(Q) = b$. Our claim
is that the hypothesis that we have equality in \eqref{eq:adaptable-ineq}
implies that the number of such $i$ is precisely equal to 
$$\dim V(-aP-bQ)/\spn(V(-aP-(b+1)Q),V(-(a+1)P-bQ)),$$ 
in which case we may choose the new $s_i$ to induce a basis for this space.
But 
\begin{multline*}\dim \spn(V(-aP-(b+1)Q),V(-(a+1)P-bQ))\\
= \dim V(-aP-(b+1)Q)+ \dim V(-(a+1)P-bQ) - \dim V(-(a+1)P-(b+1)Q),
\end{multline*}
so the claim follows easily by applying equality in \eqref{eq:adaptable-ineq}
to $(a,b)$, $(a,b+1)$, $(a+1,b)$ and $(a+1,b+1)$. We then check the desired
condition on pairs $(a',b')$: if $a'=a$ and $b'=b$, then appying the
cases $(a,b+1)$ and $(a+1,b)$ and $(a+1,b+1)$ from the previously chosen 
$s_i$ we see that 
the previous $s_i$ with $i \in I_{(a,b+1)} \cup I_{(a+1,b)}$ form a basis 
for $\spn(V(-aP-(b+1)Q),V(-(a+1)P-bQ)$,
and thus together with the new
$s_i$ we obtain a basis for $V(-aP-bQ)$, as desired. On the other hand,
if $a'>a$ or $a'=a$ and $b'>b$, then the new $s_i$ are irrelevant, and
we have the desired condition from the induction hypothesis. We thus
have constructed the desired $s_i$ in this case.

In the case that $a$ has decreased and $b$ is set to $a_k^{(\sE,V)}(Q)+1$,
there are no new allowable values of $i$, so it suffices to observe that
the previously chosen $s_i$ satisfy the desired condition on pairs
$(a',b')$ with either $a'>a$ or $a'=a$ and $b'\geq b$. But the latter
possibility is vacuous, since in this case $V(-a'P-b'Q)=0$ and there are 
no allowable indices $i$. On the other hand, if $a'>a$ we are in the cases
handled by the induction hypothesis. Thus, the induction hypothesis for
$(a+1,0)$ implies the desired statement for $(a,a_k^{(\sE,V)}(Q)+1)$,
as we wished to prove.
\end{proof}

\subsection{Chains of curves}
We now return to studying reducible curves, and more particularly, chains
of smooth projective curves.

\begin{defn}\label{def:chain-adapted}
Let $X$ be a curve consisting of a chain of smooth projective curves
$Y_1,\dots,Y_n$ over $\Spec F$, with $P_i,Q_i \in Y_i(F)$ for each $i$,
and the point $Q_i$ on $Y_i$ glued to $P_{i+1}$ on $Y_{i+1}$. Then a
refined Eisenbud-Harris-Teixidor limit linear series on $X$ is \textbf{chain 
adaptable} if the pair induced by restriction to each $Y_i$ ($i=2,\dots,n-1$)
is $(P_i,Q_i)$-adaptable.
\end{defn}

Since imposing particular vanishing sequences at the nodes is open on the
refined locus of Eisenbud-Harris-Teixidor limit linear series, Proposition 
\ref{prop:adapted} implies:

\begin{prop}\label{prop:chain-adapted-open} The locus of chain-adaptable 
EHT limit series is open in $\cG^{k,\EHT}_{r,d,d_{\bullet}}(X)$.
\end{prop}

More substantively, we have the following result.

\begin{prop}\label{prop:chain-adapted-basis} If 
$((\sE^i,V^i)_{i=1,\dots,n},(\vp_i)_{i=1,\dots,n-1})$ 
is a $K$-valued chain-adaptable EHT limit linear series on $X$, then 
there exist bases $s^i_j$ of $V^i$ for each 
$i=1,\dots,n$, and permutations
$\sigma_i \in S_k$ for $i=1,\dots,n-1$, such that:
\begin{ilist} 
\itm for each $i=2,\dots,n-1$, the $s^i_j$ form a $(P_i,Q_i)$-adapted 
basis of $V^i$;
\itm for each $i=1,\dots,n-1$,
and each $j=1,\dots,k$, we have 
$\ord_{Q_i}(s^i_j)+\ord_{P_{i+1}} s^{i+1}_{\sigma_i(j)}=b$,
and $s^i_j$ glues to $s^{i+1}_{\sigma_i(j)}$ under 
$\vp_i$.
\end{ilist}
\end{prop}

\begin{proof} 
The argument is by induction, but requires some additional notation. For 
each $i$, let $a^i$ be the vanishing sequence of $V^i$
at $P_i$, and let $b^i$ be the vanishing sequence at $Q_i$,
so that $b^i_j=b-a^{i+1}_{k+1-j}$ for $i=1,\dots,n-1$ and $j=1,\dots,k$.
For $i=1,\dots,n-1$, let $r^i_1,\dots,r^i_{\ell_i}$ denote the number
of repetitions in the sequence $a^{i+1}$, so that 
$\sum_{j=1}^{\ell_i} r^i_j=k$. 
It will also be convenient to write $R^i_j:= \sum_{\ell=1}^{j-1} r^i_{\ell}$ 
(with $R^i_1=0$), so that the distinct values of the sequence $a^{i+1}$ are 
given by $a^{i+1}_{R^i_1+1},\dots,a^{i+1}_{R^i_{\ell_i}+1}$, and the
distinct values of $b_i$ are given by 
$b^i_{k-R^i_{\ell_i}},\dots,b^i_{k-R^i_1}$.
For $i=1,\dots,n-1$ and $j=1,\dots,\ell_i$,
denote by $W^i_j$ the $r^i_j$-dimensional subspace of the fiber 
$\sE^i(-(b-a^{i+1}_{R^i_j+1}) Q_i)|_{Q_i}$
induced by $V^i$.
Thus, given a $(P_i,Q_i)$-adapted basis $s^i$, we have
that the restrictions of $s^i_{R^i_j+1},\dots,s^i_{R^i_j+r^i_j}$ to $Q_i$ form
a basis of $W^i_j$.

We argue by induction on $i$, showing that for each $i$, there exist
bases of $V^1,\dots,V^i$ as in the statement, and that furthermore
we have the following flexibility: for each $j=1,\dots,\ell_i$ 
there exists an ordering $\alpha_j$ of the 
$s^i_{R^i_j+1},\dots,s^i_{R^i_j+r^i_j}$ 
so that given any basis $e_1,\dots,e_{r^{i}_j}$ of $W^{i}_j$ obtained from the
restrictions to $Q_i$ of $s^i_{R^i_j+1},\dots,s^i_{R^i_j+r^i_j}$ via
a change of basis which is upper triangular under the ordering $\alpha_j$,
without modifying the permutations $\sigma_j$
we may modify our choices of the $s^{i'}$ for $i' \leq i$ so that the 
restriction of $s^i_{R^{i}_j+\ell}$ to $Q_i$ is $e_{\tau(\ell)}$ for 
$\ell=1,\dots,r^{i}_j$. That is, we claim that under a suitable ordering,
we can realize arbitrary upper triangular change of bases to the images of
the $s^i_j$ in the spaces $W^{i}_\ell$. In the base case $i=1$, there is 
nothing to show, as in fact the $s^1_j$ may be chosen to give arbitrary 
bases of the $W^1_{\ell}$, so any ordering suffices. Suppose now that the 
desired statement is true 
for $i-1$, and we wish to show it for $i$. For $j=1,\dots,\ell_{i-1}$, let 
$\alpha'_j$ be the ordering on 
$s^{i-1}_{R^{i-1}_j+1},\dots,s^{i-1}_{R^{i-1}_j+r^{i-1}_j}$ 
supplied by induction.

In order to prove the inductive statement, we first consider all 
possibilities for $(P_i,Q_i)$-adapted bases of $V^i$. If 
$\tilde{s}^i_1,\dots,\tilde{s}^i_k$ is one choice of $(P_i,Q_i)$-adapted
basis, then in order to maintain the required vanishing orders at both $P_i$
and $Q_i$, any change of basis matrix is required to be upper and lower 
block triangular, with block sizes determined by 
$r^{i-1}_1,\dots,r^{i-1}_{\ell_{i-1}}$ and $r^i_1,\dots,r^i_{\ell_i}$. 
Considering first the fiber at $P_i$, we see that for any $W^{i-1}_j$, the 
allowable changes of basis for the images of 
$\tilde{s}^i_{R^{i-1}_j+1},\dots,\tilde{s}^i_{R^{i-1}_j+r^{i-1}_j}$
are block triangular. On the other hand, by the induction hypothesis 
we can also achieve upper triangular changes of basis for the images in
$W^{i-1}_j$ of
$s^{i-1}_{R^{i-1}_j+1},\dots,s^{i-1}_{R^{i-1}_j+r^{i-1}_j}$,
so by Lemma \ref{lem:flags-basis} below we conclude that there exists some
choices of bases $s^1,\dots,s^i$ and permutations 
$\sigma_1,\dots,\sigma_{i-1}$ as desired for the statement.
It remains to
prove the stronger inductive statement on the flexibility for restrictions 
to $Q_i$.

For each $j$, we thus describe an ordering $\alpha_j$ on
$s^i_{R^i_j+1},\dots,s^i_{R^i_j+r^i_j}$ with the desired property. First,
the ordering respects vanishing order at $P_i$, so that if
$a^i_m < a^i_{m'}$ for $m,m' \in [R^i_j+1,\dots,R^i_j+r^i_j]$, then
then $\alpha_j$ places $s^i_m$ ahead of $s^i_{m'}$. On the other hand,
if $a^i_m=a^i_{m'}$, then the ordering is determined by the appropriate
$\alpha'_{j'}$ and $\sigma_{i-1}$, as follows: since $a^i_m=a^i_{m'}$,
we have $m,m' \in [R^{i-1}_{j'}+1,\dots,R^{i-1}_{j'}+r^{i-1}_{j'}]$ for some
$j'$.  Then $\alpha_j$ places $s^i_m$ ahead of $s^i_{m'}$ if 
$s^{i-1}_{\sigma_{i-1}^{-1}(m)}$ is ahead of 
$s^{i-1}_{\sigma_{i-1}^{-1}(m')}$ under $\alpha'_{j'}$ (note that 
$\sigma_{i-1}$ necessarily preserves vanishing order at $Q_{i-1}$, so this
makes sense). It thus remains to verify that under the ordering $\alpha_j$, 
we can realize any upper triangular change of basis for the images of
$s^i_{R^i_j+1},\dots,s^i_{R^i_j+r^i_j}$ in $W^i_j$. The first observation
is that by construction, such changes of basis preserve orders of vanishing
at both $P_i$ and $Q_i$, so they certainly maintain the condition of
being a $(P_i,Q_i)$-adapted basis. Now, given such a change of basis,
we need to verify that we can also modify the choices of $s^1,\dots,s^{i-1}$
to maintain the gluing condition at each node. But again using the
construction of the ordering, for any $j'$ such that 
$R^{i-1}_{j'}+1,\dots,R^{i-1}_{j'}+r^{i-1}_{j'}$ overlaps with
$R^i_j+1,\dots,R^i_j+r^i_j$, a change of basis of
$s^i_{R^i_j+1},\dots,s^i_{R^i_j+r^i_j}$ which is upper triangular with 
respect to $\alpha_j$ induces in $W^{i-1}_{j'}$ a change of basis which
is upper triangular with respect to $\alpha'_{j'}$ (taking into account
also the reordering dictated by $\sigma_{i-1}$). Note here that although
the change of basis also involves $s^i_m$ with $m>R^{i-1}_{j'}+r^{i-1}_{j'}$,
this does not affect the image in $W^{i-1}_{j'}$ because the vanishing order
at $P_i$ is strictly greater. Thus, the induction hypothesis allows us
modify the $s^1,\dots,s^{i-1}$ to achieve the desired change of basis, and
we conclude the statement of the proposition. 
\end{proof}

The following lemma is surely a standard fact from linear algebra, but
for the sake of completeness we include a brief proof.

\begin{lem}\label{lem:flags-basis}
Let $V$ be a $d$-dimensional vector
space, and $U_0 \subseteq U_1 \subseteq \dots \subseteq U_d$ and
$W_0 \subseteq W_1 \subseteq \dots \subseteq W_d$ complete flags in
$V$. Then there exists a basis $v_1,\dots,v_d$ of $V$ such that every
$U_i$ and every $W_i$ is the span of some subset of the $v_j$.
\end{lem}

\begin{proof} The proof is by induction on $d$, with the case $d=1$
being trivial. Set $V'=U_{d-1}$. Then the chain of subspaces
$$W_0 \cap V' \subseteq W_1 \cap V' \subseteq \dots W_d \cap V'$$
still has every quotient with dimension at most $1$, and there must
be precisely one index $i_0 \leq d-1$ such that 
$W_{i_0} \cap V' = W_{i_0+1} \cap V'$.
Further, we have $W_i \subseteq V'$ if and only if $i \leq i_0$.
We thus obtain new complete flags $U_i'$ and $W'_i$ in $V'$ by setting
$U_i'=U_i$ for $i=0,\dots,d-1$, setting $W'_i=W_i \cap V'$ for $i=0,\dots,i_0$
and $W'_i=W_{i+1} \cap V'$ for $i=i_0+1,\dots,d-1$. According to the 
induction hypothesis, we obtain a basis $b'_1,\dots,b'_{d-1}$ of $V'$ 
such that every $W'_i$ and $U'_i$ is the span of some subset of the $b'_i$.
It follows that $U_i$ has the property for all $i<d$, as does $W'_i$ for
$i \leq i_0$. Finally, set $b_i=b'_i$ for $i=1,\dots,d-1$, and set $b_d$ to 
be any vector in $W_{i_0+1} \smallsetminus V'$; we see easily that the
$b_i$ form a basis, and that $W_i = \spn (W_i', b_d)$ for all $i>i_0$, giving
us the desired property.
\end{proof}

\begin{rem}
As mentioned in the proof, the conditions of the statement of Proposition 
\ref{prop:chain-adapted-basis} place strong limitations on
the permutations $\sigma_i$: specifically, they imply that if 
$\sigma_i(j) \neq j$, we must still have 
$\ord_{P_{i+1}} s^{i+1}_{\sigma(j)}=\ord_{P_{i+1}} s^{i+1}_j=b$.
Nonetheless, it is not difficult to find examples for which non-trivial
permutations are necessary in order to find bases of the desired form.

In the inductive statement in the proof, it is precisely due to the 
necessity of these permutations that we are required to introduce the 
orderings with respect to which we can take upper triangular changes of
basis. Thus, consideration of some such orderings is necessary. On the
other hand, we are in fact being unnecessarily restrictive by looking only 
at upper triangular matrices. The natural level of flexibility is expressed 
in terms of block upper triangular matrices, or equivalently compatibility 
with partial flags. However, by working with the coarser statement we
simplify notation considerably.
\end{rem}

\begin{cor}\label{cor:chain-adapted-constrained} 
A chain-adaptable EHT limit linear series on $X$ is constrained.
\end{cor}

\begin{proof} 
Let $((\sE^i,V^i)_{i=1,\dots,n},(\vp_i)_{i=1,\dots,n-1})$
be a chain-adaptable EHT limit linear series on $X$.
Then our main claim is that for any $w \in V(\bar{G}_{\II})$
such that $V^i(w)\neq 0$ for each $i$,
the kernel of \eqref{eq:gluing-ker2} has dimension exactly $k$.

For $i=1,\dots,n$, following Notation \ref{not:twisting-divisors},
write $D_{w,i}=a_i P_i + b_i Q_i$ for some $a_i,b_i$;
thus, we have by definition that $b_i+a_{i+1}=b$ for $i=1,\dots,n-1$.
Also, denote by $a^i$ the vanishing sequence of $V^i$ at $Q_i$,
so that the vanishing sequence of $V^{i+1}$ at $P_{i+1}$ is given by
$b-a^i_k,\dots,b-a^i_1$. For each $i$ from $1$ to
$n-1$, choose
$\ell_i$ minimal with $a^i_{\ell_i} \geq b_i$ and and $m_i$ maximal 
with $a^i_{m_i} \leq b_i$. 
Then we have
$$\dim V^1(w)=\#\{j:a^1_j \geq b_1\}=k+1-\ell_1,\text{ and}$$
$$\dim V^n(w)=\#\{j:b-a^{n-1}_j \geq a_n= b-b_{n-1}\}=m_{n-1}.$$
By Proposition \ref{prop:adapted}, for each $i$ with $1<i<n$, we have
\begin{align*} \dim V^i(w) 
& =\#\{j:b-a^{i-1}_j \geq a_i, a^i_j \geq b_i\} \\
& =\#\{j:a^{i-1}_j \leq b_{i-1},a^i_j \geq b_i\} \\
& =\max\{0,m_{i-1}+1-\ell_i\} \\
& =m_{i-1}+1-\ell_i,
\end{align*}
where the last equality holds because of the hypothesis that $V^i(w) \neq 0$.
Adding these up, we find
$$\sum_{i=1}^n \dim V^i(w)=k+\sum_{i=1}^{n-1}(m_i+1-\ell_i).$$
At the same time, at a given node $P_i$, we have that $V^i(w)$ and
$V^{i+1}(w)$ each have restriction to $P_i$ of dimension $m_i+1-\ell_i$,
so the rank of \eqref{eq:gluing-ker2} is at least (in fact, exactly)
$$\sum_{i=1}^{n-1}(m_i+1-\ell_i),$$
and we conclude that the kernel has dimension (at most, hence equal to)
$k$, as claimed.

Now, with notation as in the statement of Proposition 
\ref{prop:chain-adapted-basis}, we see that for $i=1,\dots,k$, taking 
$s^1_i,s^2_{\sigma_1(i)},s^3_{\sigma_2(\sigma_1(i))},\dots,
s^n_{\sigma_{n-1}(\dots \sigma_1(i))}$, we obtain a multidegree 
$w_i \in V(\bar{G}_{\II})$
and a vector $v_i$ in the kernel of the corresponding map 
\eqref{eq:gluing-ker2}.
Moreover, by definition for each $j=1,\dots,n$ the images of the $v_i$
in $V^j$
are simply a permutation of $s^j_1,\dots,s^j_k$, so form a basis for
$V^j$. We thus conclude that the given limit linear series is constrained.
\end{proof}

Combining Theorem \ref{thm:compare} with Corollary 
\ref{cor:chain-adapted-constrained}, we obtain the following comparison
result, which may be viewed as the main result of the present paper.

\begin{cor}\label{cor:compare} In Situation \ref{sit:basic-ii}, the 
forgetful morphism 
$$\cG^{k,\II}_{r,d,d_{\bullet}}(X,\theta_{\bullet}) \to \cG^{k,\EHT}_{r,d,d_{\bullet}}(X)$$
is an isomorphism on the preimage of the chain-adaptable locus.
\end{cor}

To illustrate how Corollary \ref{cor:compare} may be applied, if
we drop the hypothesis that $B$ is a point, and use 
Theorem \ref{thm:foundation}, we conclude the following smoothing result.

\begin{cor}\label{cor:chain-adapted-smoothing}
In Situation \ref{sit:basic-ii}, if $\Gamma$ is a chain, and for some fiber
$X_0$ the space $\cG^{k,\EHT}_{r,d,d_{\bullet}}(X_0)$ has dimension exactly 
$\rho-1$ 
at a chain-adaptable point $z$, then $\cG^{k,\II}_{r,d_{\bullet}}(X/B,\theta_{\bullet})$ 
is universally open over $B$ at $z$, and has pure 
fiber dimension $\rho-1$ in a neighborhood of $z$.
\end{cor}

The above result is uninteresting in and of itself, since it is weaker
than the corresponding result proved by Teixidor i Bigas in Theorem 2.6 of
\cite{te1}, using an adaptation of the construction of Eisenbud and
Harris. However, Corollary \ref{cor:compare} is the key ingredient in 
\cite{o-t2}, in which we 
extend Corollary \ref{cor:chain-adapted-smoothing} to the case of special 
determinants. By analyzing type II linked linear series in this context,
we prove new smoothing results for the modified expected 
dimensions arising in that case, and as a consequence obtain new results
on existence of vector bundles with sections on smooth curves. This
extended version of Corollary \ref{cor:chain-adapted-smoothing} should be
seen as the natural culmination of the machinery developed in the present
paper.

\section{Complementary results}\label{sec:complement}

In this section, we discuss several complementary results to our main
theorems, treating the fixed determinant case, as well as questions on
variation of the parameter $b$, specialization in one-dimensional families, 
and stability.

\subsection{The fixed determinant case}

We now describe the case of fixed determinant. Our results hold for
arbitrary determinant line bundle $\sL$, with everything the same as
in the varying determinant case except that the dimension $\rho-1$ is 
replaced by $\rho-g$. The dimension lower bounds 
are not optimal in the case that $\sL$ is special, but the
foundational results we develop here nonetheless play a fundamental role in
the generalization to special determinants presented in \cite{o-t2}.

We first fix notation for moduli spaces of vector bundles with specified 
determinant. It will often be convenient to specify the determinant in
a given predetermined multidegree (for instance, in the case of canonical
determinant), which may or may not agree with the parity possibilities for
the degrees of the vector bundles under consideration, so we allow the
determinant multidegree to differ from our fixed $w_0$. Also, for technical 
reasons, it is convenient to work with line bundles $\sL$ not necessarily 
defined over the original base $B$ (see Proposition 6.2 of \cite{o-t2}).
We will want to have separate notation for 
twisting line bundles in the rank-$1$ case, which we set as follows:

\begin{notn}\label{not:rk-1-twist} Given $w,w' \in \ZZ^{V(\Gamma)}$ both
lying in the affine hyperplane consisting of vectors summing to $d$, 
let $\sO'_{w,w'}$ be the line bundle defined as in Notation
\ref{not:twisting-bundles}, if applied to the case $r=1$.
\end{notn}

Note that when $r=1$, the hyperplane in Notation \ref{not:rk-1-twist}
is the same as $V(G_{\II})$. Although Notation \ref{not:twisting-bundles}
applies \textit{a priori} only to vertices in $V(G_{\I})$, in fact we
can between any $w$ and $w'$ by twisting with a sequence of bundles
$\sO_{(e,v)}$, and we express the construction this way in order to
avoid imposing the additional hypothesis required for the type II 
construction.

\begin{notn}\label{not:bundles-fixed-det}
Given a multidegree $w_0$ on $\Gamma$, a $B$-scheme $B'$, 
and a line bundle $\sL$ on $X_{B'}:=X \times_B B'$ having multidegree $w$ 
with total degree equal to that
of $w_0$, let $\cM_{r,w_0,\sL}(X_{B'}/B')$ be the moduli stack of vector
bundles $\sE$ on $X_{B'}$ of multidegree $w_0$, together with an 
isomorphism 
$$\left(\det \sE\right) \otimes \sO'_{w_0,w} \risom \sL.$$

Given also $d_{\bullet}$, let
$\cM_{r,w_0,\sL,d_{\bullet}}(X_{B'}/B')$ be the open substack of
$\cM_{r,w_0,\sL}(X_{B'}/B')$ for which the underlying vector bundles
lies in $\cM_{r,w_0,d_{\bullet}}(X/B)$.
\end{notn}

Thus, $\cM_{r,w_0,\sL}(X_{B'}/B')$ is the (2-)fibered product of 
$\cM_{r,w_0}(X/B)$ with $B'$ over $\cPic^{w}(X/B)$, with the 
first map given by taking determinant and twisting by $\sO'_{w_0,w}$,
and the second map being the one
induced by $\sL$, and similarly for $\cM_{r,w_0,\sL,d_{\bullet}}(X_{B'}/B')$ 
and $\cM_{r,w_0,d_{\bullet}}(X/B)$. These stacks are algebraic, and smooth 
over $B'$ of relative dimension $(r^2-1)(g-1)$.

We can thus make the following definitions.

\begin{defn}\label{def:fixed-det-stacks}
Given a line bundle $\sL$ of degree $d$ and multidegree $w$ on
$\Gamma$, as well as a $B$-scheme $B'$, we define
the stack $\cG^{k,\I}_{r,\sL,d_{\bullet}}(X/B)$ (respectively,
$\cG^{k,\II}_{r,\sL,d_{\bullet}}(X/B,\theta_{\bullet})$, 
$\cG^{k,\EHT}_{r,\sL,d_{\bullet}}(X/B)$) to be the (2-)fibered
product of $\cG^{k,\I}_{r,d,d_{\bullet}}(X/B)$ (respectively,
$\cG^{k,\II}_{r,d,d_{\bullet}}(X/B,\theta_{\bullet})$, $\cG^{k,\EHT}_{r,d,d_{\bullet}}(X/B)$)
with $\cM_{r,w_0,\sL}(X_{B'}/B')$ over
$\cM_{r,w_0}(X/B)$.

We further define the refined and constrained loci of
$\cG^{k,\EHT}_{r,\sL,d_{\bullet}}(X/B)$ to be the preimages of the 
corresponding loci of $\cG^{k,\EHT}_{r,d,d_{\bullet}}(X/B)$, and similarly
with the simple locus of $\cG^{k,\II}_{r,\sL,d_{\bullet}}(X/B,\theta_{\bullet})$.
\end{defn}

Thus, each stack is as in the varying determinant case, except that
we add an isomorphism between $\sL$ and the appropriate twist of the
determinant of the underlying vector bundle. This isomorphism not only
specifies the determinant, but rigidifies the groupoid, which is why
the dimension in the fixed determinant case agrees with the dimension
one obtains with coarse moduli spaces.

Note that each of these groupoids lives naturally over $B'$ rather than $B$,
but that the family $X/B$ is used in the definition in order to define the
necessary twisting line bundles, because $X_{B'}$ may not be regular.

As before, we have two main theorems on these spaces.

\begin{thm}\label{thm:foundation-det} Let $X/B$ be an almost local smoothing 
family, $k,r,d,d_{\bullet}$ as in Situation \ref{sit:basic}, and
$\sL$ a line bundle of degree $d$ on $X _B'$ for
some $B$-scheme $B'$.
Then $\cG^{k,\I}_{r,\sL,d_{\bullet}}(X/B)$ is an Artin stack over $B'$,
and the natural map
$$\cG^{k,\I}_{r,\sL,d_{\bullet}}(X/B) \to 
\cM_{r,w_0,\sL,d_{\bullet}}(X_{B'}/B')$$
is relatively representable by schemes which are projective, at least
locally on the target. Moreover, formation of 
$\cG^{k,\I}_{r,\sL,d_{\bullet}}(X/B)$
is compatible with any base change $B'' \to B$ which preserves the almost 
local smoothing family hypotheses.
In particular, if $y \in B'$ is a point with $X_y$ smooth, then
the base change to $y$ parametrizes triples $(\sE,\psi,V)$ of
a vector bundle $\sE$ of rank $r$ and degree $d$ on $X_y$ together with
an isomorphism $\psi:\det \sE \risom \sL|_y$
and a $k$-dimensional vector 
space $V \subseteq H^0(X_y,\sE)$.

Under the further hypothesis of Situation \ref{sit:basic-ii}, all of the
above statements also hold for $\cG^{k,\II}_{r,\sL,d_{\bullet}}(X/B,\theta_{\bullet})$.
Moreover, the simple locus of $\cG^{k,\II}_{r,\sL,d_{\bullet}}(X/B,\theta_{\bullet})$ has 
universal relative dimension at least 
$k(d-k-r(g-1))$ over $\cM_{r,w_0,\sL,d_{\bullet}}(X/B)$, and therefore
universal relative dimension at least $\rho-g$ over $B'$.
In particular, if some fiber
$\cG^{k,\II}_{r,\sL,d_{\bullet}}(X_y/y,\theta_{\bullet})$ has dimension 
exactly $\rho-g$ at a simple point $z$, then 
$\cG^{k,\II}_{r,\sL,d_{\bullet}}(X/B,\theta_{\bullet})$ is universally
open at $z$, and has fibers of pure dimension $\rho-g$ in an open 
neighborhood of $z$. 
\end{thm}

\begin{rem} We also see immediately that in the situation of
Theorem \ref{thm:foundation-det}, in the case that $B$ is a point,
we have that $\cG^{k,\EHT}_{r,\sL,d_{\bullet}}(X/B)$ is an Artin stack over 
$B'$, and the natural map
$$\cG^{k,\EHT}_{r,\sL,d_{\bullet}}(X/B) \to 
\cM_{r,w_0,\sL,d_{\bullet}}(X_{B'}/B')$$
is relatively representable by schemes which are projective, at least
locally on the target.
\end{rem}

\begin{thm}\label{thm:compare-det} Suppose that $B$ is a point.
Then the morphism \eqref{eq:forget} induces
a surjective morphism
\begin{equation}\label{eq:forget-det}
\cG^{k,\I}_{r,\sL,d_{\bullet}}(X/B) \to \cG^{k,\EHT}_{r,\sL,d_{\bullet}}(X/B),
\end{equation}
which is an isomorphism on the preimage of the refined locus.

In addition, under the further hypothesis of Situation \ref{sit:basic-ii},
if we consider the composed forgetful morphism
\begin{equation}\label{eq:forget-comp-det} 
\cG^{k,\II}_{r,\sL,d_{\bullet}}(X,\theta_{\bullet}) \to
\cG^{k,\EHT}_{r,\sL,d_{\bullet}}(X),
\end{equation}
we find that it is an isomorphism on the preimage of the constrained locus.
\end{thm}

Theorem \ref{thm:compare-det} follows immediately from Theorem 
\ref{thm:compare}, since the relevant maps are simply obtained from base
change of those in the latter result.

\begin{proof}[Proof of Theorem \ref{thm:foundation-det}] That the groupoids
are Artin stacks follows immediately from Theorem \ref{thm:foundation}
and the description as fibered products. Similarly, the map 
$$\cG^{k,\I}_{r,\sL,d_{\bullet}}(X/B) \to 
\cM_{r,w_0,\sL,d_{\bullet}}(X_{B'}/B')$$
is obtained from the map
$$\cG^{k,\I}_{r,d_{\bullet}}(X/B) \to \cM_{r,w_0,d_{\bullet}}(X/B)$$
by base change, so its projectivity is preserved, and likewise
for $\cG^{k,\II}_{r,\sL,d_{\bullet}}(X/B,\theta_{\bullet})$. Finally, since universal
relative dimension is preserved under base change, 
the dimension statement likewise follows immediately from
Theorem \ref{thm:foundation}.
\end{proof}

\subsection{The parameter $b$}\label{sec:b}
In the classical rank-$1$ case, Eisenbud and Harris made use of the facts
that twisting acts transitively on multidegrees, and that line bundles of
negative degree on a smooth curves have no nonzero global sections, in
order to always work with linear series of degree $d$ on each component.
Both of these facts fail in higher rank, which is the reason for the
introduction of the parameters $b$ and $d_{\bullet}$ in Teixidor's 
work \cite{te1}. The $d_{\bullet}$ keeps track of the extremal degrees
(and implicitly, the congruence class of allowable multidegrees), while
$b$ measures how much we twist to get from one extremal degree to another.
Because of condition (I) in Definitions \ref{def:grd-space-i}, 
\ref{def:grd-space-ii}, and \ref{def:eht-lls}, it is possible for $b$ to
be ``too small,'' but we can always increase $b$ and $d_{\bullet}$, and
if we do so in certain ways, we obtain an open immersion from one stack
of limit linear series to the other.

\begin{prop}\label{prop:increase-b} Suppose we have $r,d,k,b,d_{\bullet}$
as in Situation \ref{sit:basic}, and we are also given $b'$ and $d'_{\bullet}$
satisfying \eqref{eq:deg-sum}, and such that for each $v \in V(\Gamma)$,
we have
$$d_v \equiv d'_v \pmod{r}$$
and
$$d_v \leq d'_v \leq d_v+r(b'-b).$$

Under the hypothesis of Situation \ref{sit:basic-ii}, we have that
$\cG^{k,\II}_{r,\sL,d_{\bullet}}(X/B,\theta_{\bullet})$ is an open substack of
$\cG^{k,\II}_{r,\sL,d'_{\bullet}}(X/B,\theta_{\bullet})$.
\end{prop}

Similar statements appear to hold for the type I and EHT cases, but
the constructions are more involved, so we do not pursue them.

\begin{proof} The first observation is that our restrictions on
$d'_{\bullet}$ and $b'$ imply that we can go from $d_{\bullet},b$ to
$d'_{\bullet}, b'$ by repeatedly carrying out the following operation:
for a fixed vertex $v_0$, increase $d_v$ by $r$ for all $v \neq v_0$,
and increase $b$ by $1$. Indeed, there exists $v_0$ such that 
$d'_{v_0}<d_{v_0}+r(b'-b)$, and for such $v_0$, it follows that
$d'_v>d_v$ for all $v \neq v_0$. If we modify $d_{\bullet}$ and $b$ as
described above, all our restrictions are still satisfied, so we can
iterate until we reach $d'_{\bullet},b'$. Thus, it is enough to consider
the situation that, for a fixed $v_0$, we have 
$$d'_v=\begin{cases} d_v: & v=v_0 \\ d_v+r: & v \neq v_0.\end{cases}$$
and $b'=b+1$.

Now, the restriction on congruence 
classes means
that the two versions of $G_{\II}$ arising from $d_{\bullet}$ and
$d'_{\bullet}$ are canonically identified.
For each $v \in V(\Gamma)$, denote by $w'_v \in V(G_{\II})$
the vertex equal to $d'_v$ in index $v$ and to $d'_{v'}-rb$ in index $v'$
for all $v' \neq v$. Since condition (I) of Definition \ref{def:grd-space-ii}
is an open condition, it is enough to check that the version imposed by
$d_{\bullet}$ is stronger than that imposed by $d'_{\bullet}$. However,
for each $v$ we have that 
$\sE_{w'_v}|_{Y_v}$ is obtained from $\sE_{w_v}|_{Y_v}$ by twisting up
by an effective divisor $D_v$ supported at nodes,
and that the degree of $D_v$ is at most $b'-b=1$. It thus follows
immediately that condition (I) for $d_{\bullet}$ is stricter than
condition (I) for $d'_{\bullet}$, as desired.
\end{proof}

The next proposition shows that even if we do not assume condition (I)
in the definition of type II linked linear series, we can always increase 
$b$ so that it will be satisfied.

\begin{prop}\label{prop:b-exists} In the situation of Definition
\ref{def:grd-space-i}, suppose that $S$ is a quasicompact $B$-scheme,
and $\sE$ is a vector bundle of rank $r$ and multidegree $w_0$ on
$X \times_B S$. Then there exist $b',d'_{\bullet}$ as in
Proposition \ref{prop:increase-b} such that $\sE$ is in 
$\cM_{r,w_0,d'_{\bullet}}(X/B)$.
\end{prop}

\begin{proof} By the quasicompactness of $S$, there exists some $N$ such
that for every $v \in V(\Gamma)$ and every $z \in S$, with image $y \in B$,
there is no nonzero map from any line bundle of
degree $N$ on $Y$ to $\sE_{w_v}|_{Y}$,
where $Y$ denotes the component of the fiber $X_z$ corresponding to 
$\cl_y(v)$. 
Equivalently, if $D$ is any divisor of degree at
least $N$ on $Y$, then $H^0(Y,\sE_{w_v}|_Y(-D))=0$. Note that if
$b+1 \geq N$, condition (I) is automatically satisfied for $\sE$. Otherwise,
if we set $b'=b+|V(\Gamma)|(N-b-1)$, and $d'_v=d_v+r|E(\Gamma)|(N-b-1)$ for
all $v \in V(\Gamma)$, then we have $\sE_{w'_v}|_Y=\sE_{w_v}|_Y(D)$ for
a divisor $D$ of degree $|E(\Gamma)|(N-b-1)$, and then
$\sE_{w_v}|_Y(D-(b'+1)\Delta')$ has no nonzero global sections, because
$$\deg \left((b'+1)\Delta'-D\right)=b+1+(|V(\Gamma)|-|E(\Gamma)|)(N-b-1)=N.$$
\end{proof}

\begin{rem}\label{rem:universal-type-ii}
The significance of Proposition \ref{prop:b-exists} is that we could
have omitted $b$ and $d_{\bullet}$, as well as condition (I), from our 
definition of type II linked linear series, and obtained a single 
canonical moduli stack for each choice of congruence class modulo $r$ of
$d_{\bullet}$. However, we have chosen to use the present definition in
order to unify the presentation as much as possible with the type I case
and the EHT case.
\end{rem}

\subsection{Specialization}
Since the moduli stack of vector bundles on a reducible curve is not
proper, we do not automatically obtain specialization results when we
have a smoothing family with one-dimensional base. However, just as was
the case with Eisenbud and Harris, and with Teixidor's generalization
thereof, we have that after base change and blowup we can always
extend to a higher-rank limit linear series. In fact, our situation is
better than these cases, because after base change and blow up, our 
extension always give a point of a relative moduli stack over the whole
family, and this can be done regardless of characteristic. Previously,
in the higher-rank case this had not been carried out in any characteristic,
and in the rank-$1$ case it had only been done in characteristic $0$.

Our result is the following:

\begin{prop}\label{prop:specialization} Suppose that $B$ is one-dimensional,
let $y \in B$ be a point, set $U:=B \smallsetminus \{y\}$, and let 
$\pi:X \to B$ be a smoothing family such that $\pi$ is smooth on the
preimage of $U$. Given $r,d,k$, suppose we have $B'$ over $B$, also
one-dimensional and regular, and with a unique point $z$ lying over $y$;
set $U'=B'\smallsetminus \{z\}$. Let $\pi':X'\to B'$ be the 
desingularization of $X \times_B B' \to B'$ obtained by blowing up the 
nodes over $z$ as necessary.
Then given a pair $(\sE,\sV)$ on $X|_{U'}$, where $\sE$ is a vector bundle 
of degree $d$, and $\sV \subseteq \pi_* \sE$ is a rank-$k$ subbundle,
there exist $d_{\bullet}$ and $b$ as in Situation \ref{sit:basic} 
such that the $U'$-valued point of 
$\cG^{k,\II}_{r,d_{\bullet}}(X'/B',\theta_{\bullet})$ determined by 
$(\sE,\sV)$ extends to a $B'$-valued point, and similarly for
$\cG^{k,\I}_{r,d_{\bullet}}(X'/B')$ and $\cG^{k,\EHT}_{r,d_{\bullet}}(X'/B')$.
\end{prop}

\begin{proof} Since we have morphisms
$\cG^{k,\II}_{r,d_{\bullet}}(X'/B',\theta_{\bullet}) \to \cG^{k,\I}_{r,d_{\bullet}}(X'/B')$ 
and 
$\cG^{k,\I}_{r,d_{\bullet}}(X'/B')\to \cG^{k,\EHT}_{r,d_{\bullet}}(X'/B')$
which are isomorphisms wherever $\pi'$ is smooth, it is enough to prove 
the claimed statement for 
$\cG^{k,\II}_{r,d_{\bullet}}(X'/B',\theta_{\bullet})$. Let $\Gamma'$ be the 
graph associated to $\pi'$.
Now, since vector bundles extend over regular surfaces, 
there exists some $w \in \ZZ^{V(\Gamma')}$ such that $\sE$ extends to
a vector bundle $\sE_w$ of multidegree $w$ on $X'$. Next, by
Proposition \ref{prop:b-exists} we see that
there exist $b$ and $d_{\bullet}$ such that 
$\sE_w$ satisfies condition (I) of Definition \ref{def:grd-space-ii}.
Now, because 
$\cG^{k,\II}_{r,d_{\bullet}}(X'/B',\theta_{\bullet})$ is proper over 
$\cM_{r,w_0,d_{\bullet}}(X'/B')$,
we conclude that we obtain the desired $B'$-valued point
of $\cG^{k,\II}_{r,d_{\bullet}}(X'/B',\theta_{\bullet})$.
\end{proof}

\subsection{Stability}
We now address conditions imposed by stability of underlying vector
bundles. There is substantial possibility for confusion, due in part to the 
fact that unlike the irreducible case, on reducible curves stability
of vector bundles is not preserved by twisting by line bundles. In fact,
in addition to the usual notion of stability on a reducible curve, which
is most useful in specialization arguments, we will use a second notion --
called \el-stability -- which is more powerful in the context of smoothing
arguments.

First, recall
that on a reducible curve, stability is not in general a canonical notion,
but rather depends on a polarization of the curve, as follows:

\begin{defn}\label{def:polarization}
Let $X_0$ be a proper nodal curve with dual graph $\Gamma$;
a \textbf{polarization} on $X_0$ consists of
a weight function $\omega:V(\Gamma) \to \QQ$ such that $\omega(v)>0$ for
all $v$, and 
$\sum_{v \in V(\Gamma)} \omega(v)=1$. Given such a polarization, a vector 
bundle $\sE$ on $X_0$ is \textbf{semistable}, if
for all nonzero subsheaves $\sF \subseteq \sE$, we have
$$\frac{\chi(\sF)}{\sum_v \omega(v) \rk \sF|_{Y_v}} 
\leq \frac{\chi(\sE)}{\rk \sE}.$$
We say $\sE$ is \textbf{stable} if we always have strict inequality above.
\end{defn}

For an irreducible curve, this is equivalent to the usual definition in
terms of degree. We see that the role of the polarization is that it
determines what replaces the notion of rank for a sheaf which has 
different rank on different components of $X_0$.

We also have the following concept, introduced in \cite{os22}.

\begin{defn}\label{def:l-stable}
Let $X_0$ be a proper nodal curve with dual graph $\Gamma$.
A vector bundle $\sE$ on $X_0$ is $\boldsymbol{\ell}$-\textbf{semistable}
if for all nonzero subsheaves $\sF \subseteq \sE$ having constant rank on
all component of $X_0$, we have
$$\frac{\chi(\sF)}{\rk \sF} 
\leq \frac{\chi(\sE)}{\rk \sE}.$$
We say $\sE$ is $\boldsymbol{\ell}$-\textbf{stable} if we always have 
strict inequality above.
\end{defn}

Thus, except in the irreducible case, \el-stability is visibly a weaker 
condition than stability. However, it is independent of polarization, and 
is also more robust, being stable under twists and gluing; see Propositions
1.3 and 1.6 of \cite{os22}. Most importantly, it is open, and thus suffices 
for existence arguments using limit linear series.

We make the following definitions for (semi)stability of a limit linear series.

\begin{defn}\label{def:stable-lls}
Let $\pi:X \to B$ be an almost-local smoothing family with
associated graph $\Gamma$, and let $\cG$ be any of
$\cG^{k,\I}_{r,d_{\bullet}}(X/B)$, $\cG^{k,\II}_{r,d_{\bullet}}(X/B,\theta_{\bullet})$,
or $\cG^{k,\EHT}_{r,d_{\bullet}}(X/B)$.
Then a $K$-valued point of $\cG$ with image $y$ in $B$
is \textbf{\el-semistable} (respectively, \textbf{\el-stable}) 
if the vector bundle $\sE_{w_0}$ is \el-semistable (respectively, \el-stable).
Similarly, if $\omega$ is a polarization on $\Gamma$,
a $K$-valued point of $\cG$ with image $y$ in $B$
is \textbf{semistable} (respectively, \textbf{stable}) with respect to $w$ 
if there exists some $w_1 \in V(G_{\II})$ such that the induced vector 
bundle $\sE_{w_1}$ is semistable (respectively, stable) with respect to 
$\omega_y$, where $\omega_y$ is the polarization on the dual graph
$\Gamma_y$ induced by $\omega$ and $\cl_y$.

We make the same definitions for the case of fixed determinant.
\end{defn}

Note that, unlike the case of (semi)stability with respect to a polarization,
because \el-(semi)stability for vector bundles is invariant under twisting, 
it would be equivalent to define it above in terms of existence of $w_1$ 
such that $\sE_{w_1}$ is \el-(semi)stable.
 
We first give a precise statement of the behavior of stability
under smoothing.

\begin{prop}\label{prop:stable-open} Let $\pi:X \to B$ be an almost-local
smoothing family. Then \el-stability
and \el-semistability are both open conditions on each 
of $\cG^{k,\I}_{r,d_{\bullet}}(X/B)$, 
$\cG^{k,\II}_{r,d_{\bullet}}(X/B,\theta_{\bullet})$,
or $\cG^{k,\EHT}_{r,d_{\bullet}}(X/B)$, and the same holds in the case
of fixed determinant.
\end{prop}

This is an immediate consequence of Proposition 1.2 of \cite{os22}.
Thus, if we produce \el-stable limit linear series, under smoothing we
will obtain \el-stable (hence stable) vector bundles on nearby smooth 
curves.

Next, we consider behavior under specialization.

\begin{prop}\label{prop:stable-specialize} In the situation of
Proposition \ref{prop:specialization}, suppose further that $\sE$ is
semistable over the generic point of $U'$.
Then given any polarization $\omega$ on the dual graph $\Gamma'$ 
of the special fiber of $X'$,
we may choose the extension given by
Proposition \ref{prop:specialization} to be semistable
with respect to $\omega$.
\end{prop}

This is immediate from Proposition 4.1 of \cite{os22}.

Here, the flexibility
in choice of of multidegree in the definition of semistability for limit
linear series weakens the
conclusion of Proposition \ref{prop:stable-specialize}, in that it gives
us only a single multidegree on which we have semistability. On the other
hand, it is easy to see that this is the most we can hope for, since 
semistability on reducible curves imposes constraints on multidegrees
(see for instance (\texttt{*}) in the proof of Proposition 2.1 of \cite{te8}).

\appendix
\section{Prelinked Grassmannians}\label{app:pre-lg}

In this appendix, we develop a partial generalization of the linked
Grassmannians introduced in Appendix A of \cite{os8} for graphs consisting
of chains of vertices. Linked Grassmannians act as degenerations of 
Grassmannians, and we expect that the full theory generalizes to graphs
of the sort appearing for type I and II linked linear series, and indeed
more generally. However, such a theory will be difficult, and does not
yield immediate applications in the higher-rank case, so for the present
we consider instead a notion of ``prelinked Grassmannians,'' developing
easy properties which suffice for our purposes.

\subsection{Definitions}

We begin with the definitions underlying prelinked Grassmannians.

\begin{sit}\label{sit:lg-basic} Let $G$ be a finite directed graph, connected
by (directed) paths, $d$ an integer, 
and $S$ a scheme. Suppose we are given data $\sE_{\bullet}$, consisting of vector 
bundles $\sE_v$ of rank $d$ for each vertex $v \in V(G)$, and morphisms 
$f_e:\sE_v \to \sE_{v'}$ for each edge $e \in E(G)$, where $v$ and $v'$
are the tail and head of $e$, respectively. For a (directed) path $P$ in
$G$, denote by $f_P$ the composition of the morphisms $f_e$ for each
edge $e$ in $P$.
\end{sit}

In the above situation, we have:

\begin{defn}\label{def:pre-lg} Given an integer $r<d$, suppose further that 
the following condition is satisfied:
\begin{Ilist}
\itm For any two paths $P,P'$ in $G$ with the same head and tail, there
are scalars $s,s' \in \Gamma(S,\sO_S)$ such that
$$s \cdot f_P =s' \cdot f_{P'};$$
moreover, if $P$ (respectively, $P'$) is minimal, then $s'$ 
(respectively, $s$) is invertible.
\end{Ilist}
Then define the \textbf{prelinked Grassmannian} $\LG(r,\sE_{\bullet})$ to be the
scheme representing the functor associating to an $S$-scheme $T$ the
set of all collections $(\sF_v)_{v \in V(G)}$ of rank-$r$ subbundles 
of the $\sE_v$ satisfying the property that for all edges $e \in E(G)$,
$f_e(\sF_v) \subseteq \sF_{v'}$, where $v$ and $v'$ are the tail and
head of $e$ respectively.
\end{defn}

Note that in the condition of the definition, we allow the empty path and
consider $f_{P}=\id$ in this case. Thus, the condition implies that for
any loop $P$ in $G$, we have $f_P =s \cdot \id$ for some scalar $s$.

The fact that the functor defining $\LG(r,\sE_{\bullet})$ is representable is
easy, and in fact we have the following.

\begin{prop}\label{prop:pre-lg} $\LG(r,\sE_{\bullet})$ is a projective scheme over
$S$, and compatible with base change.
\end{prop}

\begin{proof} We construct $\LG(r,\sE_{\bullet})$ as a closed subscheme of the
natural product of Grassmannians $G:= \prod_{v \in V(G)} G(r,\sE_v)$. 
For each $v$, let $\widetilde{\sE}_v$ denote the pullback of $\sE_v$ to $G$, 
and let $\widetilde{\sF}_v$ denote the pullback of the universal subbundle 
on $G(r,\sE_v)$ to $G$. Then it suffices to observe that $\LG(r,\sE_{\bullet})$
is cut out as the locus on which, for every edge $e \in E(G)$, if $v$
and $v'$ are the tail and head of $e$, then the composed map of vector 
bundles
$$\widetilde{\sF}_v \to \widetilde{\sE}_{v'} \to 
\widetilde{\sE}_{v'}/\widetilde{\sF}_{v'}$$
is zero. 

That $\LG(r,\sE_{\bullet})$ is compatible with base change is evident from the
definition.
\end{proof}

In this generality, we will mainly be interested in points satisfying
a strong additional condition:

\begin{defn}\label{def:simple-pt} Let $K$ be a field over $S$, and $(F_v)_v$ 
a $K$-valued point of a prelinked Grassmannian $\LG(r,\sE_{\bullet})$.
We say that $(F_v)_v$ is \textbf{simple} if
there exist $v_1,\dots,v_r \in V(G)$ (not necessarily distinct) and
$s_i \in F_{v_i}$ for $i=1,\dots,r$ such that for every $v \in V(G)$,
there exist paths $P^{v_1},\dots,P^{v_g}$ with each $P^{v_i}$ going from 
$v_i$ to $v$, and such that $f_{P^{v_1}}(s_1),\dots,f_{P^{v_r}}(s_r)$ form 
a basis for $F_v$.
\end{defn}

Note that as a consequence of the definition of a prelinked Grassmannian,
we may always take the $P^{v_i}$ in the definition of a simple point to be
minimal paths from $v_i$ to $v$.

\subsection{Properties of simple points}

In \cite{os8}, a notion of ``exact points'' of linked Grassmannians 
constitutes the focal point of the analysis. In the special case studied in 
\textit{loc.\ cit.}, exact points are the same as simple points, which 
substantially simplifies the situation. Exact points generalize naturally to 
the cases we are interested in, but are no longer the same as simple points;
see Example \ref{ex:exact-not-simple}. We will now prove that simple points
are smooth of the expected dimension, which is the relevant statement for
our applications, and is much easier than analyzing exact points in general.

First, using Nakayama's lemma we easily conclude the following:

\begin{prop}\label{prop:simple-open} The simple points form an open subset 
of $\LG(r,\sE_{\bullet})$. 
\end{prop}

More substantively, we have the following result, which can be viewed as
a generalization of Lemma A.12 (ii) and Lemma A.14 of \cite{os8}.

\begin{prop}\label{prop:simple-smooth} On the locus of simple points,
$\LG(r,\sE_{\bullet})$ is smooth over $S$ of relative dimension $r(d-r)$.
\end{prop}

\begin{proof} We first prove smoothness. Let $(F_v)_v$ be a $K$-valued 
simple point, and fix choices of $v_i$, $s_i$ and $P^{v_i}$ as in Definition 
\ref{def:simple-pt}, with each $P^{v_i}$ a minimal path.
Let $A$ be a Noetherian local ring over $S$ with residue field $K$, and 
$A'$ a quotient ring of $A$. Suppose we are given $(\sF'_v)_v$ over $A'$
specializing to the given point under restriction to $K$; we wish to
show that we can lift to $(\sF_v)_v$ over $A$. To keep notation reasonable,
for any edge $e$ or path $P$ in $G$ we will still use the notation $f_e$ or
$f_P$ for all base changes of the original morphisms $f_e$ and $f_P$ on $S$.
For each $i$, let 
$\tilde{s}_i'$ be any lift of $s_i$ in $\sF'_{v_i}$; then by Nakayama's
lemma, we have that for each $v$, the images of $\tilde{s}_i'$ under
the corresponding
$f_{P^{v_i}}$ form a basis for $\sF'_v$. 

We claim that if we let
$\tilde{s}_i$ be any lift of $\tilde{s}_i'$ to $\sE_{v_i}$ for $i=1,\dots,r$,
then the $\tilde{s}_i$ uniquely determine an $A$-valued point of
$\LG(r,\sE_{\bullet})$ extending the given $A'$-valued one. Indeed, for each 
$v \in V(G)$, we have the map
$$A^{\oplus r} \to \sE_v$$
determined by $(a_1,\dots,a_r) \mapsto \sum_i a_i f_{P^{v_i}}(\tilde{s}_i)$.
By hypothesis, this morphism has full rank at the closed point, and it
follows that its image is a subbundle of $\sE_v$, which will be our $\sF_v$.
It remains to check that the $\sF_v$ thus determined are in fact linked by
$f_e$ for each edge $e \in E(G)$. Accordingly, let $e$ be an edge, with
tail $v$ and head $v'$. It suffices to observe that for each $i$, there is 
some $s \in A$ such that
$f_e (f_{P^{v_i}}(\tilde{s}_i)) = s \cdot f_{P^{v'}_i}(\tilde{s}_i)$,
because we have chosen $P^{v'}_i$ to be a minimal path. This proves the
desired lifting statement, and we conclude smoothness.

To conclude the proof of the proposition, we carry out a tangent space 
computation in the fiber. Since dimension isn't affected by extension of
base field, we may even assume that $S=\Spec K$ (where we are still 
considering the above $K$-valued simple point). Our claim is that the 
tangent space is (given the choice of $v_i,s_i$ as above) canonically 
identified with $\bigoplus_{i=1}^r E_{v_i}/F_{v_i}$, which has dimension
$r(d-r)$, as desired. We first construct a map from 
$\bigoplus_{i=1}^r E_{v_i}$ to the tangent space, as follows: setting
$A=K[\epsilon]/\epsilon^2$, and $A'=K$, and applying our above analysis, 
we see that any choices of $\tilde{s}_i$ determine a tangent vector to
$\LG(e,E_{\bullet})$ at the given point, and conversely, for any tangent vector
$(\sF_v)_v$, of course there exists some choice of lifts of the $\tilde{s}_i$ 
contained in the given $\sF_{v_i}$, and thereby inducing $(\sF_v)_v$.
Now, the ambient spaces $\sE_v$ over $A$ are by definition equal to 
$E_v \oplus \epsilon E_v$, so a lift of $s_i$ is given uniquely by an
element of $E_{v_i}$. This gives the map from 
$\bigoplus_{i=1}^r E_{v_i}$ to the tangent space which we have seen to
be surjective,
so it remains to check that the kernel is precisely given by
$\bigoplus_{i=1}^r F_{v_i}$. It is clear that if $s_i$ is lifted to
$s_i+\epsilon w$, with $w \in F_{v_i}$, the resulting spaces $\sF_v$ are
the same as if $s_i$ were lifted to $s_i$.
Thus, the given space is contained in the kernel. Conversely, if given
lifts $s_i+\epsilon w_i$ of the $s_i$ yield the trivial deformations
$F_v \oplus \epsilon F_v$, then in particular we see that we must have
$w_i \in F_{v_i}$ for all $i$, so the kernel is as asserted, and we
conclude the tangent space has dimension $r(d-r)$.
\end{proof}

\begin{ex}\label{ex:exact-not-simple} Let $G$ be a graph consisting of
four vertices $v_1,v_2,v_3,v_4$ with edges $e_{i,j}$ from $v_i$ to $v_j$
whenever $i=1$ or $j=1$. We consider vector spaces $E_{v_i}=K^3$ for all 
$i$, and maps $f_{e_{i,j}}$ as follows:
$$f_{e_{1,2}}=
\begin{bmatrix} 0 & 0 & 0 \\ 0 & 0 & 0 \\ 0 & 0 & 1 \end{bmatrix},
\quad 
f_{e_{2,1}}=
\begin{bmatrix} 1 & 0 & 0 \\ 0 & 1 & 0 \\ 0 & 0 & 0 \end{bmatrix};$$
$$f_{e_{1,3}}=
\begin{bmatrix} 0 & 0 & 0 \\ 0 & 1 & 0 \\ 0 & 0 & 0 \end{bmatrix},
\quad 
f_{e_{3,1}}=
\begin{bmatrix} 1 & 0 & 0 \\ 0 & 0 & 0 \\ 0 & 0 & 1 \end{bmatrix};$$
$$f_{e_{1,4}}=
\begin{bmatrix} 1 & 0 & 0 \\ 0 & 0 & 0 \\ 0 & 0 & 0 \end{bmatrix},
\quad 
f_{e_{4,1}}=
\begin{bmatrix} 0 & 0 & 0 \\ 0 & 1 & 0 \\ 0 & 0 & 1 \end{bmatrix}.$$
Then we set subspaces
$$V_1=\left<(1,1,0),(1,0,1)\right>, \quad
V_2=\left<(1,1,0),(0,0,1)\right>,$$
$$V_3=\left<(0,1,0),(1,0,1)\right>, \quad
V_4=\left<(1,0,0),(0,1,-1)\right>.$$
One sees that this is an exact point, but not simple. Indeed, the
images of $V_2$, $V_3$ and $V_4$ in $V_1$ are three distinct lines,
and this can never happen at a simple point.

Although this graph is not explicitly of the form arising from a 
type II linked linear series, it can easily be extended to one. 
Indeed, we can see the same behavior arising directly from a linked
$\fg^1_1$ (of rank $1$) on a curve with four rational components, having 
one main component glued to each of the other three. Then the image in
multidegree $(1,0,0,0)$ of the spaces in multidegrees $(0,1,0,0)$, $(0,0,1,0)$
and $(0,0,0,1)$ consists of sections vanishing at the first, second or
third node respectively, so again gives three different lines in the 
$2$-dimensional space of global sections.
\end{ex}

\subsection{Prelinked Grassmannians over stacks}\label{sec:lg-stacks}

Our results on prelinked Grassmannians can be generalized in a routine manner
to stacks. In order to keep the situation more geometric, we restrict to
algebraic stacks, although this is not strictly necessary. Indeed, 
Definition \ref{def:pre-lg} generalizes immediately to the case that $S$
is replaced with a stack $\cS$, yielding a groupoid $\cLG(r,\sE_{\bullet})$ over
$\cS$. Compatibility with base change in Proposition \ref{prop:pre-lg}
implies that $\cLG(r,\sE_{\bullet})$ is also a stack, relatively representable over
$\cS$ by projective schemes. The definition of simple point also generalizes
to this context, and Proposition \ref{prop:simple-smooth} lets us conclude
the following:

\begin{cor}\label{cor:lg-stacks} Let $\cS$ be an algebraic stack, and
$\cLG(r,\sE_{\bullet})$ a prelinked Grassmannian over $\cS$. Then $\cLG(r,\sE_{\bullet})$
is also an algebraic stack, relatively representable by projective schemes
over $\cS$, and furthermore, on the locus of simple points,
$\cLG(r,\sE_{\bullet})$ is smooth over $\cS$ of relative dimension $r(d-r)$.
\end{cor}

\section{Generalized determinantal loci}\label{app:pushforward-detl}

In this appendix, we study properties of classical determinantal loci,
and apply the results to generalize the definition to a relative
setting. Our point of view is that frequently in applications to moduli
spaces, we are interested not in the rank of a given map, but in the size
of its kernel. We explore this point of view first in the context of
classical determinantal ideals, and then we apply our results to develop
a theory of determinal loci of pushforwards.

\subsection{Observations on determinantal ideals}

We begin with some results in the context of classical determinantal loci.

\begin{defn}\label{def:vanishing-locus}
Given locally free sheaves $\sE,\sF$ of finite rank on a scheme $B$, a
morphism $f: \sE \to \sF$, and an integer $k\geq 0$, define the 
\textbf{$k$th vanishing locus} $V_k(f)$ of $f$ to be the closed subscheme 
of $B$ cut out by the vanishing of 
$$\bigwedge^{r+1-k} f:\bigwedge^{r+1-k}\sE \to \bigwedge^{r+1-k} \sF,$$
where $r$ is the rank of $\sE$. By convention, if $k > r$, we set 
$V_k(f)=\emptyset$.
\end{defn}

Thus, this gives a scheme structure to the locus on which $f$ has kernel
of dimension at least $k$.  While statements in terms of such loci are often
obvious on a set-theoretic level, we are concerned with canonical scheme
structures, and scheme-theoretic statements require more thought. Our 
first observation is that this scheme structure depends only on the 
``universal kernel'' of $f$, in the following sense:

\begin{lem}\label{lem:detl-invariance} Suppose that $f:\sE \to \sF$ and
$f':\sE' \to \sF'$ are two morphisms of locally free sheaves of finite
rank on a locally Noetherian scheme $B$, and suppose that for every closed 
subscheme $Z \subseteq B$, we have 
$$\ker (f|_Z)=\ker (f|_{Z'}).$$
Then for each $k \geq 0$, we have $V_k(f)=V_k(f')$.
\end{lem}

Most of the argument which follows is due to David Eisenbud.

\begin{proof} The statement being local on $B$, we may assume $B=\Spec A$
is affine, with $A$ local and $V_k(f)$ and $V_k(f')$ cut out by ideals $I$ 
and $I'$, respectively. Suppose that $I \neq I'$. Then we may assume without 
loss of generality that $I \not\subseteq I'$, and then replacing $A$ by 
$A/I'$, we have that all $(r'+1-k)\times(r'+1-k)$ minors of $f'$ vanish, but
not all $(r+1-k)\times(r+1-k)$ minors of $f$ vanish, where $r$ and $r'$
are the ranks of $\sE$ and $\sE'$, respectively. Let $g$ be one of the
non-zero minors of $f$. Let $J$ be an ideal maximal among those not
containing $g$; we claim that if we mod out by $J$, we obtain a Gorenstein
Artin local ring. That we obtain an Artin local ring follows from
the Krull Intersection Theorem.
On the other hand, in the Artin case the Gorenstein property is equivalent 
to having simple socle (Proposition 21.5 of \cite{ei1}). But since the 
socle is by definition a vector space over the residue field, if it is not 
simple it contains an element linearly independent from $g$, which is not 
possible by the maximality of $J$.

We have thus reduced to the case $B=\Spec A$, with $A$ a Gorenstein Artin
local ring. We can now apply duality in this context (\S 21.1 of \cite{ei1})
to conclude that the dual sequences
$$\Hom_A(\sF(B),A) \overset{f^{\vee}}{\to} \Hom_A(\sE(B),A) \to 
\Hom_A(\ker f(B), A) \to 0$$
and
$$\Hom_A(\sF'(B),A) \overset{(f')^{\vee}}{\to} \Hom_A(\sE'(B),A) \to 
\Hom_A(\ker f'(B), A) \to 0$$
each give free presentations of modules which are
by hypothesis isomorphic. By the theory of Fitting ideals 
(Corollary-Definition 20.4 of \cite{ei1}), we conclude that
the ideals of $(r+1-k)\times(r+1-k)$ minors of $f^{\vee}$ and of
$(r'+1-k)\times(r'+1-k)$ minors of $(f')^{\vee}$ are equal. But since
the matrices in question are simply the transpose of those for $f$ and
$f'$, this gives a contradiction.
\end{proof}

Next, it will be helpful to know that one obtains morphisms between different
vanishing loci when expected. 

\begin{prop}\label{prop:detl-inclusion}
Let $f:\sE \to \sF$ and $f':\sE' \to \sF'$ be morphisms of
locally free sheaves of finite rank on $B$.
Suppose we also have morphisms
$g:\sE \to \sE'$ and $h:\sF \to \sF'$ such that $f' \circ g=h \circ f$,
and such that at every point of $B$, after restriction to fibers we have 
$g$ injective on $\ker f$. Then for any $k \geq 0$, we have that
$V_k(f)$ is a closed subscheme of $V_k(f')$.
\end{prop}

\begin{proof} 
The question being local, we may assume that $B=\Spec A$, where $A$ is a 
local ring having maximal ideal $\fm$. Using Nakayama's lemma, it is
easy to see that we may write $\sE \overset{f}{\to} \sF$ as a direct sum
a minimal and a trivial $2$-term complex (here minimal means that the
morphism is $0$ modulo $\fm$, and trivial means that the morphism is an
isomorphism). We may then replace $\sE \overset{f}{\to} \sF$ by the
minimal subcomplex without affecting our hypotheses or the vanishing loci
in question, 
and in this case we see that since $f$ is $0$ modulo $\fm$, then
$g$ is injective modulo $\fm$, and we see we have reduced to the case
that $g$ realizes $\sE$ as a subbundle of $\sE'$. 

Let $r$ and $r'$ denote the ranks of $\sE$ and $\sE'$, respectively.
Now, if we restrict to $V_k(f)$, we need only verify that $V_k(f')$ is
all of $B$, or equivalently, that the morphism
$$\bigwedge^{r'+1-k} f':\bigwedge^{r'+1-k}\sE' \to \bigwedge^{r'+1-k} \sF'$$
vanishes identically. However, one can
check this easily on bases from the hypotheses that $\bigwedge^{r+1-k} f=0$
and that $\sE'$ is a subbundle of $\sE$.
\end{proof}

\begin{rem} Note the contrast between Lemma \ref{lem:detl-invariance}
and Proposition \ref{prop:detl-inclusion} that the latter requires a
morphism of complexes, while the former does not. The morphism of
complexes is very much necessary for the validity of Proposition
\ref{prop:detl-inclusion}, as demonstrated by simple examples over a
non-reduced point. By the same token, if in the statement of the proposition
we assume that $g$ induces an isomorphism of kernels at all points, it does
not follow that the $k$th vanishing loci are the same.
\end{rem}

The final statement of this form is a bit more complicated to state, although
the proof is not difficult. 

\begin{prop}\label{prop:fiber-prod-detl} Given $f:\sE \to \sF$ and
$f_i:\sE_i \to \sF_i$ for $i=1,2$ morphisms of locally free sheaves of 
finite rank on $B$, as well as a locally free sheaf $\sK$ of rank $k$ and 
imbeddings $\sK \hookrightarrow \sF_i$ as subbundles for $i=1,2$, suppose
we have a closed subscheme $Z$ and positive integers $m_1,m_2$ such that:
\begin{Ilist}
\itm for $i=1,2$, the $m_i$th vanishing locus of 
$$f'_i:\sE_i \to \sF_i/\sK$$
contains $Z$;
\itm after restriction to any closed subscheme $Z'$ of $Z$, the kernel of
$f$ is isomorphic to the fibered product of the kernels of $f'_1$ and
$f'_2$ over $\sK$.
\end{Ilist}
Then 
$$V_{m_1+m_2-k}(f) \supseteq Z.$$
\end{prop}

\begin{proof} The hypotheses and conclusion being compatible with base change,
we may first restrict to $Z$, and thus wish to show that 
$V_{m_1+m_2-k}(f) =B$. Similarly, the question being local we may assume
that $B=\Spec A$ is affine, so that $\sK \to \sF_i$ splits, and we choose
isomorphisms $\sF_i\cong \sF_i/\sK \oplus \sK$. Making use of this 
isomorphism, we construct a morphism
$$g:\sF_1 \oplus \sF_2 \to \sF_1/\sK \oplus \sF_2/\sK \oplus \sK$$
induced by $f_1$ and $-f_2$, and we observe that the kernel of $g$ is
(universally) equal to the fibered product of the kernels of $f'_1$ and 
$f'_2$ over $\sK$. Thus, according to Lemma \ref{lem:detl-invariance},
it is enough to show that $V_{m_1+m_2-k}(g)=B$.

But this follows immediately from the definitions: we localize further so
that all modules are free, and we write $r_i$ and $s_i$ for the ranks of
$\sE_i$ and $\sF_i$, respectively, then the upper
$s_1+s_2$ rows of the matrix defining $g$ are block-diagonal, with
blocks of size $s_1 \times r_1$ and $s_2 \times r_2$. Our hypotheses
imply that the $(r_i+1-m_i) \times(r_i+1-m_i)$ minors of the $i$th
block for $i=1,2$ vanish uniformly, and the desired statement follows.
\end{proof}

\begin{rem} 
Suppose we have a finite complex of locally free sheaves of finite rank.
Flatness implies that quasi-isomorphism is preserved under base change, 
so it follows from Lemma \ref{lem:detl-invariance} that the $k$th 
vanishing locus of the first morphism in the complex is in fact a
quasi-isomorphism invariant. This is enough for most of our applications
to pushforwards, and can be proved directly, but one comparison theorem
will involve universally isomorphic kernels which are not naturally
obtained from quasi-isomorphic complexes. 

We also observe that for the quasi-isomorphism invariance, it
is not enough to have a morphism of complexes inducing an isomorphism on
the first two cohomology groups, or indeed on the first $n$ cohomology groups
for any fixed $n$. Indeed, letting $B=\AA^n_k$, consider $\cF^{\bullet}$
a trivial complex, say supported on $\cF_0$ and $\cF_1$, and let 
$\cG^{\bullet}$ be the Koszul complex, so that the cohomology of 
$\cG^{\bullet}$ is nonzero only in the $n$th and final place.
Then there is a morphism $\cF^{\bullet} \to \cG^{\bullet}$ inducing
isomorphisms on $H^0$ through $H^{n-1}$, but the $k$th vanishing loci of 
$\cF^{\bullet}$ for $k>0$ are empty, while those for $\cG^{\bullet}$ are
not.
\end{rem}

\subsection{Subbundles of pushforwards}
We now turn towards pushforwards. It turns out that it is useful to
have generalizations of determinantal loci not only for maps of pushforwards,
but also incorporating subbundles of pushforwards. We thus begin by
recalling the following definition.

\begin{defn}\label{def:subbundle} Let $\pi: X \rightarrow B$ be a 
proper morphism, locally of finite presentation,
and $\sE$ a quasicoherent sheaf on $X$, locally finitely presented and flat 
over $B$. A subsheaf 
$\sV$ is defined to
be a \textbf{subbundle} of $\pi_* \sE$ if $\sV$ is locally free of finite 
rank, and for any $S \rightarrow B$, the natural map
$\sV_S \rightarrow \pi_{S*} \sE_S$ remains injective.
\end{defn}

\begin{rem}\label{rem:perfect-complex}
If $B$ is locally Noetherian, our hypotheses simplify to
$\pi$ being proper, and $\sE$ being coherent and flat over $B$. However,
we want to allow for $B$ not locally Noetherian, since we don't want to
restrict to the category of locally Noetherian schemes in defining our
moduli stacks.

A basic fact about the situation of Definition \ref{def:subbundle} is
that on any affine open subset of $B$, there exists a finite complex
consisting of locally free sheaves of finite rank which is a representative
for $R \pi_* \sE$, and remains so after arbitrary base change. In the
Noetherian case, this is Theorem 6.10.5 (see also Remark 6.10.6) of
\cite{ega32}, and the general case follows by Noetherian approximation
and our finite presentation hypotheses, noting also that flatness
implies that formation of $R \pi_* \sE$ commutes with base change. 
\end{rem}

The above definition of subbundle has a number of desirable properties,
which we explore in the following lemma.

\begin{lem}\label{lem:grd-sub} In the situation of Definition 
\ref{def:subbundle}, we have the following statements.
\begin{ilist}
\itm Suppose that
$$\sF^{\bullet}=\sF^0 \overset{d^0}{\to} \dots \overset{d^{n-1}}{\to} \sF^n$$
is a finite complex of locally free sheaves
of finite rank on $B$ which is a representative of $R\pi_* \sE$.
Then subbundles of $\pi_* \sE$ in our sense are the same as 
locally free subsheaves of $\pi_* \sE$ which are subbundles of $\sF^0$ in the 
usual sense.
\itm Suppose we have $\sE$ such that $\pi_* \sE$ is locally free, and
$R^i \pi_* \sE=0$ for all $i>0$. Then our definition of subbundle of
$\pi_* \sE$ is equivalent to the usual one.
\itm If $\sV$ is locally free, a morphism
$\sV \to \pi_* \sE$ yields a subbundle if 
and only if for all points $y \in B$, the induced map
\begin{equation}\label{eq:subbundle-pt} \sV|_y \to H^0(X_y,\sE|_y)
\end{equation}
is injective. Moreover, if for some $y \in B$ we have \eqref{eq:subbundle-pt}
injective, then $\sV \subseteq \pi_* \sE$ is a subbundle on an open
neighborhood of $y$.
\itm Let $\sV_1$, $\sV_2$ be subbundles of rank $r$ of $\pi_* \sL$ in our
sense, and suppose $\sV_1 \subset \sV_2$. Then $\sV_1=\sV_2$.
\end{ilist}
\end{lem}

Note that (iii) says that our definition of subbundle is
the same as the \textit{family of $\fg^r_d$'s} used in \S IV.3 of 
\cite{a-c-g-h} to describe the functor represented by the classical $G^r_d$
space for a smooth curve.

\begin{proof} For (i), by hypothesis, we have an injection
$\pi_* \sE \to \sF^0$, so the statement makes sense. Let 
$\sV \subseteq \pi_* \sE$ be a locally free 
subsheaf, and define $\sQ$ by the exact sequence 
$$0 \to \sV \to \sF^0 \to \sQ \to 0$$
obtained by composing with the above injection.
We then wish to show that $\sV$ is a subbundle of $\pi_* \sE$ in our sense
if and only if $\sQ$ is locally free. Because the $\sF^i$ are flat over $B$,
for any $S \to B$, we have that $\sF^{\bullet}|_S$ is a representative of
$R\pi_{S*} \sE|_S$; what is relevant for our purposes is that
$(\pi_* \sE)|_S \to \sF^0|_S$ factors through the canonical morphism
$(\pi_* \sE)|_S \to \pi_{S*} \sE|_S$, with the induced morphism
$\pi_{S*} \sE|_S \to \sF^0|_S$ being injective.
Thus, the map
$\sV|_S \to \sF^0|_S$ is obtained as the composition
$\sV|_S \to \pi_{S*} \sE|_S \to \sF^0|_S$, 
and the second map is injective. We conclude that 
$\sV$ is a subbundle of $\pi_* \sE$ in our sense if and only if the induced
map $\sV|_S \to \sF^0|_S$ is injective for all $S \to B$, which is equivalent
to the condition that $\sQ$ is locally free.

For (ii), we simply observe that in this case, if we set $\sF=\pi_* \sE$,
the hypotheses of (i) are satisfied because by \cite[Cor. 6.9.9]{ega32} (see
also \cite[6.2.1]{ega32}), the natural map
$(\pi_* \sE)_S \rightarrow \pi_{S*} \sE_S$ is an
isomorphism for all $S \to B$.

For (iii), since restriction to $y \in B$ is a special
case of base change, only the ``if'' direction requires argument. Suppose
injectivity is satisfied for all $y \in B$, and let $S \to B$ be any 
morphism. Because field extensions are flat, we have that for any $s \in S$,
the map $\sV|_s \to \pi_{s*} \sE|_s$ is the base extension of
$\sV|_y \to \pi_{y*} \sE|_y$, where $y$ is the image of $s$ in $B$,
and is in particular injective. By Remark \ref{rem:perfect-complex}, we
may choose $\sF^{\bullet}$ as in (i);
then $\pi_{s*} \sE|_s$ is the kernel of 
$\sF^0|_s \to \sF^1|_s$, so we conclude by Nakayama's lemma 
that $\sV|_S \to \sF^0|_S$ is injective,
and hence that $\sV|_S \to \pi_{S*} \sE|_S$ is likewise injective, as
desired.

Similarly, if we are given only that \eqref{eq:subbundle-pt} is injective
for a single $y \in B$, since injectivity of $\sV|_y \to \sF^0|_y$ is
the complement of a determinantal locus, it holds on an open neighborhood
of $y$, and we conclude that $\sV$ is a subbundle on that neighborhood.

Finally, (iv) is straightforward: let $\sQ = \sV_2/\sV_1$, and let $y \in B$ 
be any point of $B$. If we base change to $\Spec \kappa(y)$, we get from the 
definition of subbundle that $(\sV_1)_y \hookrightarrow (\sV_2)_y$, so since 
both have dimension $r$, we get $\sQ_y=0$,
and by Nakayama's lemma we conclude $\sQ=0$ and $\sV_1=\sV_2$, as asserted.
\end{proof}

\subsection{Generalized determinantal loci}
We now give the foundational definition generalizing determinantal loci
to pushforwards.

\begin{defn}\label{def:pushforward-locus} 
Let $\pi:X \to B$ be a proper morphism, locally of finite presentation,
and let $\sE,\sF$ be quasicoherent sheaves on $X$, locally of finite 
presentation and flat over
$B$. Suppose we are given also a morphism $f:\sE \to \sF$, as well as a
subbundle $\sW \subseteq \pi_* \sF$.
Then the \textbf{$k$th vanishing locus} of the induced map
$$\pi_* \sE \to \pi_* \sF/\sW$$
is defined as follows: in the case that $B$ is affine, let $\sE^{\bullet}$
and $\sF^{\bullet}$ be finite complexes of finite locally free modules on
$B$ representing $R\pi_* \sE$ and $R\pi_* \sF$ respectively, and having
an induced morphism $f^{\bullet}:\sE^{\bullet} \to \sF^{\bullet}$.
Then we obtain an induced morphism
$$\sE^0 \to \sE^1 \oplus \sF^0/\sW,$$
and we define the desired locus to be the $k$th vanishing locus of this
morphism. For general $B$, we define the 
desired locus locally on an affine open cover.
\end{defn}

Note that we obtain a generalization of standard determinantal loci, since 
we can recover the usual
notion by setting $\pi=\id$ and $\sW=0$; then we have $\sE^{\bullet}=\sE$
and $\sF^{\bullet}=\sF$. The most basic result for our generalized version
is as follows.

\begin{prop}\label{prop:pushforward-detl-defined} In the situation of
Definition \ref{def:pushforward-locus}, the $k$th vanishing locus is a 
well-defined closed subscheme of $B$, stable under base change, and having 
as its underlying set the set of points $y \in B$ on which the map 
$$H^0(X_y, \sE|_y) \to H^0(X_y,\sF|_y)/\sW|_y$$
has kernel of dimension at least $k$.

Moreover, the vanishing locus is determined as a subscheme by the 
isomorphism classes of the kernels of
$$\pi_* (\sE|_Z) \to \pi_* (\sF|_Z)/\sW|_Z$$
as $Z$ varies over closed subschemes of $B$.
\end{prop}

\begin{proof} Because the $\sE^{\bullet}$ and $\sF^{\bullet}$ are flat,
they continue to compute $R\pi_* \sE$ and $R\pi_* \sF$ after arbitrary
base change. Thus, after arbitrary base change $S \to B$, the kernel of
$$\sE^0 \to \sE^1 \oplus \sF^0/\sW$$
is equal to the kernel of 
$$\pi_{S*} \sE|_S \to \pi_{S*} \sF|_S/\sW|_S$$
and we see immediately that the set-theoretic support is as asserted. 
We further conclude from Lemma \ref{lem:detl-invariance} that our definition
is independent of choices, and is determined by kernels as stated.
Finally, since also formation of 
determinantal ideals commutes with pullback, the construction is
stable under base change.
\end{proof}

\begin{ex} Considering the special case that $\sF=0$, 
we find that the $k$th vanishing locus is supported on the set
of points of $y \in B$ on which $H^0(X_y, \sE|_{X_y})$ has dimension at
least $k$.

In particular, if $Y$ is a variety over an algebraically closed field $F$,
and $B$ is $\Pic(Y/F)$, and we set $X=Y \times_F B$, and let $\sE$ be the
universal line bundle, then we recover a new way of looking at classical
Brill-Noether loci, which generalizes to higher-dimensional varieties.
\end{ex}

We next translate some standard statements relating kernels to 
determinantal loci into our generalized context.

\begin{prop}\label{prop:pushforward-detl-subbundles}
In the situation of Definition \ref{def:pushforward-locus}, 
suppose there exists a subbundle $\sV$ of $\pi_* \sE$ of rank $k$ which is 
contained in the kernel of $\pi_* \sE \to \pi_* \sF/\sW$. Then the $k$th 
vanishing locus is equal (scheme-theoretically) to $B$. 

Conversely, if the 
$k$th vanishing locus is equal (scheme-theoretically) to $B$, and if the 
$(k+1)$st vanishing locus is empty, then the kernel $\sK$ of 
$\pi_* \sE \to \pi_* \sF/\sW$ is a rank-$k$ subbundle of $\pi_* \sE$, and
we have moreover that the kernel commutes with base change: that is, for 
all $S \to B$, the kernel of  
$$\pi_{S*} \sE|_S \to \pi_{S*} \sF|_S/\sW|_S$$
is equal to $\sK|_S$.
\end{prop}

\begin{proof} Suppose that $\sE^0$ has rank $r$.
For the first assertion, we have by Lemma \ref{lem:grd-sub}
that $\sV$ gives a subbundle of $\sE^0$ in the usual sense.
By definition, the vanishing locus in question is equal to the 
$k$th vanishing locus of
\begin{equation}\label{eq:pushforward-detl} \sE^0 \to \sE^1 \oplus \sF^0/\sW,
\end{equation}
and by construction $\sV$ is contained in the kernel of this map. The
statement thus reduces to the standard fact that if a morphism of locally
free sheaves contains a subbundle of rank $k$ in its kernel, then the
$k$th vanishing locus is the entire scheme.

Conversely, if we have that \eqref{eq:pushforward-detl} has fixed rank
$r-k$ on $B$ (i.e., the $(r-k)$th determinantal locus is all of $B$, 
while the $(r-k+1)$st is empty) then the kernel is a subbundle of rank $k$,
which is universal under base change. Indeed, by Proposition 20.8 (see also 
comments on p.\ 407) of \cite{ei1},
we have that the image is a subbundle
of rank $r-k$, and hence the kernel is a subbundle of rank $k$, and commutes
with base change.
Moreover, since the kernel of \eqref{eq:pushforward-detl} is
equal to the kernel of $\pi_* \sE \to \pi_* \sF/\sW$ even after arbitrary
base change, we obtain the desired statement.
\end{proof}

We conclude with two corollaries translating our results on determinantal
loci to the generalized context.

\begin{cor}\label{cor:pushforward-morphisms} Suppose we are in the
situation of Definition \ref{def:pushforward-locus}, and we are given
also $\sE'$, $\sF'$, $f'$ and $\sW'$ as in the same definition, as well as 
morphisms $g:\sE \to \sE'$ and $h:\sF \to \sF'$ satisfying:
\begin{Ilist}
\itm $h \circ f = f' \circ g$;
\itm $h(\sW) \subseteq \sW'$;
\itm at every point $y \in B$, we 
have that $g$ induces an injection 
$$\ker\left(\pi_* (\sE|_{X_y}) \to \pi_* (\sF|_{X_y}) / \sW|_y\right)
\hookrightarrow 
\ker\left(\pi_* (\sE'|_{X_y}) \to \pi_* (\sF'|_{X_y}) / \sW'|_y\right).$$
\end{Ilist}
Then for every $k$, the $k$th vanishing locus of
$\pi_* \sE \to \pi_* \sF/\sW$ is a closed subscheme of the $k$th vanishing
locus of $\pi_* \sE' \to \pi_* \sF'/\sW'$.
\end{cor}

Note that the first two conditions ensure that the kernel of
$\pi_* \sE \to \pi_* \sF/\sW$ maps into the kernel of 
$\pi_* \sE' \to \pi_* \sF'/\sW'$, so the third is asserting simply that
this map is universally injective. Under these conditions, it is clear that
we have a set-theoretic inclusion of the loci in question.

\begin{proof} We may obviously assume that $B$ is affine, so that, 
continuing with the notation of Definition \ref{def:pushforward-locus}, 
the vanishing loci in question are defined by determinantal loci of 
morphisms $\sE^0 \to \sE^1 \oplus (\sF^0/\sW)$ and
$(\sE')^0 \to (\sE')^1 \oplus ((\sF')^0/\sW')$, and conditions (I) and (II)
imply that $g$ and $h$ induce
a morphism from the first complex to the second.
As in the proof of Proposition \ref{prop:pushforward-detl-defined},
we have that 
formation of our complexes commutes with restriction to any point of $B$, 
so (III) implies that the hypotheses of
Proposition \ref{prop:detl-inclusion} are satisfied, and we conclude the
desired statement.
\end{proof}

\begin{cor}\label{cor:pushforward-fibered-prod} Suppose we are in the
situation of Definition \ref{def:pushforward-locus}, and we are also
given $\sE_i$, $\sF_i$, $f_i$, and $\sW_i$ as in the same definition
for $i=1,2$. Suppose further that we have a locally free sheaf $\sW'$ of
rank $k$ on $S$, maps $\sW' \to \sW_i$ for $i=1,2$ realizing $\sW'$ as
a subbundle of each, integers $m_1,m_2$ and a closed subscheme $Z$ of 
$B$ such that:
\begin{Ilist}
\itm for $i=1,2$, the $m_i$th vanishing locus of 
$\pi_* \sE_i \to \pi_* \sF_i/\sW_i$ contains $Z$;
\itm for all closed subschemes $Z'$ of $Z$, we have that the kernel of
$\pi_{Z'*} \sE|_{Z'} \to \pi_{Z'*} \sF|_{Z'}/\sW|_{Z'}$ is isomorphic to
the fibered product of the kernels of
$\pi_{Z'*} \sE_i|_{Z'} \to \pi_{Z'*} \sF_i|_{Z'}/\sW_i|_{Z'}$ over
$\sW'|_{Z'}$.
\end{Ilist}
Then the $(m_1+m_2-k)$th vanishing locus of 
$\pi_* \sE \to \pi_* \sF/\sW$ contains $Z$.
\end{cor}

\begin{proof} Once again, the fact that the complexes in question are
compatible with base change immediately reduces the statement of the
corollary to that of Proposition \ref{prop:fiber-prod-detl}.
\end{proof}

\begin{rem} If, in the situation of Definition \ref{def:pushforward-locus},
we are given also a subbundle $\sV \subseteq \pi_* \sE$, one can make an 
analogous definition of the $k$th vanishing locus of the induced map
$$\sV \to \pi_* \sF/\sW.$$
However, in full generality we do not know 
that $\pi_* \sE$ is a subbundle of itself, so this version is not a
generalization of the stated version. We have stated the definition in the
form in which we apply it.
\end{rem}

\newpage
\section*{Index of notation and terminology}

\begin{minipage}[t]{.55\textwidth}
\noindent$v_{\I}(\e)$, Notation \ref{not:v-e-P}
\vspace{2pt}

\noindent$e_{\I}(\e)$, Notation \ref{not:v-e-P}
\vspace{2pt}

\noindent$P(w,(e_1,v_1),\dots,(e_m,v_m))$, Notation \ref{not:v-e-P}
\vspace{2pt}

\noindent$v_{\II}(\e)$, Notation \ref{not:v-e-P}
\vspace{2pt}

\noindent$P(w,v_1,\dots,v_m)$, Notation \ref{not:v-e-P}
\vspace{2pt}

\noindent$w_v$, Notation \ref{not:wv}
\vspace{2pt}

\noindent$t_{(e,v)}(w)$, Notation \ref{notn:twists} 
\vspace{2pt}

\noindent$G_{\I}$, Definition \ref{def:graphs} 
\vspace{2pt}

\noindent$G_{\II}$, Definition \ref{def:graphs} 
\vspace{2pt}

\noindent$\bar{G}_{\II}$, Definition \ref{def:bar-g}
\vspace{2pt}

\noindent$\sO_{(e,v)}$, Notation \ref{not:twisting-bundles} 
\vspace{2pt}

\noindent$\sO_v$, Notation \ref{not:twisting-bundles-ii} 
\vspace{2pt}

\noindent$\sO_{w,w'}$, Notation \ref{not:twisting-bundles} and \ref{not:twisting-bundles-ii} 
\vspace{2pt}

\noindent$\sO'_{w,w'}$, Notation \ref{not:rk-1-twist}
\vspace{2pt}

\noindent$\sE_w$, Notation \ref{not:twisting-bundles} and \ref{not:twisting-bundles-ii} 
\vspace{2pt}

\noindent$f_{\e}$, Notation \ref{not:fe} and \ref{not:type-ii-maps}
\vspace{2pt}

\noindent$f_{P}$, Notation \ref{notn:path-maps} 
\vspace{2pt}

\noindent$D_{w,v}$, Notation \ref{not:twisting-divisors} 
\vspace{2pt}

\noindent$\sE^v(w)$, Notation \ref{not:twisted-subspaces} 
\vspace{2pt}

\noindent$V^{v}(w)$, Notation \ref{not:twisted-subspaces} 
\end{minipage}
\begin{minipage}[t]{.45\textwidth}
\noindent$\cG^{k,\I}_{r,d,d_{\bullet}}(X/B)$, Definition \ref{def:grd-space-i}

\noindent$\cG^{k,\II}_{r,d,d_{\bullet}}(X/B,\theta_{\bullet})$, Definition \ref{def:grd-space-ii}

\noindent$\bar{\cG}^{k,\II}_{r,d,d_{\bullet}}(X/B,\theta_{\bullet})$, Definition \ref{def:bar-g}

\noindent$\cG^{k,\EHT}_{r,d,d_{\bullet}}(X)$, Notation \ref{not:eht}

\noindent$\cG^{k,\EHT,\refn}_{r,d,d_{\bullet}}(X)$, Notation \ref{not:eht}

\noindent$\cG^{k,\EHT}_{r,d,d_{\bullet}}(X)$, Definition \ref{defn:eht-stack}

\noindent$\cG^{k,\EHT}_{r,d,d_{\bullet},a^{\Gamma}}(X)$, Definition \ref{def:eht-var}

\noindent$\widetilde{\cG}^{k,\EHT}_{r,d,d_{\bullet}}(X)$, Definition \ref{def:tildeG}

\noindent$\cG^{k,\I}_{r,\sL,d_{\bullet}}(X/B)$, Definition \ref{def:fixed-det-stacks}

\noindent$\cG^{k,\II}_{r,\sL,d_{\bullet}}(X/B,\theta_{\bullet})$, Definition \ref{def:fixed-det-stacks}

\noindent$\cG^{k,\EHT}_{r,\sL,d_{\bullet}}(X/B)$, Definition \ref{def:fixed-det-stacks}

\noindent$\cM_{r,w}(X/B)$, Definition \ref{def:multidegree}

\noindent$\cM_{r,w_0,d_{\bullet}}(X/B)$, Notation \ref{not:bundles-substack} 

\noindent$\cM_{r,w_0,\sL}(X_{B'}/B')$, Notation \ref{not:bundles-fixed-det}

\noindent$\cM_{r,w_0,\sL,d_{\bullet}}(X_{B'}/B')$, Notation \ref{not:bundles-fixed-det}

\noindent$\cP^k_{r,d_{\bullet}}(X)$, Notation \ref{not:pkrd}

\noindent$\widetilde{\cP}^k_{r,d_{\bullet}}(X)$, Notation \ref{not:tildePkrd} 

\noindent$V_k(f)$, Definition \ref{def:vanishing-locus}

\noindent$\LG(r,\sE_{\bullet})$, Definition \ref{def:pre-lg}

\noindent$\cLG(r,\sE_{\bullet})$, \S \ref{sec:lg-stacks}
\end{minipage}
\vspace{10pt}

\noindent$(P,Q)$-adapted, Definition \ref{def:adapted}

\noindent$(P,Q)$-adaptable, Definition \ref{def:adapted}

\noindent almost local, Definition \ref{def:almost-local}

\noindent chain adaptable, Definition \ref{def:chain-adapted}

\noindent constrained, Definitions \ref{def:constrained} and \ref{def:fixed-det-stacks}

\noindent limit linear series, Eisenbud-Harris-Teixidor (or EHT), Definition \ref{def:eht-lls}

\noindent linked linear series, type-I, Definition \ref{def:grd-space-i}

\noindent linked linear series, type-II, Definition \ref{def:grd-space-ii}

\noindent multidegree (in families), Definition \ref{def:multidegree}

\noindent polarization, Definition \ref{def:polarization}

\noindent prelinked Grassmannian, Definition \ref{def:pre-lg}

\noindent refined, Definitions \ref{def:eht-lls} and \ref{def:fixed-det-stacks}

\noindent simple (linked linear series), Definitions \ref{def:simple} and \ref{def:fixed-det-stacks}

\noindent simple (point of linked Grassmannian), Definition \ref{def:simple-pt}

\noindent smoothing family, Definition \ref{def:smoothing-family}

\noindent semistable (vector bundle), Definition \ref{def:polarization}

\noindent stable (vector bundle), Definition \ref{def:polarization}

\noindent semistable (limit linear series), Definition \ref{def:stable-lls}

\noindent stable (limit linear series), Definition \ref{def:stable-lls}

\noindent $\ell$-semistable (vector bundle), Definition \ref{def:l-stable}

\noindent $\ell$-stable (vector bundle), Definition \ref{def:l-stable}

\noindent $\ell$-semistable (limit linear series), Definition \ref{def:stable-lls}

\noindent $\ell$-stable (limit linear series), Definition \ref{def:stable-lls}

\noindent subbundle (of a pushforward), Definition \ref{def:subbundle}

\noindent $k$th vanishing locus (of a morphism of vector bundles), Definition \ref{def:vanishing-locus}

\noindent $k$th vanishing locus (of the pushforward of a morphism), Definition \ref{def:pushforward-locus}

\bibliographystyle{amsalpha}
\bibliography{gen}

\end{document}